\documentclass[12pt,a4paper,draft]{article}
\usepackage[catalan,spanish,english]{babel}
\usepackage[psamsfonts]{amssymb}
\usepackage{cmmib57}
\usepackage{exscale}
\usepackage{amsmath}
\usepackage{amscd}
\usepackage{latexsym}
\usepackage{graphicx}

\usepackage{pdfsync}
\numberwithin{equation}{section}
\numberwithin{subsection}{section}
\newtheorem{thm}{Theorem}[section]
\newtheorem{cor}[thm]{Corollary}
\newtheorem{lem}[thm]{Lemma}
\newtheorem{prop}[thm]{Proposition}

\newtheorem{rem}{Remark}[section]
\def\1{{\rm 1 }\hskip -0.21truecm 1} 
\newcommand{\IR}{\mathbb{R}}
\newcommand{\rk}{{\mathbb{R}}^k}
\newcommand{\tf}{\mathcal{F}}
\newcommand{\T}{\mathbb{T}_0}
\newcommand{\cO}{\mathcal{O}}

\newcommand{\hacd}{{\mathcal{H}}^d}
\newcommand{\hac}{\mathcal{H}}
\newcommand{\kac}{\mathcal{K}}
\newcommand{\sT}{\sf T}
\newcommand{\ep}{\varepsilon}
\newcommand{\beq}{\begin{equation}}
\newcommand{\eeq}{\end{equation}}
\newcommand{\beqn}{\begin{equation*}}
\newcommand{\eeqn}{\end{equation*}}
\newcommand{\qed}{\quad{$\square$}}

\begin{document}
\begin{titlepage}
\null 
\begin{center}
{\Large\bf Hitting probabilities
for non-linear systems\\
 of
stochastic waves}
\medskip

by\\
\vspace{7mm}

\begin{tabular}{l@{\hspace{10mm}}l@{\hspace{10mm}}l}
{\sc Robert C. Dalang}$\,^{(\ast)}$  &and&{\sc Marta Sanz-Sol\'e}$\,^{(\ast\ast)}$\\
{\small Institut de Math\'ematiques}          &&{\small Facultat de Matem\`atiques}\\
{\small Ecole Polytechnique F\'ed\'erale}  &&{\small Universitat de Barcelona}\\
{\small Station 8 }                                      &&{\small Gran Via 585}\\
{\small 1015 Lausanne, Switzerland}        &&{\small 08007 Barcelona, Spain}\\              
{\small e-mail: robert.dalang@epfl.ch}       &&{\small e-mail: marta.sanz@ub.edu}\\
\null

\end{tabular}
\end{center}


\noindent{\bf Abstract.} We consider a $d$-dimensional random field $u = \{u(t,x)\}$ that solves a non-linear system of stochastic wave equations in spatial dimensions $k \in \{1,2,3\}$, driven by a spatially homogeneous Gaussian noise that is white in time. We mainly consider the case where the spatial covariance is given by a Riesz kernel with exponent $\beta$. Using Malliavin calculus, we establish upper and lower bounds on the probabilities that the random field visits a deterministic subset of $\IR^d$, in terms, respectively, of Hausdorff measure and Newtonian capacity of this set. The dimension that appears in the Hausdorff measure is close to optimal, and shows that when $d(2-\beta) > 2(k+1)$, points are polar for $u$. Conversely, in low dimensions $d$, points are not polar. There is however an interval in which the question of polarity of points remains open.

\bigskip

\noindent{\bf Keywords:} Hitting probabilities, systems of stochastic wave equations, spatially homogeneous Gaussian noise, capacity, Hausdorff measure. 

\medskip
\noindent{\sl 2010 Mathematics Subject Classifications.} Primary: 60H15, 60J45; Secondary: 60G15, 60H07, 60G60.
\medskip



\noindent

{\footnotesize

{\begin{itemize}
\item[$^{(\ast)}$] Supported in part by the Swiss National Foundation for Scientific Research.
\item[$^{(\ast\ast)}$] Supported by the grant MTM-FEDER
2009-07203 from the \textit{Direcci\'on General de
Investigaci\'on, Ministerio de Ciencia e Innovaci\'on}, Spain, and by BE 2-00333, Generalitat de Catalunya.
\end{itemize}}
}
\end{titlepage}

\newpage

\section{Introduction}
\label{s0}
In this paper, we consider an $\IR^d$-valued random field 
$$
   \{u(t,x) = (u_1(t,x),\dots,u_d(t,x)),\ (t,x) \in [0,T] \times \rk\}
$$
that is a solution of the $d$--dimensional system of stochastic wave equations 
\begin{align}
\label{1.1}
&\left(\frac{\partial^2}{\partial t^2} - \Delta \right) u_i(t,x) = \sum_{j=1}^d\sigma_{i,j}\big(u(t,x)\big) \dot M^j(t,x) + b_i\big(u(t,x)\big),\nonumber\\ 
&u_i(0,x) = 0, \quad  \frac{\partial}{\partial t}u_i(0,x) = 0,  
\end{align} 
$i=1,\ldots,d$, where $T \in \IR_+^*$, $(t,x)\in\, ]0,T]\times\rk$, $k\in\{1,2,3\}$, and $\Delta$ denotes the  Laplacian on $\IR^k$.  We shall give several scientific motivations for studying such systems of equations at the end of this section.

   The free terms of this system are Lipschitz continuous functions $\sigma_{i,j}, b_i: \IR^d\rightarrow \IR$, $1\le i,j\le d$. The notation $\dot M(t,x)$ refers to the formal derivative of a $d$-dimensional Gaussian random field with independent components,  white in the time variable and with a spatially homogeneous covariance in space. More explicitely, on a complete probability space $(\Omega, \mathcal{G}, P)$ and for any $j\in\{1,\ldots,d\}$, we consider a family of centered Gaussian random variables
$M^j=\{M^j(\varphi),\, \varphi\in\mathcal{C}_0^\infty(\IR^{d+1})\}$, where $\mathcal{C}_0^\infty(\IR^{k+1})$ denotes the space of infinitely differentiable functions with compact support. The $M^j$ are independent and the covariance function of each $M^j$ is given by
\begin{equation}
\label{1.2}
 E(M^j(\varphi)M^j(\psi))=\int_{\IR_+}ds\int_{\IR^k}  \Gamma(dx)\left(\varphi(s,\cdot)\star\tilde\psi(s,\cdot)\right)(x),
\end{equation}
where ``$\star$" denotes the convolution operator in the spatial argument and $\tilde \psi(t,x)=\psi(t,-x)$. 

 As in \cite{d99}, formula \eqref{1.2} can be also written as
\begin{equation}
\label{1.3}
 E(M^j(\varphi)M^j(\psi))= \int_0^\infty \int_{\IR^k}\mu(d\xi)\, \tf\varphi(t)(\xi)\, \overline{\tf\psi(t)(\xi)},
 \end{equation}
 where $\mu$ denotes the spectral measure ($\mu=\tf^{-1}f$) and $\tf$ denotes the Fourier transform operator on $\IR^k$.
 
 The system \eqref{1.1} is rigorously interpreted as follows:
\begin{align}
\label{1.4}
u_i(t,x)&= \sum_{j=1}^d \int_0^t \int_{\IR^k} G(t-s,x-y) \sigma_{i,j}(u(s,y))\,M^j(ds,dy)\nonumber\\
   &\qquad  +\int_0^t \int_{\IR^k}G(t-s,x-y) b_i(u(s,y))\, ds dy,
\end{align}
$i=1,\ldots,d$, where $G$ denotes the fundamental solution of the wave equation, $M^j$ is extended to a worthy martingale measure and the stochastic integral is as in \cite{d99}. For the expressions of $G$ for $k=1,2,3$, and their Fourier transforms, see \eqref{fdt} and \eqref{fg}, respectively.

Assume that the spectral measure $\mu$ satisfies 
\begin{equation}\label{e1}
\int_{\IR^k} \frac{\mu(d\xi)}{1+|\xi|^2} <\infty,
\end{equation}
where $|\cdot |$ is the Euclidean norm.
Then, from an easy extension of Theorem 13 in \cite{d99}, there is a unique random field solution to \eqref{1.4} (and the above condition is also essentially necessary). This solution is a $d$--dimensional stochastic process $u=\{u(t,x)=(u_i(t,x))_{1\le i\le d},\, (t,x)\in [0,T]\times \IR^k\}$, adapted to the natural filtration associated with the martingale measure
 $\{M_t(B), t\in[0,T], B\in\mathcal{B}_b(\IR^k)\}$ and continuous in $L^2(\Omega)$. Here, $\mathcal{B}_b(\IR^k)$ stands for the set of bounded Borel sets of $\IR^k$.

\smallskip

 Our objective is to establish upper and lower bounds on the probabilities that  the process $u$ hits a  set $A\in\mathcal{B}_b(\IR^d)$ and thus, to develop a probabilistic potential theory for the process $u$. This will be obtained in terms of the Hausdorff measure and the Bessel-Riesz capacity of $A$, respectively. 

In \cite{dss10}, hitting probabilities for systems of Gaussian waves for any spatial dimension $k$ have been studied. This reference covers the case where, in \eqref{1.4}, $\sigma=(\sigma_{i,j})_{1\le i,j\le d}$ is an invertible matrix with constant entries, $b=(b^i,\ i=1,\ldots,d)\equiv 0$ and $k\in\mathbb{N}$. The following result was proved. Assume that the covariance $\Gamma$ is absolutely continuous with respect to Lebesgue measure and its density is given by $f(x)=|x|^\beta$,
with $\beta\in\,]0, 2\wedge k[$. Fix $t_0\in\,]0,T]$ and let $I$, $J$ be compact subsets of $[t_0,T]$ and $\rk$, each with positive Lebesgue measure. Fix $N>0$. There exist positive constants $c$ and $C$ depending on $I$, $J$, $N$, $\beta$, $k$ and $d$, such that, for any Borel set $A\subset[-N,N]^d$,
\begin{equation}
\label{1.9}
c\ {\rm{Cap}}_{d-\frac{2(k+1)}{2-\beta}}(A)\le P\{u(I\times J)\cap A\ne\emptyset\}\le C\ \mathcal{H}_{d-\frac{2(k+1)}{2-\beta}}(A).
\end{equation}
Here, the lower bound contains the $d-\frac{2(k+1)}{2-\beta}$-dimensional Bessel-Riesz capacity of $A$, and the upper bound the Hausdorff measure of $A$ of dimension $d-\frac{2(k+1)}{2-\beta}$ (at the end of the section, we will recall the definitions of these notions). The parameter $d-\frac{2(k+1)}{2-\beta}$ gives the optimal value that can be expected in a possible generalization of this result to  Equation \eqref{1.4}. 

The proof of \eqref{1.9} is carried out by applying some general criteria for hitting probabilities established in \cite{dss10} (see also \cite{blc09}).
They are described by properties of densities and joint densities of random variables in the random field $u$, and also by the regularity of its
sample paths. For the random field given by 
\begin{equation*}
u_i(t,x)= \sum_{j=1}^d \int_0^t \int_{\IR^k} G(t-s,x-y) \sigma_{i,j}\,M^j(ds,dy),\qquad i=1,\ldots,d,
\end{equation*}
checking the necessary assumptions to apply the above-mentioned  criteria is done by working directly on the explicit formulas for the (Gaussian) densities. 

Hitting probabilities for systems of non-linear stochastic partial differential equations have been studied in \cite{dn04}, \cite{dkn07}, \cite{dkn09}, \cite{dkn10} and \cite{dss10} (see also \cite{mt} for a (Gaussian) random string, \cite{dmz06} for a heat equation with reflection, and \cite{nv09} for a system of heat equations driven by an additive fractional Brownian motion). In the first reference, a system of two-parameter It\^o equations driven by a Brownian sheet has been considered. By a rotation of forty-five degrees, the system is transformed into a new one consisting of wave equations in spatial dimension $k=1$ driven by a space-time white noise. References \cite{dkn07,dkn09} are concerned with systems of heat equations in spatial dimension $k=1$ driven by space-time white noise. Systems of non-linear heat equations with $k\ge 1$ driven by a Gaussian noise of the type described before, with covariance measure $\Gamma(dx)=|x|^{-\beta}dx$, $\beta\in\,]0,2\wedge k[$ are considered in \cite{dkn10} (and uses some ideas developed here). In contrast with \cite{dss10} and \cite{blc09}, in non Gaussian cases, the existence, expression and properties of densities are obtained using Malliavin Calculus. Here, we also use this tool. 
For this, we shall consider the following conditions (which are standard for Malliavin calculus) on the coefficients of the system \eqref{1.4}:
\vskip 12pt

\noindent (P1) The functions $\sigma_{i,j}$ and $b_i$, $1\le i,j\le d$ are infinitely differentiable, with bounded partial derivatives of any positive order. 
Moreover, the $\sigma_{i,j}$, $1\le i,j\le d$, are bounded.
\vskip 12pt

\noindent (P2) The matrix-valued function $\sigma$ is uniformly elliptic. This implies that for any $v\in\IR^d$ with $\Vert v\Vert = 1$,
\beq
\label{1.5}
\inf_{x\in \IR^d}\Vert v^t \sigma(x)\Vert\ge \rho_0>0.
\eeq

\smallskip


Section \ref{s1} is devoted to establishing the upper bound in terms of Hausdorff measure. According to \cite[Theorem 2.4]{dss10}, it suffices to prove the existence of a density $p_{t,x}(z)$ for any random variable $u(t,x)$, $(t,x)\in]0,T]\times \IR^k$, and that this density in uniformly bounded (in its three variables) over compact sets. In addition, we must bound from above the $L^q$-moments of increments $u(t,x)-u(s,y)$ in terms of a power of the distance $|t-s|+|x-y|$, for any $q\in[1,\infty[$ (see the statements (R1) and (R2) in Section \ref{s1}). 

The existence of smooth densities in the $z$-variable for the random variables $u(t,x)$,  $(t,x)\in\,]0,T]\times \IR^k$, solution to a stochastic wave equation with $k=2$ and $k=3$ (\eqref{1.4} with $d=1$) has been proved in \cite{mss99} and \cite{qss05}, respectively. However, to the best of our knowledge, no written proof exists for $k=1$ and spatially homogeneous noise, and the issue of uniform boundedness has not been addressed either. The main ingredients for proving this result are uniform bounds in $(t,x)$ over compact sets for the $L^p$-moments
of the determinant of the inverse of the Malliavin matrix $\gamma_{u(t,x)}$. This is proved in Proposition \ref{p2} with an approach that is independent of the dimension $k\in \{1,2,3\}$. 

As for the analysis of $L^q$-moments of increments $u(t,x)-u(s,y)$, we consider different types of assumptions on the covariance measure $\Gamma$. 
The two cases---additive and multiplicative noise---are handled separately. In comparison with the latter, the former admits more general covariances and the estimates are sharper. For multiplicative noise, we rely on known results for $k=2$ \cite{mss99,mimo} and $k=3$ \cite{dss09}, but, not being aware of any reference, we give a complete proof for $k=1$. Finally, our result on upper bounds is stated in Theorem \ref{t1}.
\smallskip

The main part of this paper is devoted to establishing lower bounds on the hitting probabilities in terms of the Bessel-Riesz capacity of the set $A$. 
For this, we assume that $\Gamma(dx)=|x|^{-\beta}$, with $\beta\in\,]0,2\wedge k[$. This is hypothesis (C1) in Section \ref{s1}, which is more restrictive than the other conditions (C1'), (C2), (C3) listed there. 

Lower bounds on hitting probabilities are obtained by using a modification of Theorem 2.1 in \cite{dss10}, as follows. When applied to the random field solution to the system \eqref{1.4}, the two hypotheses of that theorem are
\begin{description}
\item{(a)} a {\it Gaussian-type} upper bound for the joint density of $(u(s,y), u(t,x))$, for$(s,y)\ne (t,x)$, $(s,y), (t,x)\in\,]0,T]\times \IR^k$;
\item{(b)} strict positivity of the density $p_{t,x}$ on $\IR^d$.
\end{description}
Property (b) has been recently established  in \cite{en10}. As for the upper bound in (a), we show here that this hypothesis can be significantly weakened. Indeed, the
exponential factor in the Gaussian-type bound can be replaced by a monomial (see the proof of Theorem \ref{rdlem1} and also Remark 3.1(b)), and this still leads to the same lower bounds on hitting probabilities as the Gaussian-type estimate. 

As in \cite{dkn09}, the main technical effort of the paper is to establish such an upper bound on the joint density. We follow the same general approach as in this reference, though the difficulties that we encounter and the way they are solved present significant differences. 
Let $Z=(u(s,y), u(t,x)-u(s,y))$ and denote by $\gamma_Z$ the Malliavin matrix of $Z$. The hard core of the analysis consists of obtaining the rate of degeneracy in terms of negative powers of $|t-s|+|x-y|$ of $L^p$-moments of the inverses of the eigenvalues of $\gamma_Z$. As was noticed in \cite{dkn09}, this is necessary only for the ``small" eigenvalues. For the ``large" eigenvalues, it suffices to prove that their inverses have $L^p$-bounded moments.



In Section \ref{s2}, we highlight this fact. Indeed, Theorem \ref{t2.10} provides lower bounds for the hitting probabilities, assuming that the conditions mentioned before on the large and small eigenvalues (see \eqref{2.111m}--\eqref{LEV}) are satisfied. 
Although the proof of Theorem \ref{t2.10} uses arguments from \cite{dkn09}, we would like to stress two novel points. First, the remark on condition (a) mentioned above. Indeed, in the examples previously analyzed in the literature, the exponential factor is obtained by applying the exponential martingale inequality. This is possible by eliminating the drift term via Girsanov's theorem. In our framework, this is not possible, since no suitable version of that theorem exists for $k>1$. Fortunately, the monomial factor can be obtained given that we have $L^p$-estimates of increments of $u(t,x)-u(s,y)$ in terms of powers of $|t-s|+|x-y|$. Indeed, this question appears naturally in relation with the analysis of sample path regularity through Kolmogorov's continuity criterion. The second point is that this theorem explains the relationship between the rate of degeneracy of the small eigenvalues that we alluded to before, described by some positive parameter $\rho$ (see \eqref{2.111m}) and the dimension of the Bessel-Riesz capacity in the lower bound on hitting probabilities (see \eqref{2.111m}). In particular, we observe from Theorem \ref{t2.10} that obtaining an optimal result for the lower bound requires an optimal result on the rate of degeneracy associated with the small eigenvalues!
\smallskip

Section \ref{s3} is devoted to the analysis of the eigenvalues of the Malliavin matrix $\gamma_Z$, in order to check the validity of assumptions \eqref{2.111m}--\eqref{LEV}, for some rate parameter $\rho$. The main tool for carrying out this program is \cite[Proposition 3.5]{dkn09}. However, we will obtain in Theorem \ref{p2.2m} a value of $\rho$ arbitrarily close to $3 - \beta$, whereas the optimal value should be $2-\beta$.


Theorem \ref{t5n} is devoted to establishing the bound \eqref{LEV} on the large eigenvalues. The proof of an analogous result in \cite{dkn09} relies on the semigroup property of the heat kernel and its properties, and uses repeatedly \cite[Proposition 3.5]{dkn09}. We also use this tool, but the overall structure of the proof is quite different. Given the absence of semigroup property, Lemma \ref{la10} plays a crucial role. This Lemma tells us that, for $\epsilon>0$, $x,y\in\IR^k$, $x\ne y$,
\begin{equation*}
\lim_{\frac{|x-y|}{\ep}\to\infty}Ê\int_0^\ep dr\, \langle G(r,y-\ast), G(h+r,x-\ast)\rangle_{\hac}=0,
\end{equation*} 
uniformly over $\frac{|t-s|}{\ep}\in[0,1]$. This result has 
some similitude with the classical Riemann-Lebesgue theorem.  In \cite{dkn09}, the space $\hac$ (related to the covariance of the spatially homogeneous noise) is replaced by $L^2([0,1])$. Finally, our main result on lower bounds is stated in Theorem \ref{t2}.


Section \ref{a} gathers a collection of auxiliary estimates that are used in the proofs of the previous section. These are of different types.
Lemmas \ref{l6.1-bis}, \ref{l6.2-bis}, \ref{la6-2-a1}, \ref{la6-2-a2}, \ref{la6-2-a3} concern integrated increments of the covariance function. 
They extend results proved in \cite[Chapter 6]{dss09} by expliciting the dependence on the length of the domain of integration in the time variable. Lemmas \ref{la1} and \ref{la7} establish upper and lower bounds on integrals in time of $\hac$-norms of the fundamental solution $G$ to the wave equation. Propositions \ref{la6-1} and \ref{la6-2} provide upper bounds for $L^p$-moments of integrals in time, on a domain of size $\ep$, of $\hac$-norms of stochastic processes which are a product of increments of the fundamental solution $G$ and an $L^p$-bounded process. Lemma \ref{la10} has already been mentioned. 

We end the description of the different sections of this paper with a short discussion about the discrepancy between the results obtained here in comparison with \eqref{1.9}
(the Gaussian case). The upper bound \eqref{u} is {\it almost} optimal. In fact, by taking $\Gamma(dx) = |x|^{-\beta}$, $\beta\in\,]0,2\wedge k[$, the dimension of the Hausdorff measure in the right-hand side of \eqref{u} is strictly less than (but arbitrarily close to)  $d-\frac{2(k+1)}{2-\beta}$. The same phenomenon appears in \cite{dkn09} in comparison with \cite{dkn07}. Concerning the lower bound, according to \eqref{t2.111}, if  the parameter $\rho$ in \eqref{2.111m} were equal to $2-\beta+\delta$, for some $\delta$ arbitrarily close to zero, then the dimension of the capacity in \eqref{t2.111} would be $d-\frac{2(k+1)}{2-\beta}+\eta$, with 
$\eta$ arbitrarily close to zero, and therefore arbitrarily close to the Gaussian case. We see in \eqref{2.11m} that we can take $\rho=3-\beta+\delta$, for any arbitrarily small $\delta$. This gives in Theorem \ref{t2} a capacity dimension $d\left(1+\frac{2d}{2-\beta}\right)+\delta-\frac{2(k+1)}{2-\beta}$, which is rather far from optimal. The reason for this discrepancy stems from the type of localization in time that we use to keep control of the size of the eigenvalues (see for instance \eqref{ss}, \eqref{s}). It is an open problem to determine a better localization procedure that would provide the optimal dimension in the lower bound. 
\smallskip

One objective in the study of hitting probabilities is to determine under which relationship between $d$, $k$ and $\beta$ points are polar or not. Recall that a point $z\in \IR^d$ is {\em polar} for $u$ if $P\{\exists (t,x) \in \,]0,T] \times \rk: u(t,x) = z\} = 0$, and it is {\em non-polar} otherwise. In general, points are polar for $d$ large and/or $k$ small, and are non-polar otherwise. The following is an immediate consequence of Theorems \ref{t1} and \ref{t2}. 

\begin{cor} Assume (P1) and (P2) above, and (C1) in Section \ref{s1}. Fix $\beta \in \,]0,2\wedge k[$ and let $u$ be the solution of \eqref{1.4}.

   (a) Suppose $d > \frac{2(k+1)}{2-\beta}$. Then points are polar for $u$.
   
   (b) Suppose $d (1+\frac{4d}{2-\beta}) < \frac{2(k+1)}{2-\beta}$. Then points are non-polar for $u$.
\end{cor}

\smallskip

We finish this introduction by defining certain notions and notations that will be used in the sequel.

Let $\gamma\in \IR$. The $\gamma$-dimensional {\it Hausdorff measure} of a Borel set $A\subset \IR^d$ is defined by 
$\mathcal{H}_\gamma(A)=\infty$ if $\gamma<0$, and for 
$\gamma\ge 0$,
\begin{equation*}
\mathcal{H}_\gamma(A)= \liminf_{\varepsilon\to 0^+}\left\{\sum_{i=1}^\infty(2r_i)^\gamma: A\subset\cup_{i=1}^\infty B_{r_i}(x_i),\ \sup_{i\ge 1}r_i\le \varepsilon\right\}.
\end{equation*}
The {\em Bessel-Riesz kernels} are defined by
\begin{equation*}
K_\gamma(r)=\begin{cases}
r^{-\gamma} & \text{if}\ \gamma>0,\\
\log\left(\frac{c}{r}\right) & \text{if}\  \gamma=0,\\
1 & \text{if}\  \gamma<0,\end{cases}
\end{equation*}
where $c$ is a positive constant.
 For every Borel set $A\subset \IR^d$, let $\mathcal{P}(A)$  be the set of probability measures on $A$. For $\mu\in\mathcal{P}(A)$, we set
\begin{equation*}
\mathcal{E}_\gamma(\mu)=\int_{A} \int_{A}K_{\gamma}(\Vert x-y\Vert)\, \mu(dx)\mu(dy).
\end{equation*}
The {\it Bessel-Riesz capacity} of a Borel set $A\subset \IR^d$ is defined as follows:
\begin{equation*}
{\rm Cap}_\gamma(A)=\left[\inf_{\mu\in\mathcal{P}(A)}\mathcal{E}_\gamma(\mu)\right]^{-1},
\end{equation*}
with the convention that $1/\infty=0$.
\smallskip

There is a Hilbert space naturally associated with the spatial covariance of $M$ (see \eqref{1.3}). Indeed, let $\hac$ be the completion of the Schwartz space of real-valued test functions $\mathcal{S}(\IR^k)$ endowed with the semi-inner product
 \begin{equation*}
 \langle\varphi,\psi\rangle_\hac = \int_{\IR^k}\mu(d\xi)\, \tf\varphi(\xi) \overline{\tf\psi(\xi)}.
 \end{equation*} 
 We define 
${\hac}^d=\{h=(h^1,\ldots,h^d): h^\ell\in \hac, \ell=1,\ldots,d\}$, $\hac_t =L^2([0,t]; \hac)$, and ${\hac}^d_t=L^2([0,t]; {\hac}^d)$, $t\in\,]0,T]$.

Assume that $\varphi\in\hac$ is a signed measure with finite total variation.  Suppose also that $\Gamma(dx)= |x|^{-\beta}dx$, with $\beta\in\,]0,2\wedge k[$ and therefore $\mu(d\xi)= c_{k,\beta}\, |\xi|^{\beta-k} d\xi$ (see \cite{d99}). Then, by applying \cite[Theorem 5.2]{kx09} (see also \cite[Lemma 12.12, page 162]{mat} for the case of probability measures with compact support) and a polarization argument on the positive and negative parts of $\varphi$, we obtain
 \begin{align}
 \label{fundamental} 
 \Vert \varphi\Vert_{\hac}^2&= \int_{\IR^k} \vert \tf \varphi(\xi)\vert^2\, \frac{d\xi}{|\xi|^{3-\beta}}\nonumber\\
 &=C\int_{\IR^k\times \IR^k}\varphi(dx)\varphi(dy)\, \vert x-y\vert^{-\beta}.
 \end{align}
 This fundamental identity will be applied on several occasions.

We recall the expressions for the fundamental solution to the wave equation in spatial dimension $k\in\{1,2,3\}$:
\begin{align}
k=1:&\quad G(t,x)=\frac{1}{2}1_{\{|x|<t\}},\nonumber\\
k=2:&\quad G(t,x)=\frac{1}{2\pi}(t^2-|x|^2)_{+}^{-\frac{1}{2}},\nonumber\\
k=3:&\quad G(t,dx)=\frac{1}{4\pi t}\sigma_t(dx),\label{fdt}
\end{align}
where $\sigma_t(dx)$ denotes the uniform measure on the sphere centered at zero and with radius $t$, with total mass $4\pi t^2$. For any dimension 
$k\in\mathbb{N}^\ast$, the Fourier transform of $G(t)$ is given by
\begin{equation}
\label{fg}
\tf G(t)(\xi)=|\xi|^{-1}\sin (t|\xi|). 
\end{equation}
We refer the reader to \cite{folland} for these formulas.

Throughout the article, we will use the notation $\T$ for a generic interval $[t_1,t_2]$, with $0<t_1<t_2\le T< \infty$, $\cO$ for a compact subset of
$\IR^d$ and $K$ for a compact subset of $\IR^k$. 

    We end this section by providing several scientific motivations for studying systems of stochastic wave equations. Recall that the deterministic wave equation is one of the three fundamental second order PDEs \cite[Chapter 2]{Evans}. It arises, via Newton's equations of motion, as a generic model for the propagation of waves (which can be acoustic, electromagnetic, etc.) in a homogeneous elastic medium. For instance, a system of three wave equations can be used to describe the motion in three dimensional space of an elastic string. If this string is in fact a DNA molecule ``floating" in a fluid, then it will be subject to (random) molecular excitation, and a system of stochastic wave equations driven by random noise provides a good description of the string's motion (here, $d= 3$ and $k=1$); see \cite[Section 1]{dalmini} for a more detailed description of the problem, \cite{gonzMaddocks} for some biological motivation, and \cite{zab} for some mathematical results on this problem. Of course, describing the motion of a family of $N$ DNA molecules would involve a system of $3N$ equations!
	
	There are several situations in signal transmission \cite{biswas} and oceanography \cite{adlermuller,miller} in which spatially correlated noise, and, in particular, spatially homogeneous noise, are used. This type of random noise is also important in the mathematical literature \cite{holley,dawson,df}. The wave equation in the case $d=1$, $k=2$, arises if one is interested in the surface waves produced by raindrops falling on the water. The case $d=1$, $k=3$, arises \cite{dl2} if one is interested in the propagation into the depths of the water of the sound waves produced by these raindrops: listening devices located under water would want to filter out this noise, and for this, a model that describes this sound propagation is useful.
	
	A last example concerns the study of shock wave propagation through the interior of the sun. This would be the case $d=1$, $k=3$, with noise that is ``spherically homogeneous" rather than spatially homogeneous: see \cite{dalmini} for a more detailed description of this problem, and \cite{dl1} for some mathematical results in this direction.
	
	Finally, it is interesting to consider the damped wave operator, and systems of equations with the wave operator $\frac{\partial^2}{\partial t^2} - \Delta$ replaced by 
$$
	a\, \frac{\partial^2}{\partial t^2} + b\, \frac{\partial}{\partial t}  - c\, \Delta\, .
$$
Here, the parameters $a$, $b$ and $c$ will have a physical interpretation, such as the mass per unit length, a friction coefficient, or an elasticity coefficient. In various limiting cases, such as $a \downarrow 0$, or $b \downarrow 0$, which correspond to the zero mass limit and the frictionless limit, respectively, one obtains either the heat operator or the wave operator (see \cite{cerraif} and the references therein for convergence of the solutions of the corresponding equations, and several motivations for this study). Thus, the study of systems of wave equations such as \eqref{1.1} can also be viewed as an important step towards the study of more general classes of systems of SPDEs.



\section{Upper bounds on hitting probabilities}
\label{s1}

We will assume that the covariance measure $\Gamma$ of the Gaussian process governing the system \eqref{1.4} is absolutely continuous with respect to Lebesgue measure on $\mathbb{R}^k$,
and we will consider the following set of hypotheses on its density $f$:
\begin{description}
\item{(C1)}  The covariance measure is of the form $\Gamma(dx) = f(x) dx$, where
$f(x)=|x|^{-\beta}$ if $x\ne 0$, with $\beta\in\,]0,2\wedge k[$.
\item{(C1')} The spectral measure $\mu$ of the noise satisfies
\beqn
\int_{\IR^k} \frac{\mu(d\xi)}{(1+|\xi|^2)^\alpha} < \infty,
\eeqn
for some $\alpha\in\, ]0,1[$.
\item {(C2)} The covariance measure is of the form $\Gamma(dx) = f(\vert x\vert) dx$, where $f:\mathbb{R}_+\to \mathbb{R}_+$ satisfies
$\int_{0^+} r^{1-\eta} f(r)dr <\infty$, for some $\eta\in\, ]0,1[$. 
\item{(C3)} The covariance measure is of the form $\Gamma(dx) = f(x) dx$, where
$f(x)=\varphi(x)|x|^{-\beta}$ if $x\ne 0$, with $\beta\in\,]0,2\wedge k[$ and $\varphi$ is a bounded and positive function belonging to $\mathcal{C}^1(\mathbb{R}^k)$ with $\nabla\varphi\in\mathcal{C}_b^\mu(\mathbb{R}^k)$ (the space of bounded and H\"older-continuous functions with exponent $\mu\in\, ]0,1]$).
\end{description}

\noindent Notice that (C1') is satisfied when $\Gamma(dx)=|x|^{-\beta}$, with $\beta \in\,]0,2\wedge k[$ and $\alpha\in\,]\beta/2,1[$. Also, (C2) is satisfied when $f:\mathbb{R}_+\to \mathbb{R}_+$ is defined by $f(r)= r^{-\beta}$, $\beta\in\,]0,2\wedge k[$, with $\eta<2-\beta$.
\medskip

This section is devoted to establishing the following result.
\begin{thm}
\label{t1} 
Let $\{u(t,x),\, (t,x)\in [0,T]\times \IR^k\}$ be the stochastic process given by the system (\ref{1.4}) of SPDEs with $k\in\{1,2,3\}$. Let $I$ and $K$ be compact subsets of $]0,T]$ and $\mathbb{R}^k$, respectively,
with positive Lebesgue measure.

Consider the following cases:
\begin{enumerate}
\item {\it Additive noise}: The matrix $\sigma=(\sigma_{i,j})_{1\le i,j\le d}$ has constant entries and $\det \sigma\ne 0$; the functions $b_i$, $1\le i\le d$, are bounded, infinitely differentiable with bounded partial derivatives of any order. The covariance measure satisfies (C1) or (C1'). 
\item {\it Multiplicative noise}: The coefficients $\sigma_{i,j}$ and $b_i$, $1\le i,j\le d$ satisfy the assumptions (P1) and (P2). For $k=2$, we assume that the covariance measure satisfies (C2). For $k\in\{1,3\}$, we assume (C3).
\end{enumerate}
Fix $\zeta\in\,]0,d[$. Then there exists a positive constant $c=c(I,K,\Gamma,k,d, \zeta)$ such that, for any Borel set $A\subset \IR^d$,
\beq
\label{u}
  P\{u(I\times K)\cap A\ne \emptyset\}\le c\mathcal{H}_{d-\zeta-\frac{k+1}{\delta_i}}(A), \, i=1,2,
  \eeq
  with the following values of $\delta_i$. In Case 1, under (C1), $\delta_1=\frac{2-\beta}{2}$, and under (C1'), $\delta_1 = 1 - \alpha$. In Case 2, for $k=1$, $\delta_2 = \frac{2-\beta}{2}$, for $k=2$, $\delta_2 = \frac{\eta}{2}$; and for $k=3$, 
  $\delta_2=\frac{2-\beta}{2}\wedge \frac{1+\mu}{2}$.
 \end{thm}
 

According to Theorem 2.4 of \cite{dss10}, in order to prove Theorem \ref{t1}, we must prove two facts.
\begin{description}
\item{(R1)} For any $(t,x)\in\T\times K$, the random vector $u(t,x)$ has a density $p_{t,x}$, and
\beq
\label{2}
\sup_{z\in{\cO}}\sup_{(t,x)\in \T\times K}p_{t,x}(z)\le C.
\eeq
\item{(R2)} There exists $\delta\in\, ]0,1]$ and a constant $C$ such that for any $q\in[1,\infty[$, $s,t \in \T$ and $x,y\in K$,
\beq
\label{3}
E(\vert u(t,x)-u(s,y)\vert^q) \le C \left(|t-s|+\vert x-y\vert\right)^{q\delta}.
\eeq
\end{description}
These two facts imply \eqref{u} with $\delta_i$ replaced by $\delta$. We note that because $\zeta$ is arbitrarily small, we will obtain \eqref{u} with $\delta_i$ provided we establish \eqref{3} for any $\delta < \delta_i$.

\subsection{Proof of (R2)}
\label{s1.1}

The results of this section can be proved under the assumption that the functions $\sigma_{i,j}$, $b_i$, $1\le i,j\le d$, are merely Lipschitz continuous, instead of supposing
the stronger assumptions (P1) and (P2).
\bigskip

\noindent{\bf Case 1: additive noise} 



\begin{prop}
\label{p1}
Let $\{u(t,x),\, (t,x)\in [0,T]\times \IR^k\}$, $k\in\{1,2,3\}$, be the solution of the system of SPDEs defined by (\ref{1.4}). Under Assumption (C1) (respectively (C1')), property (R2) holds with $\delta=\frac{2-\beta}{2}$ (respectively $\delta< 1 - \alpha$).
\end{prop}
\noindent{\em Proof.} Without loss of generality, we may assume in this proof that $d=1$ and therefore, we omit the index $i$. Define
\beqn
v(t,x) = \int_0^t \int_{\IR^k} G(t-r,x-z) M(dr, dz).
\eeqn
Fix $q\in[1,\infty[$ and consider $x,y\in K$, $t\in[0,T]$ and $h>0$ such that $t+h\in [0,T]$. Then
\beq\label{3'}
E\left(\left\vert u(t,x)-u(t+h,x+y)\right\vert^q\right)\le C(q)\left(T_1+T_2\right),
\eeq
with
\begin{align*}
T_1&= E\left(\left\vert v(t,x)-v(t+h,x+y)\right\vert^q\right),\\
T_2&=E\left(\left\vert \int_0^t dr \int_{\IR^k} dz\  G(t-r,x-z) b(u(r,z))\right.\right.\\
&\qquad -\left.\left.\int_0^{t+h} dr \int_{\IR^k} dz\ 
G(t+h-r,x+y-z)b(u(r,z))\right\vert^q\right).
\end{align*}
Since the process $\{v(t,x),\, (t,x)\in [0,T]\times \IR^k\}$ is Gaussian,
\beqn
T_1\le C(q)\left[E\left(\left\vert v(t,x)-v(t+h,x+y)\right\vert^2\right)\right]^{\frac{q}{2}}.
\eeqn
Then, by Proposition 4.1 of \cite{dss10}, under (C1), we obtain
\beq
\label{4}
T_1\le C(q) \left(h+\vert y\vert\right)^{q\delta},\qquad\mbox{with } \delta = \frac{2-\beta}{2},
\eeq
and under (C1'), it is a consequence of the calculations carried out in \cite[Corollary 7.3 and Proposition 7.4]{cd08} that
\beq
\label{4'}
T_1\le C(q) \left(h+\vert y\vert\right)^{q\delta},\qquad\mbox{with } \delta < 1-\alpha.
\eeq

For the analysis of the term $T_2$, we first write
$$
\int_0^t dr \int_{\IR^k} dz\, G(t-r,x-z) b(u(r,z))
 =\int_0^t dr \int_{\IR^k} G(r,dz)\,b(u(t-r,x-z)),
$$
and
\begin{align*}
&\int_0^{t+h} dr \int_{\IR^k} dz\, G(t+h-r,x+y-z)b(u(r,z))\\
&\qquad\qquad = \int_0^{t+h} dr \int_{\IR^k} G(r,dz)\, b(u(t+h-r,x+y-z)),
\end{align*}
and consider the terms
\begin{align*}
T_{2,1}&= E\Big(\Big\vert \int_0^t dr \int_{\IR^k} G(r,dz)\\
&\qquad\qquad\times \left[b(u(t-r,x-z))-b(u(t+h-r,x+y-z))\right]\Big\vert^q\Big),\\
T_{2,2}&= E\left(\left\vert \int_t^{t+h} dr \int_{\IR^k} G(r,dz)\, b(u(t+h-r,x+y-z))\right\vert^q\right).
\end{align*}
It has been proved in \cite{d99} that the solution to the stochastic wave equation with vanishing initial conditions is stationary in space.
Hence, by applying H\"older's inequality with respect to the finite measure $G(r,dz) dr$, using the Lipschitz property of $b$ and the change of variable $r\mapsto t-r$, we obtain
\begin{align*}
T_{2,1}&\le C \int_0^t dr \int_{\IR^k} G(r,dz)\\
&\qquad\times E\left(\left\vert u(t-r,x-z)-u(t+h-r,x+y-z)\right\vert^q\right)\\
&\le  C \int_0^t dr E\left(\left\vert u(t-r,x)-u(t+h-r,x+y)\right\vert^q\right)\\
& = C \int_0^t dr E\left(\left\vert u(r,x)-u(r+h,x+y)\right\vert^q\right).
\end{align*}
H\"older's inequality with respect to Lebesgue measure and the linear growth of the function $b$ imply that
\begin{align*}
T_{2,2}&\le C(q) h^{q-1}\int_t^{t+h} dr \int_{\IR^k} G(r,dz) \left(1+\sup_{t,x}E\left(\vert u(t,x)\vert^q\right)\right)\\
&\le C h^q.
\end{align*}
In the last inequality, we have used (see \cite[Theorem 13]{d99}) that
\beq\label{eq2.7}
\sup_{(t,x)\in [0,T]\times \IR^k}E\left(\vert u(t,x)\vert^q\right) < \infty.
\eeq

Consequently
\beq
\label{5}
T_2\le C\left(\int_0^t dr E\left(\left\vert u(r,x)-u(r+h,x+y)\right\vert^q\right) + h^q\right).
\eeq
The estimates (\ref{3'})--(\ref{5}), along with Gronwall's lemma, yield
\beq
\label{6}
E\left(\vert u(t,x)-u(t+h,x+y)\vert^q\right) \le C \left( \vert h\vert + |y|\right)^{q\delta},
\eeq
with $\delta = \frac{2-\beta}{2}$ under (C1) and $\delta < 1-\alpha$ under (C1'). 
\hfill\qed


\bigskip

\noindent{\bf Case 2: multiplicative noise}
\medskip

These are the known results in dimensions $k\in \{1,2,3\}$.
\medskip

\noindent{\sl Case $k=2$}.
Assume that the covariance measure $\Gamma$ of the noise satisfies (C2).
Then (R2) holds with $\delta\in[0,\frac{\eta}{2}[$ (see \cite[Proposition 1.4]{mss99} and its improvement in   \cite[Theorem 2.2]{mimo}).  
\medskip

\noindent{\sl Case $k=3$}.
Suppose that $\Gamma$ satisfies (C3). For $d=1$, Theorem 4.11 in \cite{dss09}  yields (R2) with $\delta\in\left]0,\frac{2-\beta}{2}\wedge\frac{1+\mu}{2}\right[$. The extension to systems of wave equations is straightforward. 
\medskip

\noindent{\sl Case $k=1$}.
Suppose that $\Gamma$ satisfies (C3). Then (R2) is obtained with $\delta = \frac{2-\beta}{2}$ from a straightforward calculation, much simpler than for $k=2$ or $k=3$, of moments of increments $u(t,x) - u(s,y)$. Since this calculation does not seem to be available in the literature, we sketch it for convenience of the reader.

   Since we are looking for an upper bound, we can work separately with the time-increments and the spatial increments. We only consider the spatial increments, so we assume that $s=t$. Then
$$
    u(t,x) - u(t,y) = T_1 + T_2,
$$
where
\begin{align*}
   T_1&= \int_0^t \int_{\IR^k}(G(t-r,x-z) - G(t-r,y-z))\, \sigma(u(r,z)) \, M(dr,dz),\\
   T_2&= \int_0^t \int_{\IR^k}(G(t-r,x-z) - G(t-r,y-z))\, b(u(r,z)) \, dr,dz,
\end{align*}
Without loss of generality, we assume that $d=1$. Then $E(\vert u(t,x) - u(t,y)\vert^q)$ is bounded by a constant times $E(\vert T_1 \vert^q) + E(\vert T_2 \vert^q)$. Let us consider the term with $T_1$. 

   By Burkholder's inequality,
\begin{align*}
   E(\vert T_1 \vert^q) &\leq C E\left[\left(\int_0^t dr\, \Vert (G(t-r,x-*) - G(t-r,y-*)) \sigma(u(r,*))\Vert_{\hac}^2 \right)^{q/2}\right] \\
      &\leq C_1 \int_0^t dr\, E [\Vert (G(t-r,x-*) - G(t-r,y-*)) \sigma(u(r,*))\Vert_{\hac}^{q}].
\end{align*}
Since $k=1$, $G(r,z) = \frac12 1_{\{\vert z\vert < r \}}$, so H\"older's inequality implies that the $\Vert \cdots \Vert_{\hac}^q$ is bounded by
\begin{align*}
   &\Big(\int_{\IR} dz \int_{\IR} dw\, \vert G(t-r,x-z) - G(t-r,y-z) \vert \, \vert \sigma(u(r,z))\vert\, \varphi(z-w) \\
   &\qquad\qquad  \times  \vert z-w\vert^{-\beta} \vert G(t-r,x-w) - G(t-r,y-w) \vert \, \vert \sigma(u(r,w))\vert\Big)^{\frac{q}{2}-1}\\
   &\qquad \times \int_{\IR} dz \int_{\IR} dw\, \vert G(t-r,x-z) - G(t-r,y-z) \vert \, \varphi(z-w)\\
   &\qquad\qquad \qquad \times \vert z-w\vert^{-\beta} \vert G(t-r,x-w) - G(t-r,y-w) \vert\\
   &\qquad\qquad \qquad\qquad \times E( \vert \sigma(u(r,z))\vert^{q/2} \, \vert \sigma(u(r,w))\vert^{q/2}).
\end{align*}
Since $\varphi$ is bounded and $\sigma$ is Lipschitz, we use the Cauchy-Schwarz inequality to see that
\begin{align*}
E(\vert T_1 \vert^q) &\leq C \int_0^t dr\, \left(1+\sup_{z \in \IR} E(\vert u(r,z)\vert^q)\right)\\
 &\qquad \times \Big(\int_{\IR} dz \int_{\IR} dw\,  \vert G(t-r,x-z) - G(t-r,y-z) \vert \\
 &\qquad\qquad \times \vert z-w\vert^{-\beta} \vert G(t-r,x-w) - G(t-r,y-w) \vert \Big)^{q/2}.
\end{align*}
Assume without loss of generality that $x \leq y$. It is useful to distinguish the cases $t-r < \frac{y-x}{2}$ and $t-r \geq \frac{y-x}{2}$. We consider only the latter case. For $t-r \geq \frac{y-x}{2}$, the integral in parentheses is equal to
\beq\label{eq2.9}
   \int_{A(x,y,t,r)} dz \int_{A(x,y,t,r)} dw \, \vert z-w\vert^{-\beta},
\eeq
where
$$
   A(x,y,t,r) = [x-(t-r), y-(t-r)]\cup [x+t-r,y+t-r].
$$
This integral can be written as the sum of four integrals, the first two of which are
\begin{align*}
   I_1 &= \int_{x-(t-r)}^{y-(t-r)} dz \int_{x-(t-r)}^{y-(t-r)} dw \, \vert z-w\vert^{-\beta},\\
   I_2 &= \int_{x-(t-r)}^{y-(t-r)} dz \int_{x+t-r}^{y+t-r} dw \, \vert z-w\vert^{-\beta}.
\end{align*}
Notice that since $\beta \in \,]0,1[$,
\begin{align*}
   I_1 &= \int_{x-(t-r)}^{y-(t-r)} dz \left(\int_{x-(t-r)}^z dw \, (z-w)^{-\beta}+ \int_z^{y-(t-r)} dw \, (w - z)^{-\beta} \right), 
\end{align*}
and an explicit integration shows that
$$
   I_1 = c\, (y-x)^{2-\beta}.
$$

   In order to bound $I_2$, we distinguish two further cases: $t-r \geq y-x$ and $y-x > t-r \geq \frac{y-x}{2}$. When $t-r \geq y-x$,
$$
   I_2 \leq \int_{x-(t-r)}^{y-(t-r)} dz \int_{x+t-r}^{y+t-r} dw \, \vert y-x\vert^{-\beta} = c\vert y-x\vert^{2-\beta}.
$$
When $y-x > t-r \geq \frac{y-x}{2}$, a direct integration shows that
\begin{align*}
   I_2 &= c [(y-x+2(t-r))^{2-\beta} - 2 (2(t-r))^{2-\beta} + (x-y +2(t-r))^{2-\beta}]\\
   & \leq C \vert y-x\vert^{2-\beta}.
\end{align*}
The other two integrals arising from \eqref{eq2.9} are handled by symmetry. We conclude that
$$
   E(\vert T_1 \vert^q) \leq C \int_0^t dr\,\left(1+ \sup_{z \in \IR} E(\vert u(r,z)\vert^q)\right)\, \vert y-x\vert^{\frac{2-\beta}{2}q}.
$$
Using \eqref{eq2.7}, we see that
$$
   E(\vert T_1 \vert^q) \leq C\, \vert y-x\vert^{\frac{2-\beta}{2}q}.
$$
We leave to the reader to check that a similar bound holds for $E\left(\vert T_2\vert^q\right)$. This leads to the bound
$$
   E(\vert u(t,x) - u(t,y)\vert^q)  \leq C\, \vert y-x\vert^{\frac{2-\beta}{2}q}.
$$
 
 This completes the proof of (R2) with $\delta = \frac{2-\beta}{2}$ under assumption (C3) in the case where $k=1$.
\hfill $\Box$
\bigskip

\subsection{Proof of (R1)}
\label{s1.2}

In this section, we consider Equation (\ref{1.4}) with multiplicative noise. We assume that $k\in\{1, 2, 3\}$, and that the functions $\sigma_{i,j}$ and $b_i$, $i,j=1,\ldots,d$, satisfy (P1) and (P2). 

We shall use Malliavin calculus in the context of  \cite{ss} (see chapters 3 and 7 there for the notions and notation). According to Proposition 3.4 in \cite{dkn09}, property (R1) will be ensured by the following properties:
\begin{enumerate}
\item $u(t,x)\in \mathbb{D}^\infty(\mathbb{R}^d)$.
\item For all $p\in\,]1,\infty[$ and $\ell\ge 1$, there exists a positive constant $c_1=c_1(\ell,p)$ such that
\beqn
\sup_{(t,x)\in \T\times K}E\left(\left\Vert D^{\ell}u_i(t,x)\right\Vert^p_{\mathcal{H}_T^{\otimes \ell}}\right)\le c_1, \ i=1,\ldots,d.
\eeqn
\item For all $p\in\,]1,\infty[$, there exists a positive constant $c_2=c_2(p)$ such that
\beqn
\sup_{(t,x)\in \T\times K}E\left[\left(\det \gamma_{u(t,x)}\right)^{-p}\right]\le c_2,
\eeqn
where $\gamma_{u(t,x)}$ denotes the Malliavin matrix of $u(t,x)$.
\end{enumerate}
In the same two cases as in Theorem \ref{t1}, Properties 1 and 2 above follow from an easy extension of Theorem 7.1 in \cite{ss} to systems of equations. Hence, we shall focus on the proof of Property 3 above. 



For its further use, we write the equation satisfied by the Malliavin derivative of the $i$-th component of $D u(t,x)$.
Remember that this is an $\mathcal{H}^d_T$-valued random variable. According to \cite[Theorem 7.1]{ss} (see also \cite{nq}), for $t\le r$, $D_{r,z}u(t,x)=0$, and for $0\le r\le t$, $j=1,\ldots, d$, 
\begin{align}
&D_{r,z}^{(j)} u_i(t,x)\nonumber \\
 &\qquad =G(t-r,x-z)\sigma_{i,j}(u(r,z))\nonumber\\
&\qquad\qquad +\sum_{\ell=1}^d\int_0^t \int_{\IR^k} G(t-s,x-y)\nabla \sigma_{i,\ell}(u(s,y))\cdot D_{r,z}^{(j)}u(s,y) M^\ell(ds,dy)\nonumber\\
&\qquad\qquad+\int_0^t ds \int_{\IR^k} dy\, G(t-s,x-y) \nabla b_i(u(s,y))\cdot D_{r,z}^{(j)} u(s,y),
\label{9}
\end{align}
where the symbol ``$\cdot$" denotes the Euclidean inner product in $\IR^d$, and the notation $D^{(j)}$ stands for the Malliavin derivative with respect to the $j$-th component of the noise.

\begin{prop}
\label{p2} 
Let  $k\in \{1, 2, 3\}$. Fix, as in Section \ref{s0}, a compact subset $\T\times K\subset \,]0,T]\times \IR^k$.  Assume that the functions $\sigma_{i,j}$, $b_i$, $i,j=1,\ldots,d$, satisfy (P1) and (P2). In the same two cases as in Theorem \ref{t1}, for $p>0$,
\beqn
\sup_{(t,x)\in \T\times K} E\left[\left( \det \gamma_{u(t,x)}\right)^{-p}\right] <+\infty.
\eeqn
\end{prop}

\noindent{\it Proof}: Observe that
\beqn
\det \gamma_{u(t,x)} \ge \left(\inf_{|v|=1} v^{\sT}\gamma_{u(t,x)} v\right)^d.
\eeqn
Set $q=pd$. It suffices to check that
\beqn
\sup_{(t,x)\in \T\times K} E\left[ \left(\inf_{|v|=1} v^{\sT}\gamma_{u(t,x)} v\right)^{-q}\right] <+\infty.
\eeqn

By definition,
\begin{align}
\label{10}
v^{\sT}\gamma_{u(t,x)}v&=\sum_{i,j=1}^d v_i\, \langle Du_i(t,x),Du_j(t,x)\rangle_{\mathcal{H}^d_T}\, v_j\nonumber\\
&=\left\langle \sum_{i=1}^d v_iDu_i(t,x),\ \sum_{j=1}^d v_jDu_j(t,x)\right\rangle_{\mathcal{H}^d_T}\nonumber\\
&=\int_0^t dr \left\Vert \sum_{i=1}^d v_i D_{r,*}u_i(t,x)\right\Vert^2_{\mathcal{H}^d}\nonumber\\
&\ge\int_{t-\delta}^t dr \left\Vert \sum_{i=1}^d v_i D_{r,*}u_i(t,x)\right\Vert^2_{\mathcal{H}^d},
\end{align}
where $\delta>0$ is such that $t-\delta>0$. Now,
$$
   \left\Vert \sum_{i=1}^d v_i D_{r,*}u_i(t,x)\right\Vert^2_{\mathcal{H}^d} = \sum_{j=1}^d \left\Vert\sum_{i=1}^d v_i D^{(j)}_{r,*}u_i(t,x)\right\Vert^2_{\mathcal{H}}.
$$

Set 
\begin{align*}
a_{i,j}(r,z;t,x)&= \sum_{\ell = 1}^d
\int_0^t \int_{\IR^k} G(t-s,x-y)\\
   & \qquad\qquad\qquad \times \nabla \sigma_{i,\ell}(u(s,y))\cdot D_{r,z}^{(j)}u(s,y) \, M^\ell(ds,dy)\\
&\qquad +\int_0^t ds \int_{\IR^k} dy\, G(t-s,x-y) \nabla b_i(u(s,y))\cdot D_{r,z}^{(j)} u(s,y).
\end{align*}
Applying twice the inequality $(a+b)^2 \geq \frac12 a^2 - b^2$ yields
\begin{align*}
\sum_{j=1}^d \left\Vert\sum_{i=1}^d v_iD_{r,*}^{(j)} u_i(t,x)\right\Vert^2_{\mathcal{H}}
&\ge \frac{1}{2} \sum_{j=1}^d \left\Vert\sum_{i=1}^d v_iG(t-r,x-*) \sigma_{i,j}(u(r,*))\right\Vert^2_{\mathcal{H}}\\
&\qquad -\sum_{j=1}^d\left\Vert\sum_{i=1}^d v_i a_{i,j}(r,*;t,x)\right\Vert^2_{\mathcal{H}}\\
&=\frac{1}{2}\sum_{j=1}^d\left\Vert G(t-r,x-*) \left(v^{\sT}\sigma(u(r,*))\right)_j\right\Vert^2_{\mathcal{H}}\\
&\qquad -\sum_{j=1}^d\left\Vert\sum_{i=1}^d v_i a_{i,j}(r,*;t,x)\right\Vert^2_{\mathcal{H}}\\
&\ge \frac{1}{4} \sum_{j=1}^d\left\Vert G(t-r,x-*) \left[v^{\sT}\sigma(u(r,x))\right]_j\right\Vert^2_{\mathcal{H}}\\
&\qquad -T_1(t,x; r) - T_2(t,x; r),
\end{align*}
where
\begin{align*}
T_1(t,x; r)&=\frac{1}{2} \sum_{j=1}^d
\left\Vert G(t-r,x-*)\left[v^{\sT}\left(\sigma(u(r,*))-\sigma(u(r,x))\right)\right]_j\right\Vert^2_{\mathcal{H}},\\
T_2(t,x; r)&=\sum_{j=1}^d\left\Vert\sum_{i=1}^d v_i a_{i,j}(r,*;t,x)\right\Vert^2_{\mathcal{H}}.
\end{align*}
By property (P2), 
\beqn
\sum_{j=1}^d \left\Vert\sum_{i=1}^d v_iD_{r,*}^{(j)} u_i(t,x)\right\Vert^2_{\mathcal{H}}
 \ge \frac{1}{4} \rho^2 \left\Vert G(t-r,x-*)\right\Vert^2_{\mathcal{H}}-T_1(t,x; r) - T_2(t,x; r).
\eeqn

   Define
\beq
\label{g}
g(\delta):=\int_0^\delta\left\Vert G(r,x-*)\right\Vert^2_{\mathcal{H}}dr.
\eeq
Then
$$
   v^{\sT}\gamma_{u(t,x)} v \geq c g(\delta) - \int_{t-\delta}^t (T_1(t,x; r) + T_2(t,x; r)) dr.
$$
We note from Lemma \ref{la1} that under Hypothesis (C1),
$$
   g(\delta) = c \delta^{3 - \beta},
$$
while by \cite[Lemma 8.6, equation (8.35)]{ss}, under (C1'), (C2) and (C3), for small $\delta >0$,
\beq
\label{11}
g(\delta)\ge C \delta^3.
\eeq

   We now bound the moments of $T_1(t,x; r)$. We apply first \eqref{fundamental} and then H\"older's inequality with respect to the finite measure on $[0,\delta]\times\IR^k\times\IR^k$ with density (with respect to Lebesgue measure)
given by $G(r,x-y) f(y-z) G(r,x-z)$. The Lipschitz continuity of $\sigma$ along with the Cauchy-Schwarz inequality yield
\begin{align*}
&E\left[\left(\int_{t-\delta}^t T_1(t,x;r) dr \right)^q\right]\\
&\qquad \le C\left[g(\delta)\right]^{q-1} \int_0^\delta dr \int_{\IR^k} dy \int_{\IR^k} dz\, G(r,x-y) f(y-z) G(r,x-z)\\
&\qquad\qquad\times \left(E\left[\left\vert u(t-r,y)-u(t-r,x)\right\vert^{2q}\right]\right)^{\frac{1}{2}} \\
&\qquad\qquad\times \left(E\left[\left\vert u(t-r,z)-u(t-r,x)\right\vert^{2q}\right]\right)^{\frac{1}{2}}.
\end{align*}
We now apply the results already mentioned on $L^p$-estimates of increments of the process $u(t,x)$ (see the proof of (R2) for multiplicative noise), yielding 
\begin{align*}
&E\left[\left(\int_{t-\delta}^t T_1(t,x;r) dr \right)^q\right]\\
&\qquad \le C\left[g(\delta)\right]^{q-1} \int_0^\delta dr \int_{\IR^k} dy \int_{\IR^k} dz\, G(r,x-y) f(y-z) G(r,x-z)\\
&\qquad\qquad\qquad\qquad\qquad\qquad\qquad\qquad \times \left(|y-x|^{\zeta q}\right)\left(|z-x|^{\zeta q}\right),
\end{align*}
where the values of $\zeta$ depend on the assumptions on the covariance of the noise, as indicated in Cases 1 and 2 of Section \ref{s1.1}: if (C1) holds, then $\zeta \in \,]0,\frac{2-\beta}{2}[$; if (C1') holds, then $\zeta \in \,]0,1-\alpha[$; if $k=2$ and (C2) holds, then $\zeta\in\,]0, \frac{\eta}{2}[$;
if $k=1$ and (C3) holds, then $\zeta\in\,]0,\frac{2-\beta}{2}[$; if $k=3$ and (C3) holds, than $\zeta \in\, ]0,\frac{2-\beta}{2}\wedge\frac{1+\mu}{2}[$.

Next, we notice that the support of $G(r,x-\cdot)$ is $\{y\in \IR^k: |x-y|\le r\}$ for $k\in\{1,2\}$ and $\{y\in \IR^k: |x-y|= r\}$ for $k=3$, respectively, which implies
\beqn
E\left[\left(\int_{t-\delta}^t T_1(t,x;r)\, dr \right)^q\right]\le C \left[g(\delta)\right]^q \delta^{2\zeta q},
\eeqn
and finally
\beq
\label{12}
\sup_{(t,x)\in\mathbb{T}_0\times K}E\left[\left(\int_{t-\delta}^t T_1(t,x;r) dr \right)^q\right]\le C \left[g(\delta)\right]^q \delta^{2\zeta q}.
\eeq

Next, we proceed with the study of the contribution of $T_2(t,x;r)$. Firstly, we consider the term involving a stochastic integral. 
By noticing that $D^{(j)}_{r,*} u(s,y)=0$ unless $0\le r\le s$, we have,
\begin{align*}
&E\left[\left(\int_{t-\delta}^t dr \sum_{j=1}^d
\left\Vert \sum_{i=1}^d v_i \right.\right.\right.\\
&\left.\left.\left.\qquad\times \sum_{\ell=1}^d \int_0^t \int_{\IR^k} G(t-s,x-y)\nabla \sigma_{i,\ell}(u(s,y))\cdot D^{(j)}_{r,*} u(s,y) M^\ell(ds,dy)\right\Vert^2_{\mathcal{H}}\right)^q\right]\\
&\le C  \sup_{1\le i,j,\ell \le d} E\Big[\Big(\int_{t-\delta}^t dr \\
&\qquad\times\left\Vert\int_0^t \int_{\IR^k} G(t-s,x-y)\nabla \sigma_{i,\ell}(u(s,y))\cdot D^{(j)}_{r,*} u(s,y) M^\ell(ds,dy)\right\Vert^2_{\mathcal{H}}\Big)^q\Big]\\
&\le C  \sup_{1\le i,j,\ell\le d} E\Big[\Big(\int_{t-\delta}^t dr \\
&\qquad\times\left\Vert\int_{t-\delta}^t \int_{\IR^k} G(t-s,x-y)\nabla \sigma_{i,\ell}(u(s,y))\cdot D^{(j)}_{r,*} u(s,y) M^\ell(ds,dy)\right\Vert^2_{\mathcal{H}}\Big)^q\Big].
\end{align*}

This is equal to
\begin{align*}
&C  \sup_{1\le i,j,\ell\le d} E\Big(\Big\Vert \int_{t-\delta}^t \int_{\IR^k} G(t-s,x-y)\\
& \qquad\qquad\qquad\qquad\qquad\times\nabla \sigma_{i,\ell}(u(s,y))\cdot D^{(j)}_{\cdot,*} u(s,y) M^\ell(ds,dy)\Big\Vert^{2q}_{\mathcal{H}_\delta}\Big),
\end{align*}
where $\mathcal{H}_\delta=L^2([t-\delta,t];\mathcal{H})$. We use the change of variables $s = t-\delta + \tilde s$ to write this as
\begin{align}\nonumber
  &C  \sup_{1\le i,j,\ell\le d} E\Big(\Big\Vert \int_{0}^\delta \int_{\IR^k} G(\delta-s,x-y) \nabla \sigma_{i,\ell}(u(t-\delta+s,y))\\
& \qquad\qquad\qquad\qquad\cdot D^{(j)}_{\cdot,*} u(t-\delta+s,y) M^\ell(d(t-\delta+s),dy)\Big\Vert^{2q}_{\mathcal{H}_\delta}\Big).
\label{eq2.14}
\end{align}
Using the notation of \cite[Theorem 6.1]{ss}, we set $S(s,y) = G(\delta - s, x-y)$ and
$$
   K_{i,j,\ell}(s,y) = \nabla \sigma_{i,\ell}(u(t-\delta+s,y)) \cdot D^{(j)}_{\cdot,*} u(t-\delta+s,y),
$$
which belongs a.s.~to the Hilbert space $\kac = \mathcal{H}_\delta$. Applying \cite[Theorem 6.1]{ss}, we bound the expression in \eqref{eq2.14} by
$$
   C  [g((\delta)]^{q} \sup_{1\le i,j,\ell\le d}\  \sup_{(s,y)\in[0, \delta]\times\IR^k}E\left(\Vert K_{i,j,\ell}(s,y)\Vert^{2q}_{\mathcal{H}_\delta}\right).
$$
Using the fact that the factor $\nabla \sigma_{i,\ell}(u(t-\delta+s,y))$ does not involve the variables $(\cdot,*)$ that appear in the formula for $\Vert \cdot \Vert_{\mathcal{H}_\delta}$, together with assumption (P1), we see that this expression is bounded by
\begin{align}\nonumber
&C [g((\delta)]^{q} \sup_{1\le j\le d}\sup_{(s,y)\in[t-\delta,t]\times\IR^k}E\left(\Vert D^{(j)}u(s,y)\Vert^{2q}_{\mathcal{H}^d_\delta}\right)\\
&\qquad \le C  [g(\delta)]^{2q},
\label{eq2.15}
\end{align}
where we have used \cite[Lemma 8.2]{ss}.

   Secondly, we consider the pathwise integral. With arguments similar to those applied to the stochastic integral, we obtain
\begin{align*}
&E\Big[\Big(\int_{t-\delta}^t dr \sum_{j=1}^d\\
&\qquad \times\Big\Vert \sum_{i=1}^d v_i \int_0^t ds \int_{\IR^k} dy\, G(t-s,x-y)\nabla b_i(u(s,y))\cdot D^{(j)}_{r,*} u(s,y) \Big\Vert^2_{\mathcal{H}}\Big)^q\Big]\\
&\le C  \sup_{1\le i,j\le d}\\
&\qquad E\left(\left\Vert \int_{0}^\delta ds \int_{\IR^k} dy\, G(s,x-y)\nabla b_i(u(t-s,y))\cdot D^{(j)}_{\cdot,*} u(t-s,y) \right\Vert^{2q}_{\mathcal{H}_\delta}\right).
\end{align*}
Minkowski's inequality applied to the norm $\Vert \cdot\Vert_{{\mathcal{H}}_\delta}$and then H\"older's inequality with respect to the measure on $[0,\delta]\times\IR^k$ given by  $G(s,dy) ds$, yield
the following upper bound for the last expression:
\begin{align}\nonumber
&C \sup_{1\le i,j\le d}
E\left[\left(\int_{0}^\delta ds \int_{\IR^k} dy\, G(s,x-y)\left\Vert D^{(j)}_{\cdot,*} u(t-s,y) \right\Vert_{{\mathcal{H}}_\delta}\right)^{2q}\right]\\ \nonumber
&\le C  \left(\int_{0}^\delta ds \int_{\IR^k} dy\, G(s,x-y)\right)^{2q}\\ \nonumber
&\qquad\qquad\times \sup_{1\le j\le d}\ \sup_{(s,y)\in[t-\delta,t]\times \IR^k} E\left(\left\Vert D^{(j)}u(s,y)\right\Vert^{2q}_{{\mathcal{H}}_\delta}\right)\\
&\le C  \delta^{4q} [g(\delta)]^q,
\label{eq2.16}
\end{align}
by Lemmas 8.6 and 8.2 in \cite{ss}.

    We conclude from \eqref{eq2.15} and \eqref{eq2.16} that
\beq\label{eq2.17}
 \sup_{(t,x)\in\mathbb{T}_0\times K}E\left[\left(\int_{t-\delta}^t T_2(t,x;r) dr \right)^q\right]\le C  \left[g(\delta)\right]^q ([g(\delta)]^q + \delta^{4 q}).
\eeq

   For small $\varepsilon >0$, choose $\delta_\varepsilon$ such that $g(\delta_\varepsilon) = \varepsilon$. By \eqref{12} and \eqref{eq2.17},  for any $\varepsilon >0$ sufficiently small (smaller than some $\varepsilon_0 >0$, say),
$$
   \inf_{\vert v\vert =1} v^t\gamma_{u(t,x)}v \ge c \varepsilon - \int_{t - \delta_\varepsilon}^t T(t,x;r)\, dr
$$
where
$$
   \sup_{(t,x)\in\mathbb{T}_0\times K}E\left[\left( \int_{t-\delta}^t T(t,x;r) dr \right) ^q\right] \le C \varepsilon^q \delta_\varepsilon^{2\zeta  q} + C  \varepsilon^q\, (\varepsilon^q + \delta_\varepsilon^{4q}).
$$
By \eqref{11}, $\delta_\varepsilon \le \varepsilon^{1/3}$, so the right-hand side is less than or equal to $\varepsilon^{(1+\rho)q}$, with $\rho >0$. Therefore, we can apply \cite[Proposition 3.5]{dkn09} to conclude the proof of Proposition \ref{p2}. 
\hfill\qed
\medskip

\noindent{\it Proof of Theorem \ref{t1}}. As has been already mentioned, we apply   \cite[Theorem 2.4]{dss10}. Then the conclusion in Case 1 follows from Propositions \ref{p1} and \ref{p2}, while those in Case 2 are a consequence of
the results stated in Section \ref{s1} under the title {\it Case 2: multiplicative noise,} together with the three properties stated at the beginning of Section \ref{s1.2}.

\hfill\qed

\section{Conditions on eigenvalues of the Malliavin matrix for lower bounds}
\label{s2}


In this section, we consider the system of equations given in (\ref{1.4}). We restrict ourselves to a covariance measure $\Gamma(dx)=\vert x\vert^{-\beta}dx$, $\beta\in\,]0,2\wedge k[$ (hypothesis (C1) in Section \ref{s1}), with $k \in \{ 1,2,3\}$. 

We recall the equation \eqref{9} satisfied by the Malliavin derivative of the components $u_i(t,x)$, $i=1,\ldots,d$, of the random field $u(t,x)$:
\beq
D_{r,z}^{(j)}u_i(t,x)=G(t-r,x-z)\sigma_{i,j}(u(r,z))+a_{i,j}(r,z;t,x),
\eeq
with
\begin{align}
a_{i,j}(r,z;t,x)&= \sum_{\ell=1}^d\int_0^t \int_{\IR^k} G(t-s,x-y)\nabla \sigma_{i,\ell}(u(s,y))\cdot D_{r,z}^{(j)}u(s,y) M^\ell(ds,dy)\nonumber\\
&\qquad+\int_0^t ds \int_{\IR^k} dy\, G(t-s,x-y) \nabla b_i(u(s,y))\cdot D_{r,z}^{(j)} u(s,y).\label{as}
\end{align}

The Malliavin covariance matrix of 
\beq
\label{Z}
Z:=(u_1(s,y),\ldots,u_d(s,y), u_1(t,x)-u_1(s,y),\ldots,u_d(t,x)-u_d(s,y)),
\eeq
is defined by
\beqn
\gamma_Z=\langle DZ,DZ\rangle_{\mathcal{H}_T^d}.
\eeqn

As has been proved in \cite{dkn09}, Lemma 6.8, there exist at least $d$ orthonormal eigenvectors $\xi^{i_1},\cdots,\xi^{i_d}$ of $\gamma_Z$, with corresponding eigenvalues
$\alpha_{i_1}, \dots,\alpha_{i_d}\ge \alpha_0$, with $\alpha_0$ strictly positive and deterministic. We fix a set $K\subset\{1,\dots,2d\}$ with ${\rm card}\ K=d$ and set
\beq
\label{evl}
A_K=\cap_{i\in K}\{\alpha_i\ge \alpha_0\}.
\eeq
We term {\it large eigenvalues} the $\alpha_i$ with $i\in K$.

The goal here is to prove that lower bounds for the hitting probabilities can be derived from two properties on the eigenvalues of the matrix $\gamma_Z$: estimates on the negative $L^p$-moments of the smallest eigenvalue of the matrix $\gamma_Z$ and boundedness of the same type of moments for the large eigenvalues. Establishing such properties on $\gamma_Z$ is postponed to the next section. We aim at proving the following theorem.

\begin{thm}
\label{t2.10} Fix $k\in\{1,2,3\}$. 
Let $I=[a,b]\subset [0,T]$, with $a>0$, $K$ a compact subset of $\rk$ and fix $N>0$. Assume (P1), (P2) and (C1). Suppose also that there exists $\rho\in\,]0,\infty[$ such that for any $p\in[1,\infty[$, there is a positive constant $C$ such that for all $s,t\in I$ and $x,y\in K$ with
$(s,y)\ne(t,x)$, $\vert t-s \vert \leq 1$ and $\vert x-y\vert \leq 1$, the following properties hold:
\begin{align}
&E\left[\left(\inf_{\vert \xi\vert=1}\xi^T \gamma_Z\xi\right)^{-p}\right]\le C \left(|t-s|+|x-y|\right)^{-p\rho}, \label{2.111m}\\
&E\left(1_{A_K}\left(\Pi_{i\in K}(\xi^i)^T\gamma_Z\xi^i\right)^{-p}\right) \le C.\label{LEV}
\end{align}

\begin{description}
\item{(1)} Fix $\delta>0$. There exists a positive constant $c=c(I,K,N,\beta,k,d,\rho,\delta)$ such that, for any
Borel set $A\subset[-N,N]^d$,
\beq
\label{t2.111}
P\left\{u(I\times K)\cap A\ne\emptyset\right\}\ge c\ {\rm Cap}_{d\left(1+\frac{4d(\rho-(2-\beta))}{2-\beta}\right)+\delta-\frac{2(k+1)}{2-\beta}}(A).
\eeq

\item{(2)} Fix $\delta>0$ and $x\in K$. There exists a positive constant $c=c(I,x,N,\beta,k,d,\rho,\delta)$ such that, for any
Borel set $A\subset[-N,N]^d$,
\beq
\label{t2.1111}
P\left\{u(I\times \{x\})\cap A\ne\emptyset\right\}\ge c\ {\rm Cap}_{d\left(1+\frac{4d(\rho-(2-\beta))}{2-\beta}\right)+\delta-\frac{2}{2-\beta}}(A).
\eeq
\item{(3)} Fix $\delta>0$ and $t\in I$. There exists a positive constant $c=c(t,K,N,\beta,k,\rho,d)$ such that, for any
Borel set $A\subset[-N,N]^d$,
\beq
\label{t2.1112}
P\left\{u(\{t\}\times K)\cap A\ne\emptyset\right\}\ge c\ {\rm Cap}_{d\left(1+\frac{4d(\rho-(2-\beta))}{2-\beta}\right)+\delta-\frac{2k}{2-\beta}}(A).
\eeq
\end{description}
\end{thm}

  The proof of this theorem is presented in Section \ref{ss2.2}.
  


\subsection{Upper bounds for joint densities}
\label{ss2.1}

   This section establishes several preliminary results, in particular Theorems \ref{tH} and \ref{rdlem1}. Indeed, according to Theorem 2.1 in \cite{dss10} (see also Hypothesis H2 in \cite{dn04}), Theorem \ref{rdlem1} is a fundamental step toward establishing lower bounds on the hitting probabilities.

Under suitable regularity and non-degeneracy conditions on a random vector, the integration by parts formula of Malliavin calculus provides an expression for its density. 
We may refer for instance to Corollary 3.2.1 in \cite{n2}. Consider the random vector $Z$ defined in \eqref{Z} and denote by $p_{s,y;t,x}(z_1,z_2)$ its density at 
$(z_1,z_1-z_2)$. By applying H\"older's inequality, we have
\begin{align}
\label{2.1}
p_{s,y;t,x}(z_1,z_2)&\le \Pi_{i=1}^d\left(P\{\vert u_i(t,x)-u_i(s,y)\vert>\vert z_1^i-z_2^i\vert\}\right)^{\frac{1}{2d}}\nonumber\\
&\quad\times \Vert H_{(1,\ldots,2d)}(Z,1)\Vert_{L^2(\Omega)}
\end{align}
(see \cite{dkn09}, page 395). 
In \eqref{2.1}, $H_{(1,\ldots,2d)}(Z,1)$ is the random variable defined recursively as follows. For any $i=1,\ldots, 2d$, and any random variable $G$,
\beq
\label{hh1}
H_{(i)}(Z,G)= \sum_{j=1}^{2d}\delta\left(G\left(\gamma^{-1}_Z \right)_{i,j}DZ^j\right),
\eeq
and then, for every integer $j\ge 2$,
\beq
\label{hh2}
H_{(1,\ldots,j)}(Z,1)= H_{(j)}(Z,H_{(1,\ldots,j-1)}(Z,1)),
\eeq
where $\delta$ denotes the divergence operator (Skorohod integral) and $\gamma_Z$ the Malliavin covariance matrix (see for instance \cite{ss} for a definition of these notions).
\bigskip

As previously, we fix compact sets $I\subset\, ]0,T]$, $K\subset \rk$
and assume that $s,t\in I$ and $x,y\in K$. 
One of the main objectives of this section is to prove the following result.
\begin{thm}
\label{tH}
The hypotheses are the same as in Therorem \ref{t2.10}. Then, for every $\delta>0$, there exists a positive constant $C$ such that
\beq
\label{H}
 \Vert H_{(1,\ldots,2d)}(Z,1)\Vert_{L^2(\Omega)}\le \frac{C}{(|t-s|+|x-y|)^\gamma},
 \eeq
with $\gamma=\frac{d}{2}(2-\beta+4d(\rho-(2-\beta)))+\delta$ and $\rho$ as in \eqref{2.111m}.
\end{thm}
\bigskip

Theorem \ref{tH} is the analogue in the context of 
systems of wave equations of the estimate (6.3) in \cite{dkn09} for systems of heat equations in spatial dimension $k=1$. In the next paragraphs, we shall describe briefly the method of the proof, following \cite{dkn09}. We refer the reader to this reference for further details. 

By the definition \eqref{hh1}--\eqref{hh2}, we see that  $\Vert H_{(1,\ldots,2d)}(Z,1)\Vert_{L^2(\Omega)}$ can be estimated by applying $L^p$-bounds for the Skorohod integral and H\"older's inequality for the Watanabe-Sobolev norms $\Vert\cdot\Vert_{m,p}$, $m\ge 1$, $p\in[1,\infty[$ (see \cite[Proposition 3.2.1]{n1} and \cite[Proposition 1.10, p.50]{w1}). By doing so, we obtain
 \begin{align*}
 &\Vert H_{(1,\ldots,2d)}(Z,1)\Vert_{L^2(\Omega)}\le C \Vert H_{(1,\ldots,2d-1)}(Z,1)\Vert_{1,4}\\
 &\quad\times \left\{\sum_{j=1}^d\left(\left\Vert\left(\gamma^{-1}_Z\right)_{2d,j}\right\Vert_{1,8}\left\Vert DZ^j\right\Vert_{1,8}\right)
 +\sum_{j=d+1}^{2d} \left(\left\Vert\left(\gamma^{-1}_Z\right)_{2d,j}\right\Vert_{1,8}\left\Vert DZ^j\right\Vert_{1,8}\right)\right\}.
 \end{align*}
 By iterating this procedure $2d-1$ more times, we go down from $H_{(1,\ldots,2d-1)}(Z,1)$ to $H_{(0)}(Z,1):=1$. Eventually, we get the estimate 
 \beq
 \label{hh3}
 \Vert H_{(1,\ldots,2d)}(Z,1)\Vert_{L^2(\Omega)}\le C  \prod_{\ell=1}^{2d} \mathcal{Z}_\ell,
 \eeq
where $\mathcal{Z}_\ell$, $\ell=1,\ldots,2d$, is given by
 \beq
 \label{z}
\sum_{j=1}^d\left(\left\Vert\left(\gamma^{-1}_Z\right)_{\ell,j}\right\Vert_{m_\ell,k_\ell}\left\Vert DZ^j\right\Vert_{m_\ell,k_\ell}\right)
 +\sum_{j=d+1}^{2d} \left(\left\Vert\left(\gamma^{-1}_Z\right)_{\ell,j}\right\Vert_{m_\ell,p_\ell}\left\Vert DZ^j\right\Vert_{m_\ell,p_\ell}\right),
 \eeq
 for some $m_\ell\ge 1$ and $p_\ell\in[1,\infty[$.
 
Therefore, \eqref{H} will be a consequence of the next two propositions. 
 \begin{prop}
\label{pz1}
Assume (P1). Let $Z(s,y;t,x)=(u(s,y), u(t,x)-u(s,y))$. For any integer $m\ge 1$ and any $p\in[1,\infty[$, the following statements hold.
\begin{enumerate}
\item For any $j\in\{1,\ldots, d\}$, $\sup_{(s,y), (t,x)\in[0,T]\times \IR^k}\left\Vert Z^j(s,y;t,x)\right\Vert_{m,p}<\infty$.
\item For any $j\in\{d+1,\ldots, 2d\}$, and any $\alpha\in\, ]0,\frac{2-\beta}{2}[$, there exists a positive constant $C$ (depending on $I$ and $K$), such that for any $s,t\in I$ and $x,y\in K$,
\beqn
\left\Vert Z^j(s,y;t,x)\right\Vert_{m,p}\le C (|t-s|+|x-y|)^\alpha.
\eeqn

\end{enumerate}
\end{prop}
\noindent{\it Proof.} In the particular case $d=1$, the conclusion of part 1 is Theorem 7.1 of \cite{ss}. The extension to
arbitrary $d$ is straightforward.
The statement in part 2 follows from an extension of the results proved in \cite{dss09}. Indeed, as has been pointed out in \cite{ss08}, the results of \cite{dss09} can be extended to the Hilbert-space-valued solutions of the equations satisfied by the Malliavin derivatives.
\hfill\qed
\vskip 16pt

For the next statement, we introduce the following notations:
 \begin{align*} 
(\rm 1)&=  \{1,2,\ldots,d\}^2,\\
(\rm 2)&=\{1,2,\ldots,d\} \times \{d+1,\ldots,2d\},\\
(\rm 3)&=\{d+1,\ldots,2d\} \times \{1,2,\ldots,d\},\\
(\rm 4)&=\{d+1,\ldots,2d\}^2.
\end{align*}
These four sets form a partition of $I_d:=\{1,2,\ldots,2d\}^2$.

\begin{prop}
\label{pz2}
The hypotheses are the same as in Theorem \ref{t2.10}. For any $m\ge 1$, $p\in[1,\infty[$ and every $\eta>0$, there exists a constant $C>0$ such that for any $s,t\in I$ and $x,y\in K$,
\beq
\label{gammaz}
\left\Vert\left(\gamma^{-1}_Z\right)_{i,j}\right\Vert_{m,p}\le C
\begin{cases}
(|t-s|+|x-y|)^{-d(\rho-(2-\beta))-\eta}, & (i,j)\in (1),\\
(|t-s|+|x-y|)^{-d(\rho-(2-\beta))-\left(\frac{2-\beta}{2}+\eta\right)}, & (i,j)\in (2),  (3),\\
(|t-s|+|x-y|)^{-d(\rho-(2-\beta))-2\left(\frac{2-\beta}{2}+\eta\right)}, & (i,j)\in (4).
\end{cases}
\eeq
\end{prop}

Assuming the validity of this proposition, let us deduce \eqref{H}.
\bigskip

\noindent{\em Proof of Theorem \ref{tH}.}
In the next arguments, the values of the parameters $\eta,\bar\eta>0$ may vary from one inequality to another; they denote arbitrary small positive numbers.

From \eqref{z}, Propositions \ref{pz1} and \ref{pz2}, we obtain
\begin{align*}
&\Vert H_{(1,\ldots,2d)}(Z,1)\Vert_{L^2(\Omega)}\le C\left(\prod_{\ell=1}^d\mathcal{Z}_\ell\right)\left(\prod_{\ell=d+1}^{2d}\mathcal{Z}_\ell\right)\\
&\qquad\le C\Big[\left(|t-s|+|x-y|\right)^{-d(\rho-(2-\beta))-\eta}\\
&\qquad\qquad\qquad  + \left(|t-s|+|x-y|\right)^{-d(\rho-(2-\beta))-\left(\frac{2-\beta}{2}+\eta\right)+\left(\frac{2-\beta}{2}-\bar\eta\right)}\Big]^d\\
& \qquad\qquad\times\left[\left(|t-s|+|x-y|\right)^{-d(\rho-(2-\beta))-\left(\frac{2-\beta}{2}+\eta\right)}\right.\\
&\left.\qquad\qquad\qquad + \left(|t-s|+|x-y|\right)^{-d(\rho-(2-\beta))-2\left(\frac{2-\beta}{2}+\eta\right)+\left(\frac{2-\beta}{2}-\bar\eta\right)}\right]^d,
\end{align*}
where we have applied Propositions \ref{pz1} and \ref{pz2} above. 
This is bounded by $\left(|t-s|+|x-y|\right)^{-2d^2(\rho-(2-\beta))-d\frac{2-\beta}{2}-\eta}$, and this yields \eqref{H}. \hfill\qed
\medskip

Before giving the proof of Proposition \ref{pz2}, let us go into the structure of the term $\gamma^{-1}_Z$. By its very definition, 
\beq
\label{gammaz1}
\left\Vert\left(\gamma^{-1}_Z\right)_{i,j}\right\Vert_{m,p}=\left\{E\left(\left\vert\left(\gamma_Z^{-1}\right)_{i,j}\right\vert^p\right)+
\sum_{\ell=1}^m E\left[\left\Vert D^\ell\left(\gamma_Z^{-1}\right)_{i,j}\right\Vert^p_{{\hacd}^{\otimes k}}\right]\right\}^{\frac{1}{p}},
\eeq
and
\beq
\label{gammaz2}
\left(\gamma_Z^{-1}\right)_{i,j}=\left(\det \gamma_Z\right)^{-1}\left(A_Z\right)_{i,j},
\eeq
where $A_Z$ is the cofactor matrix of $\gamma_Z$.

\begin{prop}
\label{pz3}
Assume (P1). Then for any 
$p\in[1,\infty[$ and every
$\alpha\in\left]0,\frac{2-\beta}{2}\right[$, there exists a constant $C>0$ such that for any $s,t\in I$ and $x,y\in K$,
\beqn
\left\Vert\left(A_Z\right)_{i,j}\right\Vert_{L^p(\Omega)}\le C
\begin{cases}
(|t-s|+|x-y|)^{2d\alpha}, & (i,j)\in(1),\\
(|t-s|+|x-y|)^{(2d-1)\alpha}, & (i,j)\in(2), (3),\\
(|t-s|+|x-y|)^{(2d-2)\alpha}, & (i,j)\in(4).
\end{cases}
\eeqn
\end{prop}
\noindent{\it Proof.} We follow the method of the proof of Proposition 6.5 in \cite{dkn09} and we see that the result is a consequence of
the estimates given in Proposition \ref{pz1} along with the upper bounds for the $L^p$-moments of $\left(A_Z\right)_{i,j}$ for the different sets of indices $(i,j)$ shown in \cite[p.400-401]{dkn09}. \hfill\qed
\medskip

\begin{prop}
\label{pz4}
The hypotheses are those of Theorem \ref{t2.10}. 
Then for any 
$p\in[1,\infty[$, there exists a constant $C>0$ such that for any $s,t\in I$ and $x,y\in K$ with $(t,x)\ne(s,y)$,
\beq
\label{gammaz3}
\left\Vert\left(\det \gamma_Z\right)^{-1}\right\Vert_{L^p(\Omega)} \le C\left(|t-s|+|x-y|\right)^{-d\rho}.
\eeq
\end{prop}
\noindent{\it Proof.} Quoting Lemma 6.8 and inequality (6.11) in \cite{dkn09}, we write
\begin{align*} 
&\left\Vert\left(\det \gamma_Z\right)^{-1}\right\Vert_{L^p(\Omega)}\le \sum_{K\subset\{1,2,\ldots,2d\}\atop |K|=d}
\left(E\left[1_{A_K}\left(\Pi_{i\in K}(\xi^i)^T\gamma_Z\xi^i\right)^{-2p}\right]\right)^{\frac{1}{2p}}\\
&\qquad\qquad\qquad\qquad\qquad\qquad\times 
\left(E\left[\left(\inf_{\xi=(\lambda,\mu)\in{\mathbb R}^{2d}\atop \vert\lambda\vert^2+\vert\mu\vert^2=1}\xi^T\gamma_Z\xi\right)^{-2pd}\right]\right)^{\frac{1}{2p}},
\end{align*}
with $A_K=\cap_{i\in K}\{\alpha_i\ge\alpha_0\}$ defined in \eqref{evl}. 

According to (\ref{LEV}), the first factor on the right-hand side of the preceding inequality is bounded by a constant.  Then \eqref{gammaz3} follows from the assumption on the small eigenvalues given in \eqref{2.111m}. \hfill\qed

\begin{prop}
\label{pz5}
We assume the hypothesis (P1). For any $p\in[1,\infty[$, $\ell\ge 1$ and $\alpha\in\, ]0,\frac{2-\beta}{2}[$, there exists a constant $C>0$ such that for any $s,t\in I$ and $x,y\in K$,
\beqn
\left\Vert D^\ell\left(\gamma_Z\right)_{i,j}\right\Vert_{p,k}\le C
\begin{cases}
1, & (i,j)\in(1),\\
(|t-s|+|x-y|)^\alpha, & (i,j)\in(2), (3),\\
(|t-s|+|x-y|)^{2\alpha}, & (i,j)\in(4).
\end{cases}
\eeqn
\end{prop}
\noindent{\it Proof.} This proposition is the analogue for systems of wave equations of \cite[Proposition 6.7]{dkn09}. As in this reference, the result follows by
applying the Leibniz rule for the Malliavin derivative and Proposition \ref{pz1} (see pages 402-403 in \cite{dkn09}).\hfill\qed
\bigskip

\noindent{\em Proof of Proposition \ref{pz2}.}
We follow the same method as in the proof of Theorem 6.3 in \cite{dkn09}. 
Consider first the case $m=0$. Applying \eqref{gammaz2}, then H\"older's inequality along with Propositions \ref{pz3} and \ref{pz4}, we obtain \eqref{gammaz}.

For $m\ge 1$, we apply a recursive argument. For this, we consider first the case $m=1$ and the identity
\beqn
D\left(\gamma_Z^{-1}\right) = -\gamma_Z^{-1}\left(D\left(\gamma_Z\right)\right)\gamma_Z^{-1},
\eeqn
which is a consequence of the chain rule of Malliavin calculus and the trivial identity $\gamma_Z^{-1}\gamma_Z= Id$. The estimates \eqref{gammaz}
follow from what has been already proved for $m=0$ and Proposition \ref{pz5}. \hfill\qed

\bigskip

With \eqref{2.1} and Theorem \ref{tH}, we can now give the second main result of this section.

\begin{thm} 
\label{rdlem1}
Let $\beta\in\, ]0,2\wedge k[$. Suppose \eqref{H} holds for some $\gamma\in[k+1,\infty[$.
Fix $\alpha \in \,]0,\frac{2-\beta}{2}[$. For any $I$, $K$ compact subsets of $[0,T]$ and $\IR^k$ respectively, both with diameter $\leq 1$, and every $z_1,z_2 \in \IR^d$ with $0<|z_1 - z_2|\le N$, there exists a constant $C:=C(\alpha,\gamma,k,d,N)$ such that
\begin{equation}\label{rd3.31}
 \mathcal{I}:=\int_{I\times K} dt dx \int_{I\times K} ds dy \, p_{s,y;t,x}(z_1,z_2) \le C\, K_{\frac{\gamma - (k+1)}{\alpha}}(\vert z_1 - z_2 \vert).
\end{equation}
\end{thm}
{\it Proof}. 
For $i=1,\dots,d$, let $\eta^i = z_1^i - z_2^i$. Assume without loss of generality that $0< \sup_i\, \vert \eta^i \vert = \vert \eta^1 \vert$. Let $\rho_0 >0$ be such that $I\times K \subset B_{\frac{\rho_0}{2}}(0)$. Recalling \eqref{2.1} and Theorem \ref{tH}, we see that
$$
  \mathcal{I} \leq \mathcal{I}_1 + \mathcal{I}_2,
$$
where
\begin{align*}
  \mathcal{I}_1 &= \int_{I\times K} dt dx \int_{I\times K} ds dy \, 1_{\{(|t-s|+|x-y|) \vert \eta \vert^{-\frac{1}{\alpha}}\leq \rho_0 N^{-\frac{1}{\alpha}}\}} \frac{C}{(|t-s|+|x-y|)^\gamma}\\
  &\qquad\qquad\times \left[P\{\vert u_1(t,x) - u_1(s,y) \vert > \vert \eta^1 \vert\} \right]^{\frac{1}{2d}} ,\\
\mathcal{I}_2 &= \int_{I\times K} dt dx \int_{I\times K} ds dy \, 1_{\{(|t-s|+|x-y|) \vert \eta \vert^{-\frac{1}{\alpha}} > \rho_0 N^{-\frac{1}{\alpha}}\}} \frac{C}{(|t-s|+|x-y|)^\gamma}.
\end{align*} 

  Apply Chebychev's inequality to see that for $p>0$, 
\begin{align*}
  \mathcal{I}_1 &\leq \int_{I\times K} dt dx \int_{I\times K} ds dy \, 1_{\{(|t-s|+|x-y|) \vert \eta \vert^{-\frac{1}{\alpha}}\leq \rho_0 N^{-\frac{1}{\alpha}}\}} \frac{C}{(|t-s|+|x-y|)^\gamma}\\
  &\qquad\qquad\times \left[\frac{E(\vert u_1(t,x) - u_1(s,y) \vert^p)}{\vert \eta^1 \vert^p} \right]^{\frac{1}{2d}} ,
\end{align*} 
By the result of \cite[Theorem 4.1]{dss09} on H\"older continuity,
$$
\mathcal{I}_1 \leq \int_{I\times K} dt dx \int_{I\times K} ds dy \, 1_{\{(|t-s|+|x-y|) \vert \eta \vert^{-\frac{1}{\alpha}}\leq \rho_0 N^{-\frac{1}{\alpha}}\}} \frac{(|t-s|+|x-y|)^{-\gamma+\frac{\alpha p}{2d}}}{\vert \eta^1 \vert^{\frac{p}{2d}}} .
$$
In the second integral, we do the change of variables $\theta = ((t,x)- (s,y))\vert \eta \vert^{-1/\alpha}$, so ``there is no longer dependence on the first integral", in order to see that
$$
\mathcal{I}_1 \leq C \frac{\vert \eta \vert^{\frac{k+1-\gamma}{\alpha}+\frac{p}{2d}}}{\vert \eta^1\vert^{\frac{p}{2d}}} \int_{B_{\rho_0 N^{-\frac{1}{\alpha}}}(0)} d\theta\, \vert \theta \vert^{-\gamma +\frac{\alpha p}{2d}}.
$$
Since $\frac{\vert \eta \vert}{\vert \eta^1 \vert} = \frac{\vert \eta \vert}{\vert \eta \vert_\infty} \leq 1$, we choose $p > (\gamma - k - 1)\frac{2d}{\alpha}$ so that the integral is finite, and we see that
\begin{equation}\label{rd3.32}
  \mathcal{I}_1 \leq C \vert \eta \vert^{\frac{k+1-\gamma}{\alpha}}.
\end{equation}

   Concerning $\mathcal{I}_2$, we use the same change of variables as above to see that
\begin{equation}\label{rd3.33}
  \mathcal{I}_2 \leq C \vert \eta \vert^{\frac{k+1-\gamma}{\alpha}} \int_{B_{\rho_0 |\eta|^{-\frac{1}{\alpha}}}(0)\setminus B_{\rho_0 N^{-\frac{1}{\alpha}}}(0)} d\theta \, \vert \theta \vert^{-\gamma}
\end{equation}

  In the case where $\gamma > k+1$, we bound the integral by
$$
 \int_{\rho_0 N^{-\frac{1}{\alpha}}}^\infty d\rho \, \rho^{-\gamma+k},
$$
which is finite since $\gamma > k+1$. Combining \eqref{rd3.32} and \eqref{rd3.33}, we obtain \eqref{rd3.31} in this case. 

  In the case where $\gamma = k+1$, we bound the integral in \eqref{rd3.33} by
\begin{equation}\label{rd3.34}
  \int_{\rho_0 N^{-\frac{1}{\alpha}}}^{\rho_0|\eta|^{-\frac{1}{\alpha}}} d\rho \, \rho^{-1} = \frac{1}{\alpha} \ln\left(\frac{1}{\vert\eta\vert}\right) + \frac{1}{\alpha} \ln N.
\end{equation}
Combining \eqref{rd3.32}, \eqref{rd3.33} and \eqref{rd3.34}, we obtain \eqref{rd3.31} in this case.

 The theorem is proved.
\hfill $\Box$
\vskip 16pt    

\begin{rem} (a) In the case where $\sigma$ is constant, it is possible to take $\alpha = \frac{2-\beta}{2}$ in Theorem \ref{rdlem1}, because of Proposition \ref{p1}.

  (b) A sufficient condition on the density that could replace the hypothesis (1) in \cite[Theorem 2.1]{dss10} would be
$$
  p_{x,y}(z_1,z_2) \leq \frac{1}{\vert x-y\vert^\gamma}\left[ \frac{\vert x-y\vert^\alpha}{\vert z_1 - z_2 \vert} \wedge 1\right]^p ,
$$
where $p > (\gamma - m)\frac{2d}{\alpha}$. The conclusion of \cite[Theorem 2.1]{dss10} would then be
$$
   P\{v(I) \cap A \neq \emptyset \} \geq C\, {\rm Cap}_{\frac{1}{\alpha}(\gamma - m)}(A)
$$
instead of (12) there.
\end{rem}
\bigskip

\subsection{Proof of Theorem \ref{t2.10}}
\label{ss2.2}

 We begin with the proof of (1). Set $\gamma=\frac{d}{2}(2-\beta+4d(\rho-(2-\beta)))+\delta$. Since $\delta >0$ is arbitrary, it suffices to consider the case where
$$
   {\rm Cap}_{\frac{2(\gamma -(k+1))}{2-\beta} + \delta}(A) >0,
$$
in which case $\frac{2(\gamma -(k+1))}{2-\beta} < d$. Since $P\{u(I\times K) \cap A \neq \emptyset\}$ decreases if we replace $I$ by a subset of $I$, we can assume that diam$(I) \leq 1$ and diam$(K) \leq 1$. We assume first that $A$ is compact and we consider two different cases.
\bigskip

\noindent {\bf Case 1: $ \gamma<k+1$.}
By the definition of the capacity, ${\text {Cap}}_{\frac{2(\gamma-(k+1))}{2-\beta}}(A)=1$. Thus, it suffices to check that
\beqn
P\{u(I\times K)\cap A \ne \emptyset\}\ge c,
\eeqn
for some positive constant $c$. For this, we proceed in a manner similar to the proof of \cite[Theorem 2.1]{dss10}.  Fix $z\in A$, $\ep\in\, ]0,1[$ and set
\beqn
J_\ep(z)= \frac{1}{(2\ep)^d}\int_{I\times K}dt dx\, 1_{B_\ep(z)}(u(t,x)).
\eeqn
We will prove that $E(J_\ep(z))\ge c_1$ and $E[(J_\ep(z))^2]\le c_2$ for some positive constants $c_1$, $c_2$.
Then, by using Paley-Zygmund inequality,
\beqn
P\{u(I\times K)\cap A^{(\ep)} \ne \emptyset\}\ge P\{J_\ep(z)>0\}\ge \frac{[E(J_\ep(z))]^2}{E[(J_\ep(z))^2]}\ge c,
\eeqn
where $A^{(\ep)}=\{x\in\mathbb{R}^d: d(x,A)<\ep\}$.

The lower bound for $E(J_\ep(z))$ is a direct consequence of the results proved in \cite{en10} (see also \cite{chss01}) on the positivity of the density for the solution of \eqref{1.4}. To establish
the upper bound, we apply \eqref{2.1} along with \eqref{H} and obtain
\begin{align*}
E[(J_\ep(z))^2]&=\frac{1}{(2\ep)^{2d}}\int_{I\times K}dt dx \int_{I\times K}ds dy\int_{B_\ep(z)\times B_\ep(z)} dz_1 dz_2\, p_{s,y;t,x}(z_1,z_2)\\
&\le C\int_{I\times K}dt dx \int_{I\times K}ds dy\, \frac{1}{(|t-s|+|x-y|)^\gamma}.
\end{align*}
We may assume that $I\times K$ is included in the $k+1$-dimensional ball $B_{r_0}(0)$, for some $r_0>0$. Thus,
changing to polar coordinates, we obtain
\begin{align*}
E[(J_\ep(z))^2]&\le C\int_{B_{r_0}(0)\times B_{r_0}(0)} dz dz^\prime\,  \frac{1}{|z-z^\prime|^\gamma}\\
&\le C \int_0^{2r_0} \rho^{k-\gamma} d\rho <\infty.
\end{align*}
This ends the proof of the theorem when $\gamma<k+1$.
\bigskip



\noindent{\bf Case 2: $0\le\frac{2(\gamma-(k+1))}{2-\beta}<d$.}
In the analysis of Case 1, we have considered a rough upper bound for the first
factor on the right-hand side of \eqref{2.1}, namely $1$. Here we will keep these factors
and apply Theorem \ref{rdlem1}. 


Let $Q$ be a probability measure with support in the set $A$. Let $g_\ep=\frac{1}{(2\ep)^d}\1_{B_\ep(0)}$, where $B_\ep(0)$ denotes the $d$-dimensional
ball centered at $0$ and with radius $\ep$. Define
\begin{align*}
J_\ep(Q)&=\frac{1}{(2\ep)^d}\int_{I\times K} dt dx \int_A Q(dz)\, \1_{B_\ep(0)}(u(t,x)-z)\\
&=\int_{I\times K} dt dx\, (g_\ep\ast Q)(u(t,x)).
\end{align*}
As in Case 1, by applying results in \cite{en10}, we obtain $E(J_\ep(Q))\ge c_1$, for some positive constant $c_1$.

We will prove that
\beq
\label{2.10}
E\left[\left(J_\ep(Q)\right)^2\right] \le C \mathcal{E}_{\frac{\gamma-(k+1)}{\alpha}}(Q),
\eeq
for any $\alpha\in\left]0,\frac{2-\beta}{2}\right[$. 
Indeed, using Fubini's theorem and the existence of the joint density of the random vector $(u(s,y), u(t,x))$ at any points $(s,y)\neq (t,y)\in\, ]0,T]\times\IR^k$, we have
\begin{align*}
E\left[\left(J_\ep(Q)\right)^2\right]&\le \int_{A^\ep\times A^\ep} dz_1 dz_2\, (g_\ep\ast Q)(z_1)(g_\ep\ast Q)(z_2)\\
&\qquad\times \int_{I\times K} dt dx \int_{I\times K} ds dy\, p_{s,y;x,t}(z_1,z_2).
\end{align*}

Fix $\alpha\in\left]0,\frac{2-\beta}{2}\right[$. By Theorem \ref{tH}, \eqref{H} holds for $\gamma$, and $\gamma \geq k+1$ since we are in Case 2. Therefore, Theorem \ref{rdlem1} along with Theorem B.1 in \cite{dkn07} imply that
\begin{align*}
E\left[\left(J_\ep(Q)\right)^2\right]&\le C \int_{A^\ep\times A^\ep} dz_1 dz_2\, (g_\ep\ast Q)(z_1)(g_\ep\ast Q)(z_2)\\
&\qquad \times K_{\frac{\gamma-(k+1)}{\alpha}}(\vert z_1-z_2\vert)\\
&=C\mathcal{E}_{\frac{\gamma-(k+1)}{\alpha}}(g_\ep\ast Q)\\
&\le C\mathcal{E}_{\frac{\gamma-(k+1)}{\alpha}}(Q).
\end{align*}
Since $\alpha$ can be chosen arbitrarily close to $\frac{2-\beta}{2}$, this finishes the proof of part (1) of Theorem \ref{t2.10}, following the arguments of the proof of Theorem 2.1 in \cite{dss10}, including the extension to the case where $A$ is Borel but not compact.

In parts (2) and (3), we are considering hitting probabilities of stochastic processes indexed by parameters of dimension $1$ and $k$, respectively. The same arguments used in the proof of part 1 
give \eqref{t2.1111} and \eqref{t2.1112}, respectively. 

This completes the proof of Theorem \ref{t2.10}.
 \hfill\qed
\bigskip

\section{Study of the eigenvalues of the Malliavin matrix and lower bounds}
\label{s3}

In this section, we prove that the solution of the system of wave equations \eqref{1.4} does satisfy properties \eqref{2.111m}\ and \eqref{LEV}  with
$\rho=\delta+3-\beta$, where $\delta>0$ is arbitrarily small. We assume that the covariance measure $\Gamma$ satisfies (C1). We start with the contribution of the {\it small eigenvalues} of $\Gamma_Z$.

\begin{thm}
\label{p2.2m}
Assume (P1), (P2) and (C1). For any $p\in[1,\infty[$ and every small $\delta>0$, there exists a positive constant $C$ such that for all $s,t\in I$ and $x,y\in K$ with
$(s,y)\ne(t,x)$, $\vert t-s\vert \leq 1$ and $\vert x-y\vert \leq 1$, 
\beq 
\label{2.11m}
E\left[\left(\inf_{\Vert \xi\Vert=1}\xi^{\sT} \gamma_Z\xi\right)^{-p}\right]\le C \left(|t-s|+|x-y|\right)^{-p(\delta+3-\beta)}.
\eeq
\end{thm}

\noindent{\it Proof}. Assume without loss of generality that $s \leq t$. We write
\beq
\label{2.11g}
\xi^{\sT}\gamma_Z\xi = J_1+J_2,
\eeq
with $\xi=(\lambda_1,\ldots,\lambda_d,\chi_1,\ldots,\chi_d)\in \IR^{2d}$, and
\begin{align*}
J_1&=\int_0^s dr\, \Big\Vert\sum_{i=1}^d(\lambda_i-\chi_i)\left[G(s-r,y-\ast)\sigma_{i,.}(u(r,\ast))+a_{i,.}(r,\ast;s,y)\right]\\
&\qquad\qquad\qquad +\chi_i\left[G(t-r,x-\ast)\sigma_{i,.}(u(r,\ast))+a_{i,.}(r,\ast;t,x)\right]\Big\Vert^2_{\hacd},\\
J_2&=\int_s^t dr \left\Vert\sum_{i=1}^d \chi_i\left[G(t-r,x-\ast)\sigma_{i,.}(u(r,\ast))+a_{i,.}(r,\ast;t,x)\right]\right\Vert^2_{\hacd}.
\end{align*}
We shall use the notation $\lambda=(\lambda_1,\ldots,\lambda_d)$, $\chi=(\chi_1,\ldots,\chi_d)$.
\bigskip

\noindent{\bf Case 1:} $s,t\in[a,b]$, with $a\in\,]0,1[$, $1 \geq t-s>0$; $|x-y|\le t-s$.
\medskip

For any $\ep\in\,]0,s\wedge(t-s)[$, we write 
\begin{align}
J_1^\ep&=\int_{s-\ep}^s dr\, \Big\Vert\sum_{i=1}^d(\lambda_i-\chi_i)\left[G(s-r,y-\ast)\sigma_{i,.}(u(r,\ast))+a_{i,.}(r,\ast;s,y)\right]\nonumber\\
&\qquad\qquad\qquad+\chi_i\left[G(t-r,x-\ast)\sigma_{i,.}(u(r,\ast))+a_{i,.}(r,\ast;t,x)\right]\Big\Vert^2_{\hacd},\label{ss}\\
J_2^\ep&=\int_{t-\ep}^t dr \left\Vert\sum_{i=1}^d \chi_i\left[G(t-r,x-\ast)\sigma_{i,.}(u(r,\ast))+a_{i,.}(r,\ast;t,x)\right]\right\Vert^2_{\hacd}.\label{s}
\end{align}

Consider positive real numbers $\nu, \theta_1, \theta_2$  whose values will be specified later on. Then, for any 
$\ep\in\left]0,a^{\frac{1}{\theta_1}}\wedge a^{\frac{1}{\theta_2}}\wedge(t-s)^{\frac{1}{\theta_2}}\wedge\left(\frac{1}{4}\right)^{\frac{1}{\nu}}\right[$, we clearly have
\begin{align*}
\inf_{\vert \xi\vert=1}\xi^{\sT} \Gamma_Z\xi&= \min\left(\inf_{\vert \xi\vert=1, \atop 0<\vert\chi\vert<\ep^{\frac{\nu}{2}}}(J_1+J_2),
\inf_{\vert \xi\vert=1, \atop \ep^{\frac{\nu}{2}}\le\vert\chi\vert\le 1}(J_1+J_2)\right)\\
&\ge\min\left(\inf_{\vert \xi\vert=1,  \atop0<\vert\chi\vert<\ep^{\frac{\nu}{2}}}J_1^{\ep^{\theta_1}},
\inf_{\vert \xi\vert=1, \atop\ep^{\frac{\nu}{2}}\le\vert\chi\vert\le 1}J_2^{\ep^{\theta_2}}\right).
\end{align*}
We will prove the following:
\bigskip

\noindent (A)
\beqn
\inf_{\vert \xi\vert=1,  \atop0<\vert\chi\vert<\ep^{\frac{\nu}{2}}}J_1^{\ep^{\theta_1}}\ge C\ep^{\theta_1(3-\beta)}-Y_{1,\ep},
\eeqn
with $Y_{1,\ep}$ satisfying
\beq
\label{2.13m}
E\left(\left\vert Y_{1,\ep}\right\vert^p\right)\le C \ep^{\min(\theta_1(5-2\beta)^-, \theta_1+\nu)p},
\eeq
for any $p\in[1,\infty[$. Here and throughout the paper, given a real number $r>0$, we denote by $r^-$ any $\bar r\in\,]0,r[$ (usually taken close to $r$).
\bigskip

\noindent (B)
\beqn
\inf_{\vert \xi\vert=1, \atop\ep^{\frac{\nu}{2}}\le\vert\chi\vert\le 1}J_2^{\ep^{\theta_2}}\ge C\ep^{\nu+\theta_2(3-\beta)}-Y_{2,\ep},
\eeqn
with $Y_{2,\ep}$ satisfying
\beq
\label{2.14m}
E\left(\left\vert Y_{2,\ep}\right\vert^p\right)\le C \ep^{\theta_2(5-2\beta)^-p},
\eeq
for any $p\in[1,\infty[$.

We shall also prove that, for a fixed $\beta\in\, ]0,2[$, we can choose $\nu, \theta_1, \theta_2>0$ such that the following conditions are satisfied:
\begin{description}
\item{(i)} $\min(\theta_1(5-2\beta)^-,\theta_1+\nu)>\theta_1(3-\beta)$,
\item{(ii)} $\theta_2(5-2\beta)^->\nu+\theta_2(3-\beta)$.
\end{description}
Then, by setting $\ep_0:=a^{\frac{1}{\theta_1}}\wedge a^{\frac{1}{\theta_2}}\wedge(t-s)^{\frac{1}{\theta_2}}\wedge \left(\frac{1}{4}\right)^{\frac{1}{\nu}}$,
according to \cite[Proposition 3.5]{dkn09}, we will deduce that
\beq
\label{2.12m}
E\left[\left(\inf_{\Vert \xi\Vert=1}\xi^{\sT} \gamma_Z\xi\right)^{-p}\right]\le C \ep_0^{-p\alpha},
\eeq
for any $p\in[1,\infty[$, with $\alpha=\max(\theta_1(3-\beta),\nu+\theta_2(3-\beta))$, and this will lead to (\ref{2.11m}) (see \eqref{ux} and \eqref{2.30m}).
\medskip

\noindent{\it Proof of} (A).
To simplify the notation, we write $\theta$ instead of $\theta_1$. By the triangular inequality, for any norm $\Vert\cdot\Vert$,

 \begin{equation}\label{rd4.7}
   \Vert a + b \Vert ^2 \geq \frac12 \Vert a\Vert^2 - \Vert b \Vert^2,
\end{equation}
therefore
\begin{equation*}
J_1^{\ep^\theta}\ge \frac{1}{2}\int_{s-\ep^\theta}^s dr \left\Vert\sum_{i=1}^d(\lambda_i-\chi_i)G(s-r,y-\ast)\sigma_{i,.}(u(s,y))\right\Vert^2_{\hacd} -Y_{1,\ep},
\end{equation*}
with
\begin{align*}
Y_{1,\ep}&=\int_{s-\ep^\theta}^s dr\, \Big\Vert\sum_{i=1}^d(\lambda_i-\chi_i) \{G(s-r,y-\ast)\left[\sigma_{i,.}(u(s,y))-\sigma_{i,.}(u(r,\ast))\right]\\
&\qquad\qquad\qquad\qquad +a_{i,.}(r,\ast;s,y)\}\\
&\qquad\qquad\qquad+\chi_i\left[G(t-r,x-\ast)\sigma_{i,.}(u(r,\ast))+a_{i,.}(r,\ast;t,x)\right]\Big\Vert^2_{\hacd}.
\end{align*}
In this case, we are assuming that $\vert \chi\vert^2\le \ep^\nu$. Since $\vert \lambda\vert^2+\vert \chi\vert^2=1$, there exists a positive constant $c$,
depending on $\ep_0$, such that $\vert \lambda-\chi\vert\ge c$. Therefore, by Lemma \ref{la1},
\begin{align}
\label{2.15m}
&\int_{s-\ep^\theta}^s dr \left\Vert\sum_{i=1}^d(\lambda_i-\chi_i)G(s-r,y-\ast)\sigma_{i,.}(u(s,y))\right\Vert^2_{\hacd}\nonumber\\
&\qquad \ge \rho_0^2 \vert \lambda-\chi\vert^2 \int_0^{\ep^\theta} dr\, \Vert G(r,y-\ast)\Vert^2_{\hac}\nonumber\\
&\qquad\ge C \ep^{\theta(3-\beta)}.
\end{align}

Next, we establish an upper bound for the $L^p(\Omega)$ norm of $Y_{1,\ep}$, for any $p\in[1,\infty[$.
After the change of variable $s:=s-r$, we apply H\"older's inequality. Along with the Lipschitz property of $\sigma$ and \eqref{fundamental}, this yields 
\begin{align}
\label{rev1}
&E\left(\left( \int_{s-\ep^\theta}^s dr \left\Vert\sum_{i=1}^d(\lambda_i-\chi_i)G(s-r,y-\ast)\left[\sigma_{i,.}(u(s,y))-\sigma_{i,.}(u(r,\ast))\right]\right\Vert^2_{\hacd}\right)^p\right)\nonumber\\
&\quad \le C\left(\int_0^{\ep^\theta} dr\, \Vert G(r,,y-\ast)\Vert^2_\hac\right)^{p-1}\nonumber\\
&\qquad \qquad \times \int_0^{\ep^\theta} dr\int_{\IR^k\times\IR^k} G(r,dz)G(r,dw)\, |z-w|^{-\beta}
\nonumber\\
&\qquad \qquad \qquad \times E\left(|u(s,y)-u(s-r,y-z)||u(s,y)-u(s-r,y-w)|\right)^p\nonumber\\
&\quad \le C \ep^{\theta(5-2\beta)^-p}.
\end{align}
The last inequality, follows from the results of Section \ref{s1.1} and Lemma \ref{la1}. Indeed, 
under (C1), we recall from Proposition \ref{p1} that the inequality \eqref{3} holds with $\delta\in\,]0,\frac{2-\beta}{2}[$. 
\smallskip


By H\"older's inequality and since $\sigma$ is bounded, we have 
\begin{align*}
&E\left(\left\vert \int_{s-\ep^\theta}^s dr \left\Vert\sum_{i=1}^d\chi_i G(t-r,x-\ast)\sigma_{i,.}(u(r,\ast))\right\Vert^2_{\hacd}\right\vert^p\right)\\
&\qquad \le C \Vert\sigma\Vert_{\infty}^{2p}\left(\int_{s-\ep^\theta}^s dr\Vert G(t-r,x-\ast)\Vert_{\hac}^2\right)^p|\chi|^{2p}.
\end{align*}
Here  $\vert\chi\vert^2\le \ep^\nu$. Then, by Lemma \ref{la1}(c), we see that this is bounded by $C\ep^{(\nu+\theta)p}$.


Finally, by Lemma \ref{la4}  we obtain, respectively
\begin{align*}
&E\left(\left\vert \int_{s-\ep^\theta}^s dr \left\Vert\sum_{i=1}^d(\lambda_i-\chi_i)a_{i,.}(r,\ast;s,y)\right\Vert^2_{\hacd}\right\vert^p\right)\le C\ep^{\theta(3-\beta)2p},\\
&E\left(\left\vert \int_{s-\ep^\theta}^s dr \left\Vert\sum_{i=1}^d\chi_ia_{i,.}(r,\ast;t,x)\right\Vert^2_{\hacd}\right\vert^p\right)\le C\ep^{[\nu+\theta]p}.
\end{align*}
Thus, we have proved \eqref{2.13m} with $\theta_1:=\theta$. Along with \eqref{2.15m}, this ends the proof of statement (A).
\bigskip

\noindent{\it Proof of} (B).
As before, we write $\theta$ instead of $\theta_2$. Applying \eqref{rd4.7}, we obtain
\begin{equation*}
J_2^{\ep^\theta}\ge \frac{1}{2} \int_{t-\ep^\theta}^t dr \left\Vert\sum_{i=1}^d \chi_i G(t-r,x-\ast)\sigma_{i,.}(u(t,x))\right\Vert^2_{\hacd} -Y_{2,\ep},
\end{equation*}
with
\begin{align*}
Y_{2,\ep}&=\int_{t-\ep^\theta}^t dr \Big\Vert\sum_{i=1}^d \chi_i\left[G(t-r,x-\ast)[\sigma_{i,.}(u(r,\ast))-\sigma_{i,.}(u(t,x))\right]\\
& \qquad\qquad\qquad\qquad +a_{i,.}(r,\ast;t,x)]\Big\Vert^2_{\hacd}.
\end{align*}
We are now assuming $\vert \chi\vert\ge \ep^{\frac{\nu}{2}}$. Thus, by (P2) and Lemma \ref{la1},
\begin{align}
\label{2.15mm}
 &\int_{t-\ep^\theta}^t dr \left\Vert\sum_{i=1}^d \chi_i G(t-r,x-\ast)\sigma_{i,.}(u(t,x))\right\Vert^2_{\hacd}\nonumber\\
  &\qquad\qquad \ge \rho_0^2\vert \chi\vert^2\int_0^{\ep^\theta} dr\, 
\Vert G(r, x-\ast)\Vert^2_{\hac} 
   \ge C \ep^{\nu+\theta(3-\beta)}.
\end{align}

Similarly to \eqref{rev1}, we have
\begin{align*}
&E\left(\left(\int_{t-\ep^\theta}^t dr \left\Vert\sum_{i=1}^d \chi_i G(t-r,x-\ast)\left[\sigma_{i,.}(u(r,\ast))-\sigma_{i,.}(u(t,x))\right]\right\Vert^2_{\hacd}\right)^p\right)\\
&\qquad \le C \ep^{\theta(5-2\beta)^-p}.
\end{align*}
Moreover, by Lemma \ref{la4} with $s=t$,
\beqn
E\left(\left(\int_{t-\ep^\theta}^t dr \left\Vert\sum_{i=1}^d \chi_ia_{i,.}(r,\ast;t,x)\right\Vert_{\hacd}^2\right)^p\right)\le C \ep^{\theta(3-\beta)2p}.
\eeqn
Thus, we have proved \eqref{2.14m} with $\theta_2:=\theta$. Along with \eqref{2.15mm}, this ends the proof of (B).
\medskip

To end the analysis of Case 1, we must check that one can find indeed positive real numbers $\nu, \theta_1, \theta_2$ satisfying the restrictions (i) and (ii) above.

For this, we remark that for $\beta\in\, ]0,2[$, $\theta_1(5-2\beta)>\theta_1(3-\beta)$, so condition (i) is equivalent to 
$$
   \nu>\theta_1(2-\beta).
$$
As for condition (ii), it is equivalent to 
 $$\nu < \theta_2(2-\beta).$$
 Hence, we may consider $0<\theta_1<\theta_2$ and then $\nu\in\,]\theta_1(2-\beta),\theta_2(2-\beta)[$, to obtain \eqref{2.12m} with
 $\ep_0=a^{\frac{1}{\theta_1}}\wedge (t-s)^{\frac{1}{\theta_2}}$ and $\alpha=\nu+\theta_2(3-\beta)$. From here, we obtain \eqref{2.11m} by the following arguments.
 
 Clearly,
$$
   \ep_0^{-p[\nu+\theta_2(3-\beta)]}= (t-s)^{-p[\frac{\nu}{\theta_2}+(3-\beta)]}\vee a^{-p[\frac{\nu+\theta_2(3-\beta)}{\theta_1}]}.
$$
Since $t-s\le 1$, for any positive constant $C$, we have $C \leq C(t-s)^{-\frac{1}{\theta_2}}$,
and we can apply this with
$$
 C = a^{-\frac{1}{\theta_1}}.
$$
Thus, up to a positive constant $C$, the upper bound in \eqref{2.12m} may be replaced by $(t-s)^{-p[\frac{\nu}{\theta_2}+3-\beta]}$. 
 
 Finally, \eqref{2.11m} is obtained as follows. Fix $\delta\in\, ]0,2-\beta[$, then choose $0<\theta_1<\theta_2$ satisfying $\theta_1(2-\beta)<\theta_2\delta$ and $\nu \in \,]\theta_1(2-\beta),\theta_2 \delta[$. For any $p\in[1,\infty[$, we obtain
 \beq
 \label{ux}
  E\left[\left(\inf_{\vert\xi\vert=1}\xi^{\sT}\Gamma_Z\xi\right)^{-p}\right] \le C (t-s)^{-p[\delta+3-\beta]}.
 \eeq
 This clearly yields \eqref{2.11m}, since we are assuming that $|x-y|\le t-s$.
 \bigskip
 
\noindent{\bf Case 2:} $s,t\in[a,b]$, $a\in\,]0,1[$, $0\le t-s\le |x-y|$.
Let $\delta_0>0$ be a real number to be determined later on, and fix $\ep\in\,\left]0, \frac{|x-y|}{\delta_0}\wedge a\right[$. We split the analysis of this case into two subcases.
\bigskip
 
\noindent{\bf Subcase 2.1:} $0\le t-s<\ep$.
We start by writing the obvious lower bound
 \begin{align*}
 \xi^T\Gamma_Z\xi &\ge \int_{s-\ep}^s dr\, \Big\Vert\sum_{i=1}^d(\lambda_i-\chi_i)\left[G(s-r,y-\ast)\sigma_{i,.}(u(r,\ast))+a_{i,.}(r,\ast;s,y)\right]\\
&\qquad\qquad\qquad+\chi_i\left[G(t-r,x-\ast)\sigma_{i,.}(u(r,\ast))+a_{i,.}(r,\ast;t,x)\right]\Big\Vert^2_{\hacd}.
\end{align*}
By the inequality \eqref{rd4.7}, the last expression is greater than or equal to $\frac{1}{2}L_\ep(s;x,y)-Y_{3,\ep}$, where we have used the notation
\begin{align*}
 L_\ep(s;x,y):&=  \int_{s-\ep}^s dr\, \Big\Vert\sum_{i=1}^d(\lambda_i-\chi_i) G(s-r,y-\ast)\sigma_{i,.}(u(s,y))\\
 &\qquad\qquad\qquad +\chi_i G(t-r,x-\ast)\sigma_{i,.}(u(s,x))\Big\Vert^2_{\hacd},\\
 Y_{3,\ep}:&= \sum_{l=1}^4U_\ep^l (s;x,y),
 \end{align*}
 with
 \begin{align*}
 U_\ep^1(s;x,y)&= \int_{s-\ep}^s dr\\
 &\quad \times \left\Vert\sum_{i=1}^d(\lambda_i-\chi_i)G(s-r,y-\ast)\left[\sigma_{i,.}(u(s,y))-\sigma_{i,.}(u(r,\ast))\right]\right\Vert^2_{\hacd},\\
 U_\ep^2(s;x,y)&= \int_{s-\ep}^s dr \left\Vert\sum_{i=1}^d \chi_iG(t-r,x-\ast)\left[\sigma_{i,.}(u(s,x))-\sigma_{i,.}(u(r,\ast))\right]\right\Vert^2_{\hacd},\\
U_\ep^3(s;x,y)&= \int_{s-\ep}^s dr \left\Vert\sum_{i=1}^d(\lambda_i-\chi_i)a_{i,.}(r,\ast;s,y)\right\Vert^2_{\hacd},\\
U_\ep^4(s;x,y)&= \int_{s-\ep}^s dr \left\Vert\sum_{i=1}^d  \chi_ia_{i,.}(r,\ast;t,x)\right\Vert^2_{\hacd}.
 \end{align*}
 We now prove that for any $p\in[1,\infty[$,
 \beq
 \label{2.16m}
 E\left(\left\vert Y_{3,\ep}\right\vert^p\right)\le C \ep^{(5-2\beta)^-p}.
 \eeq
 Indeed, this is a consequence of $L^p(\Omega)$ estimates of each term $U_\ep^l (s;x,y)$, as follows.
 
 
 The term $E\left(\left\vert U_\ep^1(s;x,y)\right\vert^p\right)$ is the same as the first one in \eqref{rev1}, with $\theta=1$. Thus, 
 \beqn
 E\left(\left\vert U_\ep^1(s;x,y)\right\vert^p\right)\le C \ep^{(5-2\beta)^-p}.
 \eeqn
  
  With similar arguments as in \eqref{rev1}, and using the fact that $t-s < \ep$, we see that
 \begin{align*}
 & E\left(\left(\int_{s-\ep}^s dr \left\Vert\sum_{i=1}^d \chi_iG(t-r,x-\ast)\left[\sigma_{i,.}(u(s,x))-\sigma_{i,.}(u(r,\ast))\right]\right\Vert^2_{\hacd}\right)^p\right)\\
  &\quad \le C \ep^{(5-2\beta)^-p}.
  \end{align*}
  This implies that
  \beqn
 E\left(\left\vert U_\ep^2(s;x,y)\right\vert^p\right)\le C \ep^{(5-2\beta)^-p}.
 \eeqn
 From Lemma \ref{la4}, we obtain
 \beqn
 E\left(\left\vert U_\ep^3(s;x,y)\right\vert^p\right)\le C \ep^{(3-\beta)2p}.
 \eeqn
Using again that $t-s<\ep$, Lemma \ref{la4} also yields
 \begin{equation*}
 E\left(\left\vert U_\ep^4(s;x,y)\right\vert^p\right)\le C \ep^p(t-s+\ep)^{(5-2\beta)p}
 \le C \ep^{(3-\beta)2p}.
 \end{equation*}
 With this, the proof of \eqref{2.16m} is complete.
 \medskip
 
 Our next goal is to establish the lower bound
 \beq
 \label{2.17m}
 L_\ep(s;x,y)\ge c\ep^{3-\beta},
 \eeq
 for some positive constant $c$.
 For this, we apply the change of variable $r\to s-r$ and the inequality $\Vert a+b\Vert^2\ge \Vert a\Vert+\Vert b\Vert-2\vert\langle a,b\rangle\vert$, valid in any Hilbert space. Property (P2) implies that
 \begin{align*} 
 L_\ep(s;x,y)
&\ge \rho_0^2\Big[\vert \lambda-\chi\vert^2\int_0^\ep dr\, \Vert G(r,y-\ast)\Vert^2_{\hac} \\
   &\qquad\qquad + \vert\chi\vert^2\int_0^\ep dr\, \Vert
  G(t-s+r,x-\ast)\Vert^2_{\hac}\Big]\\
  &\quad-2\Big\vert\int_{s-\ep}^s dr
  \left\langle (\lambda-\chi)^{\sT}\sigma(u(s,y))G(s-r,y-\ast),\right.\\
  &\qquad\qquad\qquad\qquad\left.\chi^{\sT}\sigma(u(s,x))G(t-r,x-\ast)\right\rangle_{\hacd}\Big\vert.
\end{align*}
Set
\begin{align*}
 T_\ep(s,t;x,y)
 &= 2 \Big\vert\int_{s-\ep}^s dr
  \left\langle (\lambda-\chi)^{\sT}\sigma(u(s,y))G(s-r,y-\ast),\right.\\ 
  &\qquad\qquad\qquad\left. \mu^T\sigma(u(s,x))G(t-r,x-\ast)\right\rangle_{\hacd}\Big\vert.
\end{align*}
Remember that $\vert\lambda\vert^2+\vert\chi\vert^2=1$. Consequently, $\vert\lambda-\chi\vert^2+\vert\chi\vert^2\ge c$, for some positive contant $c$.
This fact, along with Lemma \ref{la7} below, implies that
\beqn
L_\ep(s;x,y)\ge \tilde c \ep^{3-\beta}- T_\ep(s,t;x,y).
\eeqn

   Our next aim is to prove that 
\beq
\label{2.18m}
T_\ep(s,t;x,y)\le \bar c \ep^{3-\beta}\Psi(\ep; x,y;s,t),
\eeq
where 
\beq
\label{2.19m}
\lim_{\frac{|x-y|}{\ep}\to + \infty}\Psi(\ep; x,y;s,t)=0,
\eeq
uniformly over $\frac{t-s}{\ep}\in[0,1]$.

By developing the inner product in the definition of $T_\ep(s,t;x,y)$, we easily obtain
\beqn
T_\ep(s,t;x,y)\le C \left\vert\int_0^\ep dr \left\langle G(r,y-\ast), G(t-s+r,x-\ast)\right\rangle_{\hac}\right\vert,
\eeqn
with a constant $C$ depending on $\vert\lambda\vert$, $\vert \chi\vert$, $\Vert\sigma\Vert_{\infty}$. Statements \eqref{2.18m}-\eqref{2.19m} now follow from Lemma \ref{la10}.

Property \eqref{2.19m} implies the existence of a $\delta_0>0$ such that, for any $\ep>0$
satisfying $\frac{|x-y|}{\ep}>\delta_0$, 
\beqn
\vert\Psi(\ep; x,y; s,t)\vert<\frac{\tilde c}{2\bar c}.
\eeqn
This yields
\begin{align*}
L_\ep(s;x,y)&\ge \tilde c \ep^{3-\beta}- T_\ep(s,t;x,y)\\
&\ge \frac{\tilde c}{2}\ep^{3-\beta}.
\end{align*}
Thus, for any $\ep<\frac{|x-y|}{\delta_0}$, we have proved that
\beq
\label{2.20m}
\inf_{\vert\xi\vert=1}\xi^{\sT}\gamma_Z\xi\ge C\ep^{3-\beta}-Y_{3,\ep},
\eeq
with $Y_{3,\ep}$ satisfying \eqref{2.16m}.
\bigskip

\noindent{\bf Subcase 2.2: $0\le \ep\le t-s$}. We apply the results obtained in {\bf Case 1} to conclude that
\beqn
\inf_{\vert\xi\vert=1}\xi^{\sT}\gamma_Z\xi\ge
\min\left(C\ep^{\theta_1(3-\beta)}-Y_{1,\ep}, C\ep^{\nu+\theta_2(3-\beta)}-Y_{2,\ep}\right),
\eeqn
with $Y_{1,\ep}$ and $Y_{2,\ep}$ satisfying \eqref{2.13m}, \eqref{2.14m}, respectively.
\bigskip

Combining the results proved so far, we see that
\begin{align*}
\inf_{\vert\xi\vert=1}\xi^{\sT}\gamma_Z\xi&\ge\min\left(C\ep^{\theta_1(3-\beta)}-Y_{1,\ep}1_{\{\ep\le t-s\}},\right.\\
&\left.\quad C\ep^{\max(\nu+\theta_2(3-\beta),3-\beta)}-Y_{2,\ep}1_{\{\ep\le t-s\}}-Y_{3,\ep}1_{\{\ep>t-s\}}\right).
\end{align*}
To reach the conclusion, we would like to apply Proposition 3.5 in \cite{dkn09}. For this, 
the parameters $\theta_1$, $\theta_2$ and $\nu$ must satisfy certain restrictions. Setting $\theta_2 = 1$, these restrictions become:
\begin{description}
\item{(a)} $\theta_1+\nu>\theta_1(3-\beta)$,
\item{(b)} $5-2\beta>\nu+ 3-\beta$,
\item{(c)} $5-2\beta>3-\beta$.
\end{description}
Since $\beta \in\,]0,2[$, (c) is clearly satisfied, and (a) and (b) are equivalent to
\beqn
\theta_1(2-\beta)<\nu<2-\beta. 
\eeqn
Consider $0<\theta_1<1$ and choose $\nu\in\,]\theta_1(2-\beta),2-\beta[$. According to  Proposition 3.5 in \cite{dkn09},
we obtain
\beq
\label{2.30m}
E\left[\left(\inf_{\vert\xi\vert=1}\xi^{\sT}\gamma_Z\xi\right)^{-p}\right]\le C |x-y|^{-p(\nu+3-\beta)}.
\eeq
Since $\theta_1$ can be chosen arbitrary small, this yields \eqref{2.30m} with $\nu$ replaced by any $\delta$ arbitrary close to zero.
\medskip

Finally, \eqref{2.11m} is obtained by combining the results proved in {\bf Case 1} (see \eqref{ux}) and {\bf Case 2}. This proves Theorem \ref{p2.2m}.
\hfill\qed
\bigskip

Our next efforts are devoted to the study of the {\it large eigenvalues} of $\gamma_Z$.

\begin{thm}
\label{t5n}
We assume (P1), (P2) and (C1). For any $p\in[1,\infty[$, there exists a constant $C$ depending on $p$, $I$ and $K$ such that, for any $(t,x)\ne(s,y)$ with 
$a\le s\le t\le b$ and $x,y\in K$,
\beq
\label{2.21n}
E\left[1_{A_K}\left(\Pi_{i\in K}(\xi^i)^{\sT}\gamma_Z\xi^i\right)^{-p}\right] \le C,
\eeq
with $A_K$ defined in \eqref{evl}.
\end{thm}
\noindent{\it Proof}. Fix $i_0\in\{1,\ldots,d\}$ such that $\alpha_{i_0}\ge \alpha_0$. In the sequel, we will write $\alpha$ instead of $\alpha_{i_0}$
and $\xi$ instead of $\xi_{i_0}$. By H\"older's inequality, it suffices to prove
\beq
\label{2.11n}
E\left[1_{\{\alpha\ge \alpha_0\}}\left(\xi^{\sT}\gamma_Z\xi\right)^{-p}\right] \le C,
\eeq
for any $p\in[1,\infty[$.

Fix $\theta\in\,]0,1]$, and then $\theta_1<\theta (2-\beta)/2$, so that $\theta_1 < \theta$.  In the subsequent arguments, we shall consider  
$\ep\in\,]0,a^{1/\theta_1}[$.
\smallskip

\medskip


\noindent{\bf Case 1:} $\ep^\theta< t-s$.
As in the proof of Theorem \ref{p2.2m}, we write $\xi=(\lambda,\chi)$. There exists $\alpha\in[0,1]$ and vectors $\tilde \lambda, \tilde \xi\in\mathbb{R}^d$ such that 
 $\lambda = \alpha \tilde\lambda$, $\mu = \sqrt{1-\alpha^2}\,\tilde\chi$, $\vert \tilde\lambda\vert=\vert \tilde\chi\vert =1$. For any $\eta>0$, we obviously have
\beqn
\inf_{\alpha_0\le \alpha\le 1}\left(\xi^{\sT}\gamma_Z\xi\right)=\min\left(\inf_{\alpha_0\le \alpha\le\sqrt{1-\ep^\eta}}\left(\xi^{\sT}\gamma_Z\xi\right),
\inf_{\sqrt{1-\ep^\eta}<\alpha\le1}\left(\xi^{\sT}\gamma_Z\xi\right)\right).
\eeqn
Assume $\alpha_0\le \alpha\le\sqrt{1-\ep^\eta}$. By using \eqref{2.11g}--\eqref{s}, we write
\beqn
\xi^{\sT}\gamma_Z\xi\ge J_2^{\ep^\theta}\ge\frac{1}{2} G_{1,\ep}-\bar G_{1,\ep},
\eeqn
where
\beqn
G_{1,\ep}= \int_{t-\ep^\theta}^t dr \left\Vert \sum_{i=1}^d\sqrt{1-\alpha^2}\tilde{\chi}_ iG(t-r,x-\ast)\sigma_{i,\cdot}(u(t,x))\right\Vert^2_{\hacd},
\eeqn
and
$\bar G_{1,\ep}=\bar G_{11,\ep}+\bar G_{12,\ep}$ with
\begin{align*}
\bar{G}_{11,\ep}&= \int_{t-\ep^\theta}^t dr \left\Vert \sum_{i=1}^d\sqrt{1-\alpha^2}\tilde{\chi}_ iG(t-r,x-\ast)
\left[\sigma_{i,\cdot}(u(t,\ast))-\sigma_{i,\cdot}(u(r,x))\right]\right\Vert^2_{\hacd},\\
\bar{G}_{12,\ep}&=\int_{t-\ep^\theta}^t dr \left\Vert \sum_{i=1}^d\sqrt{1-\alpha^2}\tilde\chi_ ia_{i,\cdot}(r,\ast;t,x)\right\Vert^2_{\hacd}.
\end{align*}
By Lemma \ref{la1} and hypothesis (P2), since $1-\alpha^2 \geq \ep^\eta$, we have
\beq
\label{2.23}
G_{1,\ep}\ge C\ep^{\eta+\theta(3-\beta)}.
\eeq
Similarly as in \eqref{rev1}, we have
\beqn
E\left(\left(\bar G_{11,\ep}\right)^p\right)\le C \ep^{\theta(5-2\beta)^-p}.
\eeqn
Lemma \ref{la4} with $s=t$ tells us that
\beqn
E\left(\left(\bar G_{12,\ep}\right)^p\right)\le C \ep^{\theta(3-\beta)2p}.
\eeqn
Both estimates above hold for any $p\in[1,\infty[$.
Consequently, for any $p\in[1,\infty[$,
\beq
\label{24n}
E\left(\left(\bar G_{1,\ep}\right)^p\right)\le  C \ep^{\theta(5-2\beta)^-p}.
\eeq

Let us now assume that $\sqrt{1-\ep^\eta}<\alpha\le1$. In this case, we use \eqref{2.11g}, \eqref{ss} and we consider the lower
bound
\beqn
\xi^{\sT}\gamma_Z\xi \ge J_1^{\ep^{\theta_1}}\ge \frac{1}{2}G_{2,\ep}-\bar G_{2,\ep},
\eeqn
where
\beqn
G_{2,\ep}=\int_{s-\ep^{\theta_1}}^s dr  \left\Vert \sum_{i=1}^d \alpha\tilde\lambda_iG(s-r,y-\ast)\sigma_{i,\cdot}(u(r,y))\right\Vert^2_{\hacd},
\eeqn
and $\bar G_{2,\ep}=\sum_{l=1}^4 \bar G_{2l,\ep}$, with
\begin{align*}
\bar{G}_{21,\ep}&=\int_{s-\ep^{\theta_1}}^s dr  \left\Vert \sum_{i=1}^d \alpha\tilde\lambda_iG(s-r,y-\ast)
\left[\sigma_{i,\cdot}(u(r,\ast))-\sigma_{i,\cdot}(u(r,y))\right]\right\Vert^2_{\hacd},\\
\bar{G}_{22,\ep}&=\int_{s-\ep^{\theta_1}}^s dr  \left\Vert \sum_{i=1}^d\left[\alpha\tilde\lambda_i-\sqrt{1-\alpha^2}\tilde\chi_i\right]a_{i,\cdot}(r,\ast;s,y)\right\Vert^2_{\hacd},\\
\bar{G}_{23,\ep}&=\int_{s-\ep^{\theta_1}}^s dr \Big\Vert \sum_{i=1}^d \sqrt{1-\alpha^2}\tilde\chi_i\\
&\qquad\qquad\qquad \times\left[G(t-r,x-\ast)-G(s-r,y-\ast)\right] \sigma_{i,\cdot}(u(r,\ast))\Big\Vert^2_{\hacd},\\
\bar{G}_{24,\ep}&=\int_{s-\ep^{\theta_1}}^s dr  \left\Vert \sum_{i=1}^d \sqrt{1-\alpha^2}\tilde\chi_ia_{i,\cdot}(r,\ast;t,x)\right\Vert^2_{\hacd}.
\end{align*}
We are assuming $\alpha\ge \alpha_0$. Therefore, by Lemma \ref{la1}, for some positive constant $C$,
\beq
\label{25n}
G_{2,\ep}\ge C\ep^{\theta_1(3-\beta)}.
\eeq
The support of  $G(s-r, y- \ast)$ is included in the set $\{z\in\mathbb{R}^3: |y-z|=s-r\}$. Then, similarly as in \eqref{rev1}, Lemma \ref{la1}, the Lipschitz continuity of $\sigma$ along with Proposition \ref{p1} and \eqref{3} yield
\begin{equation*}
E\left(\left(\bar G_{21,\ep}\right)^p\right)\le \ep^{\theta_1(5-2\beta)^-p},
\end{equation*}
and by Lemma \ref{la4} with $s=t$,
\begin{equation*}
E\left(\left(\bar G_{22,\ep}\right)^p\right)\le C \ep^{\theta_1(3-\beta)2p}.
\end{equation*}
The assumption on $\alpha$ implies $1-\alpha^2\le \ep^\eta$. Hence, by applying Corollary \ref{ca-1}, we obtain
\beq
E\left(\left(\bar G_{23,\ep}\right)^p\right)
\le C \ep^{[\eta+\theta_1(1-\beta)^-1_{\beta\in]0,1[}+\theta_1(2-\beta)^-1_{\beta\in[1,2[}]p}.
\eeq
By Lemma \ref{la4}, we have
\beqn
E\left(\left(\bar G_{24,\ep}\right)^p\right)\le C \ep^{\eta p}\ep^{\theta_1p}(t-s+\ep^{\theta_1})^{(5-2\beta)p}
\le C \ep^{(\eta+\theta_1)p}.
\eeqn
Summarizing the estimates obtained so far, we obtain 
\beq
\label{26n}
E\left(\left(\bar{G}_{2,\ep}\right)^p\right)\le C\ep^{\min[\theta_1(5-2\beta)^-,\eta+\theta_1(1-\beta)^-1_{\beta\in]0,1[}+\theta_1(2-\beta)^-1_{\beta\in[1,2[}]p}.
\eeq
\bigskip


\noindent{\bf Case 2:} $0\le t-s \le \ep^\theta$, $\frac{|x-y|}{\delta_0}<\ep^\theta$, with $\delta_0>0$ to be determined later (after \eqref{31r}).
According to \eqref{2.11g} and \eqref{ss}, we can write, for $\theta_2>0$ sufficiently small,
\beqn
\xi^{\sT}\gamma_Z\xi\ge J_1^{\ep^{\theta_2}}\ge \frac{1}{2}G_{3,\ep}- \bar G_{3,\ep},
\eeqn
where
\beqn
G_{3,\ep}=\int_{s-\ep^{\theta_2}}^s dr  \left\Vert \sum_{i=1}^d \alpha\tilde\lambda_iG(s-r,y-\ast)\sigma_{i,\cdot}(u(s,y))\right\Vert^2_{\hacd},
\eeqn
and $\bar G_{3,\ep}=\sum_{l=1}^4\bar {G}_{3l,\ep}$ with 
\begin{align*}
\bar{G}_{31,\ep}&=\int_{s-\ep^{\theta_2}}^s dr  \left\Vert \sum_{i=1}^d \alpha\tilde\lambda_iG(s-r,y-\ast)
\left[\sigma_{i,\cdot}(u(r,\ast))-\sigma_{i,\cdot}(u(s,y))\right]\right\Vert^2_{\hacd},\\
\bar{G}_{32,\ep}&=\int_{s-\ep^{\theta_2}}^s dr  \left\Vert \sum_{i=1}^d\left[\alpha\tilde\lambda_i-\sqrt{1-\alpha^2}\tilde\chi_i\right]a_{i,\cdot}(r,\ast;s,y)\right\Vert^2_{\hacd},\\
\bar{G}_{33,\ep}&=\int_{s-\ep^{\theta_2}}^s dr  \Big\Vert \sum_{i=1}^d \sqrt{1-\alpha^2}\tilde\chi_i\\
&\qquad\qquad\qquad \times\left[G(t-r,x-\ast)-G(s-r,y-\ast)\right]\sigma_{i,\cdot}(u(r,\ast))\Big\Vert^2_{\hacd},\\
\bar{G}_{34,\ep}&=\int_{s-\ep^{\theta_2}}^s dr  \left\Vert \sum_{i=1}^d \sqrt{1-\alpha^2}\tilde\chi_ia_{i,\cdot}(r,\ast;t,x)\right\Vert^2_{\hacd}.
\end{align*}
As has already been argued several times, by applying Lemma \ref{la1} and since $\alpha \geq \alpha_0$, we obtain
\beq
\label{27n}
G_{3,\ep} \ge C\ep^{\theta_2(3-\beta)}.
\eeq
Similarly as in \eqref{rev1}, we have
\begin{equation*}
E\left(\left(\bar G_{31,\ep}\right)^p\right)\le C \ep^{\theta_2(5-2\beta)^-p},
\end{equation*}
while thanks to Lemma \ref{la4} with $s=t$,
\begin{equation*}
E\left(\left(\bar G_{32,\ep}\right)^p\right)\le C \ep^{\theta_2(3-\beta)2p}.
\end{equation*}
These two bounds hold  for any $p\in[1,\infty[$.

From Propositions \ref{la6-1} and \ref{la6-2} and the discussion following Corollary \ref{ca-1} (see \eqref{kau}, \eqref{ep3} and \eqref{ep5}), it follows that for any $p\in[1,\infty[$, 

\begin{align}
\label{280n}
E\left(\left(\bar G_{33,\ep}\right)^p\right)\nonumber& \le C\left\{\ep^{\theta(2-\beta)^-+\theta_2}+\ep^{2\theta^-+\theta_2(1-\beta)}+\ep^{\theta^-+\theta_2(2-\beta)}\right.\nonumber\\
&\left.\qquad\quad \quad+ \ep^{\theta(1-\beta)^- + 2 \theta_2}\right\}^p 1_{\beta\in\,]0,1[}\nonumber\\
&\qquad\quad+\left\{\ep^{\theta(2-\beta)^-+\theta_2}+\ep^{\theta+\theta_2(2-\beta)}
\right.\nonumber\\
&\qquad\qquad\qquad \left.+ \ep^{\theta\frac{2-\beta}{2}+\theta_2\frac{4-\beta}{2}}
\right\}^p1_{\beta\in[1,2[}.
\end{align}
 
By applying Lemma \ref{la4}, since we are in the case where $t-s \leq \ep^\theta$, we obtain
\begin{align*}
E\left(\left(\bar G_{34,\ep}\right)^p\right)&\le C \ep^{\theta_2p}(t-s+\ep^{\theta_2})^{(5-2\beta)p}\\
&\le C \left(\ep^{[\theta_2+\theta(5-2\beta)]p}+\ep^{\theta_2(3-\beta)2p}\right)
\end{align*}
for any $p\in[1,\infty[$.

Therefore,
\beq
\label{28n}
E\left(\left(\bar{G}_{3,\ep}\right)^p\right) \le C \ep^{\min[\theta_2(5-2\beta)^-,\theta_2+\theta(5-2\beta),\pi(\theta,\theta_2,\beta)]p},
\eeq
with
\begin{align}
\label{pi}
\pi(\theta,\theta_2,\beta)&= \min\big(\theta(2-\beta)^- +\theta_2,\ 2\theta^-+\theta_2(1-\beta),\ \theta^-+\theta_2(2-\beta)\nonumber\\
&\qquad\qquad\qquad \theta (1-\beta)^- + 2 \theta_2
\big)1_{\beta\in\,]0,1[}\nonumber\\
&\qquad +  \min\Big(\theta(2-\beta)^-+\theta_2, \theta+\theta_2(2-\beta), 
\nonumber\\
&\qquad\qquad\qquad\qquad \theta\frac{2-\beta}{2}+\theta_2\frac{4-\beta}{2}
\Big)1_{\beta\in[1,2[}.
\end{align}

\medskip

\noindent{\bf Case 3:} $0 \leq t-s \leq \ep^\theta$ and  $0<\ep^\theta<\frac{|x-y|}{\delta_0}$.
Using \eqref{2.11g}--\eqref{s}, we obtain
\beqn
\xi^{\sT}\gamma_Z\xi\ge J_1^{\ep^\theta}\ge \frac{1}{2} G_{4,\ep}-\bar G_{4,\ep},
\eeqn
where
\begin{align*}
G_{4,\ep}&=\int_{s-\ep^\theta}^s dr \Big\Vert\sum_{i=1}^d \left[\left(\alpha\tilde \lambda_i - \sqrt{1-\alpha^2}\tilde\chi_i\right)
G(s-r,y-\ast)\sigma_{i,\cdot}(u(s,y))\right.\\
&\qquad\qquad\qquad\qquad \left.+\sqrt{1-\alpha^2} \tilde \chi_iG(t-r,x-\ast)\sigma_{i,\cdot}(u(s,x))\right]\Big\Vert^2_{\hacd},
\end{align*}
and $\bar G_{4,\ep}=\sum_{l=1}^4\bar G_{4l,\ep}$ with
\begin{align*}
\bar G_{41,\ep}&=\int_{s-\ep^\theta}^s dr \Big\Vert\sum_{i=1}^d \left(\alpha\tilde \lambda_i -\sqrt{1-\alpha^2}\tilde\chi_i\right)G(s-r,y-\ast)\\
&\qquad\qquad\qquad\qquad\times\left[\sigma_{i,\cdot}(u(r,\ast))-\sigma_{i,y}(u(s,y))\right]\Big\Vert^2_{\hacd},\\
\bar G_{42,\ep}&=\int_{s-\ep^\theta}^s dr \Big\Vert\sum_{i=1}^d \sqrt{1-\alpha^2} \tilde \chi_iG(t-r,x-\ast)\\
&\qquad\qquad\qquad\qquad\times\left[\sigma_{i,\cdot}(u(r,\ast))-\sigma_{i,x}(u(s,x))\right]\Big\Vert^2_{\hacd},\\
\bar G_{43,\ep}&=\int_{s-\ep^\theta}^s dr \left\Vert\sum_{i=1}^d\left(\alpha\tilde \lambda_i -\sqrt{1-\alpha^2}\tilde\chi_i\right)a_{i,\cdot}(r,\ast;s,y)\right\Vert^2_{\hacd},\\
\bar G_{44,\ep}&=\int_{s-\ep^\theta}^s dr \left\Vert\sum_{i=1}^d \sqrt{1-\alpha^2} \tilde \chi_ia_{i,\cdot}(r,\ast;t,x)\right\Vert^2_{\hacd}.
\end{align*}
We fix $p\in[1,\infty[$ and give upper bounds for the $L^p(\Omega)$-norms of each one of the terms above.

We have already seen (see \eqref{rev1}) that
\beqn
E\left(\left(\bar G_{41,\ep}\right)^p\right)\le C\ep^{\theta(5-2\beta)^-p}.
\eeqn
For the second term, we use the H\"older continuity of the sample paths, Lemma \ref{la1} and the fact that the support of $G(t-r,x-\ast)$ is included in 
$\{z\in\mathbb{R}^k: |x-z|\le t-r\}$ to see that
\begin{align*}
E\left(\left(\bar G_{42,\ep}\right)^p\right)&\le C
 (t-s+\ep^\theta)^{(2-\beta)^-p}\, 
(t-s+\ep^\theta)^{(3-\beta)p}\\
  & = C \ep^{\theta (5-2\beta)^- p},
\end{align*}
since $t-s < \ep^\theta$. For the third term, we apply Lemma \ref{la4} (with $t=s$) and we  obtain
\beqn
 E\left(\left(\bar G_{43,\ep}\right)^p\right) \leq C \ep^{\theta(3-\beta)2p}.
\eeqn
Finally, for the fourth term, we use Lemma \ref{la4}. We get 
\beqn
 E\left(\left(\bar G_{44,\ep}\right)^p\right) \leq C \ep^{\theta p}(t-s+\ep^\theta)^{(5-2\beta)p}\leq C \ep^{\theta(3-\beta)2p}.
\eeqn
From the above estimates, we conclude that
\beq
\label{30r}
E\left[\left(\sum_{l=1}^4\bar G_{4l,\ep}\right)^p\right]\le C\ep^{\theta(5-2\beta)^-p}.
\eeq
\medskip

Our next aim is to find a lower bound for $G_{4,\ep}$. For this, we apply the change of variables $r\to s-r$, and then we develop the square of the norm  on $\hacd$ to obtain
\begin{align*}
G_{4,\ep}&\ge \int_0^{\ep^\theta} dr  \left\Vert G(r,y-\ast)\right\Vert^2_{\hac}\left\vert(\alpha\tilde\lambda^T-\sqrt{1-\alpha^2}\tilde \chi^T)\sigma(u(s,y))\right\vert^2\\
&\qquad+\int_0^{\ep^\theta} dr  \left\Vert G(t-s+r,x-\ast)\right\Vert^2_{\hac}\left\vert\sqrt{1-\alpha^2}\tilde \chi^T\sigma(u(s,x))\right\vert^2\\
&\qquad-2\Big\vert \int_0^{\ep^\theta} dr \left\langle G(r,y-\ast),G(t-s+r,x-\ast)\right\rangle_\hac\\
& \qquad\qquad\qquad\times\left\langle\left(\alpha\tilde\lambda-\sqrt{1-\alpha^2}\tilde\chi\right)
\sigma(u(s,y)), \sqrt{1-\alpha^2}\tilde\chi\sigma(u(s,x))\right\rangle\Big\vert.
\end{align*}
Property (P2) along with Lemmas \ref{la1} and \ref{la7} yield
\beqn
G_{4,\ep}\ge c\rho_0^2\varphi(\alpha,\tilde\lambda,\tilde\chi)\ep^{\theta(3-\beta)}-2\left\vert\psi(\alpha,\tilde\lambda,\tilde\chi)I(\ep;t,x,s,y)\right\vert,
\eeqn
with
\begin{align*}
\varphi(\alpha,\tilde\lambda,\tilde\mu)&=\left\vert\alpha\tilde\lambda^T-\sqrt{1-\alpha^2}\tilde\chi^T\right\vert^2+\left\vert\sqrt{1-\alpha^2}\tilde\chi^T\right\vert^2,\\
\psi(\alpha,\tilde\lambda,\tilde\chi)&=\left\langle\left(\alpha\tilde\lambda-\sqrt{1-\alpha^2}\tilde\chi\right)
\sigma(u(s,y)), \sqrt{1-\alpha^2}\tilde\chi  \sigma(u(s,x))\right\rangle,\\
I(\ep;t,x,s,y)
&=\int_0^{\ep^\theta} dr \left\langle G(r,y-\ast),G(t-s+r,x-\ast)\right\rangle_\hac,
\end{align*}
and $c=\frac{1}{3-\beta}\int_{\rk}\frac{\sin^2(|w|)}{|w|^{k-\beta+2}}dw$.

One can easily check that
\beqn
\varphi_0:=\inf_{\alpha\in[\alpha_0,1]\atop \Vert\tilde\lambda\Vert=\Vert\tilde\chi\Vert=1}\varphi(\alpha,\tilde\lambda,\tilde\chi)>0.
\eeqn
Also,
\beqn
\sup_{\alpha\in[\alpha_0,1]\atop \Vert\tilde\lambda\Vert=\Vert\tilde\chi\Vert=1}\vert \psi(\alpha,\tilde\lambda,\tilde\chi)\vert\le 2\Vert\sigma\Vert_\infty^2.
\eeqn
Therefore,
\beq
\label{31r}
G_{4,\ep}\ge c\rho_0^2\varphi_0 \ep^{\theta(3-\beta)}
-2 \Vert\sigma\Vert_\infty^2 \left\vert I(\ep;t,x,s,y)\right\vert.
\eeq
Since $0\leq t-s \leq \ep^\theta$, by Lemma \ref{la10}, there exists $\delta_0>0$ such that for any $\ep>0$ satisfying
$\frac{|x-y|}{\ep^\theta}>\delta_0$, we have that
\beqn
\Vert\sigma\Vert_\infty^2 \left\vert I(\ep;t,x,s,y)\right\vert\le \frac{c\rho_0^2\varphi_0}{4}\, \ep^{\theta(3-\beta)}.
\eeqn
From \eqref{31r}, it follows that
\beq
\label{32r}
G_{4,\ep}\ge\frac{c\rho_0^2 \varphi_0}{2}\ep^{\theta(3-\beta)}.
\eeq


With the estimates \eqref{2.23}, \eqref{25n}, \eqref{27n} and \eqref{32r}, we have established the following. For any $0<\ep<a^{\frac{1}{\theta_1}}$,
\begin{align*}
\inf_{\alpha_0\le\alpha\le 1}\left(\xi^T \Gamma_Z\xi\right)&\ge \min\left(c\ep^{\eta+\theta(3-\beta)} -
\bar G_{1,\ep}, c\ep^{\theta_1(3-\beta)}-\bar G_{2,\ep},\right.\\
&\left.\quad \quad c\ep^{\theta_2(3-\beta)}- \bar G_{3,\ep}, c\ep^{\theta(3-\beta)}-\bar G_{4,\ep}\right)\\
&\ge \min\left(c\ep^{\eta+\theta(3-\beta)} - (\bar G_{1,\ep}+\bar G_{4,\ep}),\right.\\
&\left.\qquad \qquad c\ep^{(\theta_1\vee\theta_2)(3-\beta)} - (\bar G_{2,\ep}+\bar G_{3,\ep})\right).
\end{align*}
Choose $\theta_2=\theta_1$. We will apply \cite[Proposition 3.5]{dkn09}. For this, taking into account the inequalities \eqref{24n}, \eqref{30r}, \eqref{26n}, \eqref{28n}, the following conditions must be satisfied:
\begin{description}
\item{(a)} $\eta+\theta(3-\beta)<\theta(5-2\beta)$,
\item{(b)} $\theta_1(3-\beta)<\min(\theta_1(5-2\beta), \eta+\theta_1(1-\beta)1_{\beta\in]0,1[}+\theta_1(2-\beta)1_{\beta\in[1,2[})$,
\item{(c)} $\theta_1(3-\beta)<\min(\theta_1(5-2\beta), \theta_1+\theta(5-2\beta), \pi(\theta,\theta_1,\beta)) $,
\end{description}
with $\pi(\theta,\theta_1,\beta)$ defined in \eqref{pi}. It is easy to check that these constraints hold by choosing $\theta$, $\theta_1$ and $\eta$
such that
\beqn
   2\theta_1<\eta<\theta(2-\beta), \qquad 0< \theta_1<\theta.
\eeqn
Hence Theorem \ref{t5n} is proved. 
\hfill\qed
\bigskip

We now apply the results just obtained, along with those obtained in Section \ref{s2}, to establish lower bounds for the hitting probabilities of the solution of the system of equations given in \eqref{1.4}. 

\begin{thm}
\label{t2}
Let $k\in\{1,2,3\}$. Assume (P1), (P2) and (C1). Let $I=[a,b]\subset [0,T]$, with $a>0$, let $K$ be a compact subset of $\rk$ and fix $N>0$.
\begin{description}
\item{(1)} Fix $\delta>0$. There exists a positive constant $c=c(I,K,N,\beta,k,d,\delta)$ such that, for any
Borel set $A\subset[-N,N]^d$,
\beqn
P\left\{u(I\times K)\cap A\ne\emptyset\right\}\ge c\ {\rm Cap}_{d\left(1+\frac{4d}{2-\beta}\right)+\delta-\frac{2(k+1)}{2-\beta}}(A).
\eeqn

\item{(2)} Fix $\delta>0$ and $x\in K$. There exists a positive constant $c=c(I,x,N,\beta,k,d,\delta)$ such that, for any
Borel set $A\subset[-N,N]^d$,
\beqn
P\left\{u(I\times \{x\})\cap A\ne\emptyset\right\}\ge c\ {\rm Cap}_{d\left(1+\frac{4d}{2-\beta}\right)+\delta-\frac{2}{2-\beta}}(A).
\eeqn
\item{(3)} Fix $\delta>0$ and $t\in I$. There exists a positive constant $c=c(t,K,N,\beta,k,d)$ such that, for any
Borel set $A\subset[-N,N]^d$,
\beqn
P\left\{u(\{t\}\times K)\cap A\ne\emptyset\right\}\ge c\ {\rm Cap}_{d\left(1+\frac{4d}{2-\beta}\right)+\delta-\frac{2k}{2-\beta}}(A).
\eeqn
\end{description}
\end{thm}
\medskip


\noindent{\bf Proof of Theorem \ref{t2}}.
Theorems \ref{p2.2m} and \ref{t5n} show that assumptions \eqref{2.111m} and \eqref{LEV} of Theorem \ref{t2.10} hold with $\rho=\delta+3-\beta$, where $\delta$ is an arbitrarily small positive real number. \hfill\qed



\section{Appendix}
\label{a}

In this section, we gather important technical results that have been used throughout the preceding sections. Except otherwise stated, we assume (C1) and (P1).
\begin{lem}
\label{la1} 
For $k\in\mathbb{N}^\ast$ and $\beta \in\,]0,2\wedge k[$, we have the following.
\begin{description}
\item{(a)} For any $r>0$ and $x \in \IR^k$,
\beq
\label{basic}
   \Vert G(r,x-\ast)\Vert^2_{\hac} = c r^{2-\beta},
\eeq
with $c=\int_{\rk} \frac{\sin^2(|w|)}{|w|^{k-\beta+2}}\ dw$.

 \item {(b)} For any $T > 0$ and $r_0\ge0$, there is $C >0$ such that
\begin{align*}
\sup_{(t,x)\in [0,T]\times\rk}\left[\int_0^{t\wedge r_0} dr \Vert G(t-r,x-\ast)\Vert^2_{\hac}\right]&\le c\sup_{t\in[0,T]}\int_0^{t\wedge r_0}(t-r)^{2-\beta} dr\\
& \le C r_0,\\
\sup_{x\in\rk}\left[\int_0^{r_0} dr \left\Vert G(r,x-\ast)\right\Vert_{\hac}^2\right]&\le \frac{c}{3-\beta}r_0^{3-\beta}.
\end{align*}

 \item{(c)} For any $0<s<t\le T$ and each $\ep\in\, ]0,s[$, 
\beqn
\sup_{x\in\rk} \left[\int_{s-\ep}^s dr\ \Vert  G(t-r,x-\ast)\Vert^2_{\hac}\right] \le C\left((t-s+\ep)^{3-\beta}-(t-s)^{3-\beta}\right).
\eeqn
\end{description}
\end{lem}

\noindent{\it Proof.} By its very definition (see \eqref{fg}),
\beqn
\Vert G(r,x-\ast)\Vert_{\hac}^2=\int_{\rk} \frac{\sin^2(r|\xi|)}{|\xi|^{k-\beta+2}} d\xi.
\eeqn
With the change of variable $\xi\to r\xi$, we obtain \eqref{basic}. The remaining statements follow easily.
\hfill\qed


\begin{lem}
\label{la4}
Let $k \in \{1,2,3\}$ and let $\{a_{i,j}(r,z;t,x)$, $0\le r\le t$, $x,y\in\IR^k$, $1\le i,j\le d\}$ be the stochastic process defined in \eqref{as}.
Fix $0<s\leq t\le T$. Then, for any $\ep\in\,]0,s[$, $\chi\in\mathbb{R}^d$, 
and each $p\in[1,\infty[$,
\beqn
\sup_{x\in\rk}E\left(\left\vert \int_{s-\ep}^s dr \left\Vert\sum_{i=1}^d\chi_ia_{i,.}(r,\ast;t,x)\right\Vert^2_{\hacd}\right\vert^p\right)\le
 C\vert \chi\vert^{2p} \ep^p(t-s+\ep)^{(5-2\beta)p}.
 \eeqn
 In particular, when $s=t$,
 \beqn
\sup_{x\in\rk}E\left(\left\vert \int_{s-\ep}^s dr \left\Vert\sum_{i=1}^d\chi_ia_{i,.}(r,\ast;s,x)\right\Vert^2_{\hacd}\right\vert^p\right)\le
 C\vert \chi\vert^{2p} \ep^{(3-\beta)2p}.
 \eeqn
\end{lem}
\noindent{\it Proof}. Set
\begin{align*}
T_{1}^{i,j}(r,\ast,t,x)&= \sum_{\ell=1}^d \int_0^t \int_{\rk}G(t-\rho,x-z)\, \nabla\sigma_{i,\ell}(u(\rho,z)) \cdot D_{r,\ast}^{(j)}u(\rho,z)\\  
&\qquad\qquad\qquad \times  M^{\ell}(d\rho,dz),\\
T_{2}^{i,j}(r,\ast,t,x)&= \int_0^t d\rho \int_{\rk} dy\, G(t-\rho,x-z)\, \nabla b_{i}(u(\rho,z)) \cdot D_{r,\ast}^{(\cdot)}u(\rho,z). 
\end{align*}
We will prove that, for any $p\in[1,\infty[$,
\begin{align}
&\sup_{(t,x)\in[0,T]\times\rk}E\left(\left\vert \int_{s-\ep}^s dr \left\Vert\sum_{i=1}^d \chi_i T_{1}^{i,\cdot}(r,\ast,t,x)\right\Vert_{\hacd}^2
\right\vert^p\right)\nonumber\\
&\qquad\qquad\qquad\qquad\qquad\qquad \le C\vert \chi\vert^{2p} \ep^p(t-s+\ep)^{(5-2\beta)p},\label{a20}\\ \nonumber
\\
&\sup_{(t,x)\in[0,T]\times\rk}E\left(\left\vert \int_{s-\ep}^s dr \left\Vert\sum_{i=1}^d \chi_i T_{2}^{i,\cdot}(r,\ast,t,x)\right\Vert_{\hacd}^2
\right\vert^p\right)\nonumber\\
&\qquad\qquad\qquad\qquad\qquad\qquad \le C\vert \chi\vert^{2p} \ep^p(t-s+\ep)^{(6-\beta)p}.
\label{a21}
\end{align}
Appealing to \eqref{as}, the conclusion of the lemma will follow.

In order to prove \eqref{a20}, we first write
\begin{align}\label{rd6.3a}
&E\left(\left\vert \int_{s-\ep}^s dr \left\Vert\sum_{i=1}^d \chi_i T_{1,\ep}^{i,\cdot}(r,\ast,t,x)\right\Vert_{\hacd}^2\right\vert^p\right)\\ \nonumber
 &\qquad = E\left(\Big\Vert\sum_{\ell=1}^d \int_{0}^t \int_{\rk} G(t-\rho,x-z)\right.\\ \nonumber
&\left.\qquad\qquad\qquad \times  \nabla([\chi^{\sT}\sigma(u(\rho,z))]_\ell) \cdot D_{\cdot,\ast} u(\rho,z)\, M^\ell(d\rho,dz) \Big\Vert_{L^2([s-\ep, s]; \hac^d)}^{2p}\right)\\ \nonumber
&\qquad \leq E\left(\sum_{\ell=1}^d \Big\Vert \int_{s-\ep}^t \int_{\rk} G(t-\rho,x-z)\right.\\ \nonumber
&\left.\qquad\qquad\qquad \times  \nabla([\chi^{\sT}\sigma(u(\rho,z))]_\ell) \cdot D_{\cdot,\ast} u(\rho,z)\, M^\ell(d\rho,dz) \Big\Vert_{L^2([s-\ep, s]; \hac^d)}^{2p}\right)
\end{align}
and then, we apply Theorem 6.1 in \cite{ss} to the $L^2([s-\ep, s]; \hac)$--valued stochastic processes
\beqn
\left(\nabla\left(\left[\chi^{\sT} \sigma(u(\rho,z))\right]_\ell\right) \cdot D_{\cdot,\ast}^{(\cdot)}u(\rho,z),\quad (\rho,z)\in[s-\ep,t]\times\rk\right), 
\eeqn
$\ell=1,\ldots,d$. This shows that the expression in \eqref{rd6.3a} is bounded above by
\begin{align}
& C\sum_{\ell=1}^d\left(\int_{s-\ep}^t d\rho \int_{\rk}\mu(d\xi) \left\vert\tf G(t-\rho)(\xi)\right\vert^2\right)^{p-1}\nonumber\\  \nonumber
&\qquad \times \int_{s-\ep}^t d\rho\, \sup_{z\in\rk}E\left(\left\Vert \nabla\left(\left[\chi^{\sT}\sigma(u(\rho,z))\right]_\ell\right) \cdot D_{r,\ast} u(\rho,z)\right\Vert^{2p}_{L^2([s-\ep, s]; \hac^d)}\right)\\    
&\qquad\qquad\qquad \times \int_{\rk}\mu(d\xi)\left\vert\tf G(t-\rho)(\xi)\right\vert^2.\label{a22}
\end{align}
Following the arguments of the proof of Lemma 8.2 in \cite{ss}, we can prove 
\beq
\label{derivative}
E\left(\Vert Du(\rho,z)\Vert^{2p}_{L^2([s-\ep, s]; \hacd)}\right)
 \le C \left(\int_{s-\ep}^{s\wedge \rho} dr\, \Vert G(\rho-r,\ast)\Vert_{\hac}^2\right)^p,
\eeq
for $\rho\in[s-\ep,t]$.
Moreover, \eqref{basic} yields
\begin{align*}
 \int_{s-\ep}^{s\wedge \rho} dr\, \Vert G(\rho-r,\ast)\Vert_{\hac}^2 &\le \int_{(\rho-s)\vee 0}^{\rho-s+\ep} dr\, \Vert G(r,\ast)\Vert_{\hac}^2 = c \int_{(\rho-s)\vee 0}^{\rho-s+\ep} dr\, r^{2 - \beta}\\
&\le C \ep(\rho-s+\ep)^{2-\beta}.
\end{align*}

Along with (P1) this yields,
\begin{align*}
& E\left(\left\Vert \nabla \left(\left[\chi^{\sT}\sigma(u(\rho,z))\right]_\ell\right) \cdot D_{\cdot,\ast} u(\rho,z)\right\Vert^{2p}_{L^2([s-\ep, s]; \hac^d)}\right)\nonumber\\
&\qquad \le C \vert\chi\vert^{2p}\sup_{(\rho,z)\in[s-\ep,t]\times\rk}E\left(\Vert Du(\rho,z)\Vert^{2p}_{L^2([s-\ep, s]; \hacd)}\right)\nonumber\\
&\qquad \le C \vert\chi\vert^{2p}\ep^p\sup_{\rho\in[s-\ep,t]} (\rho - s +\ep)^{(2-\beta)p}\nonumber\\
&\qquad \le C \vert\chi\vert^{2p} \ep^p(t-s+\ep)^{(2-\beta)p}.
\end{align*}
Using this estimate in \eqref{a22} gives 
\begin{align*}
&E\left(\left( \int_{s-\ep}^t dr \left\Vert\sum_{i=1}^d \chi_i T_{1,\ep}^{i,\cdot}(r,\ast,t,x)\right\Vert_{\hacd}^2\right)^p\right)\\
&\quad\le C \vert\chi\vert^{2p}\ep^{p} (t-s+\ep)^{(2-\beta)p}\left(\int_{s-\ep}^t d\rho\int_{\rk}\mu(d\xi)\left\vert\tf G(t-\rho)(\xi)\right\vert^2\right)^p\\
&\quad\le C \vert\chi\vert^{2p} \ep^p (t-s+\ep)^{(5-2\beta)p},
\end{align*}
where in the last inequality, we have applied Lemma \ref{la1}. 
\medskip

We now prove \eqref{a21}. For this, we write
\begin{align*}
&E\left(\left\vert \int_{s-\ep}^s dr \left\Vert\sum_{i=1}^d \chi_i T_{2}^{i,\cdot}(r,\ast,t,x)\right\Vert_{\hacd}^2\right\vert^p\right)\\
&\quad = E\left(\sum_{j=1}^d\left\Vert\int_{s-\ep}^t d\rho \int_{\rk}  G(t-\rho,dz)\right.\right.\\
&\left.\left. \qquad\qquad\qquad \times \nabla \left(\chi^{\sT} b(u(\rho,x-z))\right) \cdot D u(\rho,x-z)\right\Vert_{L^2([s-\ep, s]; \hac^d)}^{2p}\right).
\end{align*}
By Minkowski's inequality, this last expression is bounded by 
\begin{align*}
  &E\Big[\Big(\int_{s-\ep}^t d\rho \int_{\rk} dz\, G(t-\rho,x-z)\left\vert\nabla (\chi^{\sT}  b(u(\rho,z)))\right\vert\\
  &\qquad \qquad\qquad\qquad\times \Vert D u(\rho,z)\Vert_{L^2([s-\ep, s]; \hacd)}\Big)^{2p}\Big].
\end{align*}
We bound this as in \eqref{a22}, and then we use \eqref{derivative} and assumption (P1) to find that
\begin{align*}
   &E\left(\left\vert \int_{s-\ep}^s dr \left\Vert\sum_{i=1}^d \chi_i T_{2}^{i,\cdot}(r,\ast,t,x)\right\Vert_{\hacd}^2\right\vert^p\right)\\
   & \qquad \leq  C\vert\chi\vert^{2p} \ep^p(t-s+\ep)^{(2-\beta)p}\left(\int_{s-\ep}^t d\rho\int_{\rk} G(t-\rho,dz)\right)^{2p}.
\end{align*}
Since 
\beqn
\int_{s-\ep}^t d\rho\int_{\rk} G(t-\rho,dz) = \int_{s-\ep}^t d\rho\, (t-\rho) \le C (t-s+\ep)^2,
\eeqn
we obtain \eqref{a21}.
\hfill\qed
\vskip 16pt


   In the sequel, we will use several times the property
\begin{equation}\label{a22a}
 \Vert G(r,x-\ast) Z(t,\ast) \Vert_{\hac}^2 = \int_{\IR^k\times\IR^k} G(r,du) G(r,dv) Z(t,u) \Vert u-v\Vert^{-\beta} Z(t,v),
\end{equation}
for an $L^2$-continuous process $(Z(t,x))$ such that
$$
   \sup_{(t,x) \in [0,T]\times \IR^k} E(\vert Z(t,x)\vert^2) < \infty.
$$
This property is a consequence of \eqref{fundamental} applied to the signed measure $\varphi(dx) = Z(t,x) G(r,dx)$, which a.s. has finite total variation (and even compact support).

\begin{prop}
\label{la6-1} Let $k \in \{1,2,3\}$, and let $Z=\{Z(t,x),\ (t,x)\in[0,T]\times\rk\}$ be a stochastic process satisfying the following conditions:
\begin{description}
\item{(a)}
 there exists $\rho\in\,]0,1]$ such that,  for every $p\in[1,\infty[$, there is $C < \infty$ such that for any $x,y\in\rk$,
\begin{equation*}
\sup_{t\in[0,T]}E\left(\left\vert Z(t,x)-Z(t,y)\right\vert^p\right)\le C\vert x-y\vert^{\rho p};
\end{equation*}
\item{(b)} for every $p\in[1,\infty[$, there is $C < \infty$ such that
\begin{equation*}
\sup_{(t,x)\in[0,T]\times \rk}E\left(\left\vert Z(t,x)\right\vert^p\right)\le C.
\end{equation*}
\end{description}
Then for any $p\in[1,\infty[$, there exists a constant $C>0$ such that for any $\ep\in\,]0,1]$, $s\in[\ep,T]$ and $x,y \in \rk$,
\begin{align*}
&E\left(\int_0^\ep dr \left\Vert [G(r,x-\ast)-G(r,y-\ast)] Z(s-r,\ast)\right\Vert^2_{\hac}\right)^p\\
&\quad \le C \left\{|x-y|^{2\rho p}\ep^{(3-\beta)p}+|x-y|^{(\alpha+\rho)p}\ep^{[3-(\alpha+\beta)]p}+|x-y|^{\bar\alpha p}
\ep^{[3-(\bar\alpha+\beta]p}\right\},
\end{align*}
with $\alpha\in\,]0,((2\wedge k)-\beta)\wedge1[$, $\bar\alpha\in\,]0,(2\wedge k)-\beta[$.
\end{prop}

\noindent{\it Proof}. Following the ideas in the first part of the proof of Proposition 3.5 in \cite{dss09} and using \eqref{a22a}, we write
\begin{align*}
&E\left[\left(\int_0^\ep dr \left\Vert [G(r,x-\ast)-G(r,y-\ast)] Z(s-r,\ast)\right\Vert^2_{\hac}\right)^p\right]\\
 &\quad\quad   \le C \sum_{i=1}^4 E\left(\left\vert J^\ep_i(x,y)\right\vert^p\right),
\end{align*}
where, for $i=1,\ldots,4$, 
\beqn
J^\ep_i(x,y)=\int_0^\ep dr \int_{\rk}\int_{\rk}G(r,du) G(r,dv)h_i(\ep,r,x,y,u,v)
\eeqn
with
\begin{align*}
h_1(\ep,r,x,y,u,v)&= f(y-x+v-u)(Z(s-r,x-u)-Z(s-r,y-u))\\
& \quad \times(Z(s-r,x-v)-Z(s-r,y-v)),\\
h_2(\ep,r,x,y,u,v)&= Df(v-u,x-y) Z(s-r,x-u)\\
& \quad \times(Z(s-r,x-v)-Z(s-r,y-v)),\\
h_3(\ep,r,x,y,u,v)&= Df(v-u,y-x) Z(s-r,y-v)\\
& \quad \times(Z(s-r,x-u)-Z(s-r,y-u)),\\
h_4(\ep,r,x,y,u,v)&= -D^2f(v-u,x-y)Z(s-r,y-u)Z(s-r,x-v).
\end{align*}
In these expressions, $f(x)=\vert x\vert^{-\beta}$, $\beta\in\,]0,2\wedge k[$ and
\begin{align*}
Df(u,x)& = f(u+x)-f(u),\\
D^2f(u,x)& = f(u-x)-2f(u)+f(u+x).
\end{align*}
For $k \in \{1,2,3\}$, $G(r,du) \geq 0$. Then, from  H\"older's inequality, the properties of $Z$, \eqref{a22a} and \eqref{basic},
\begin{align}
E\left(\left\vert J^\ep_1(x,y)\right\vert^p\right)&\le \sup_{t\in[0,T]} \left\{E\left(\left\vert Z(t,x)-Z(t,y)\right\vert^{2p}\right)\right\}\nonumber\\
&\quad\quad\times\left(\int_0^\ep dr \int_{\rk\times\rk}G(r,du) G(r,dv)\, f(y-x+v-u)\right)^p\nonumber\\\
&\le C|x-y|^{2\rho p}\left(\int_0^\ep dr\, \Vert G(r,\ast)\Vert^2_{\hac}\right)^p\nonumber\\\
&\le C|x-y|^{2\rho p}\ep^{(3-\beta)p}. \label{i1}
\end{align}

Set 
\beq
\label{mu2}
\mu_2(x,y)= \int_0^\ep dr \int_{\rk\times\rk} G(r,du) G(r,dv)\left\vert Df(v-u,x-y)\right\vert
\eeq
The variant of Lemma 6.1 in \cite{dss09} given in Lemma \ref{l6.1-bis} below shows that
\beqn
\mu_2(x,y)\le C|x-y|^\alpha \ep^{3-(\alpha+\beta)},
\eeqn
for any $\alpha\in\,]0,((2\wedge k)-\beta)\wedge 1[$. Thus, H\"older's inequality along with the assumptions on the process $Z$
yield
\begin{align}
E\left(\left\vert J^\ep_2(x,y)\right\vert^p\right)&\le C \left(\mu_2(x,y)\right)^p |x-y|^{\rho p}\nonumber\\
&\le C |x-y|^{(\alpha+\rho)p}\ep^{[3-(\alpha+\beta)]p}. \label{i2}
\end{align}
The same upper bound holds for $E\left(\left\vert J^\ep_3(x,y)\right\vert^p\right)$. 

For the analysis of the remaining term, we consider
\beq
\label{mu4}
\mu_4(x,y)= \int_0^\ep dr \int_{\rk\times\rk} G(r,du) G(r,dv)\left\vert D^2f(v-u,x-y)\right\vert.
\eeq
By  applying  the modified version of Lemma 6.2 in \cite{dss09} given in
Lemma \ref{l6.2-bis} below, we obtain
\begin{align}
E\left(\left\vert J^\ep_4(x,y)\right\vert^p\right)&\le C \left(\mu^4(x,y)\right)^{p}
\sup_{(t,x)\in[0,T]\times \rk}E\left(\left\vert Z(t,x)\right\vert^{2p}\right)\nonumber\\
&\le |x-y|^{\bar\alpha p}\ep^{[3-(\bar\alpha+\beta)]p},\label{i3}
\end{align}
for $\bar\alpha\in\,]0,(2\wedge k)-\beta[$.

With \eqref{i1}, \eqref{i2} and \eqref{i3} we finish the proof of Proposition \ref{la6-1}.  \hfill \qed
\bigskip

The two lemmas below were used in the proof of Proposition \ref{la6-1}. 

\begin{lem}
\label{l6.1-bis}
Let $\mu_2(x,y)$ be defined by \eqref{mu2}. For any $\alpha\in\,]0,((2\wedge k)-\beta)\wedge 1[$, there is $C>0$ such that for all $x,y\in \rk$,
\beqn
\mu_2(x,y) \le C |x-y|^{\alpha}\ep^{3-(\alpha+\beta)}.
\eeqn
\end{lem}
\noindent{\it Proof}. In comparison with \cite[Lemma 6.1]{dss09}, this result quantifies the dependence of $\mu_2(x,y)$ on the domain of integration in time, it is about the fundamental solution $G$ of the wave equation instead of the regularisation $G_n$,
and the function $f$ does not include a smooth factor $\varphi$. The proofs of both lemmas are very similar. We give the main lines.

Appealing to Lemma 2.6(a) of \cite{dss09} with $d:=k$, $b:= \alpha$, $a:=k-(\alpha+\beta)$, $u:=v-u$, $c:=|x-y|$, $x:=\frac{x-y}{|x-y|}$, we obtain, after using \eqref{a22a},
\beqn
\mu_2(x,y) \le |x-y|^\alpha\int_0^\ep dr \int_{\rk} d\xi\, \vert\tf G(r)(\xi)\vert^2\, \vert\xi\vert^{-(k-(\alpha+\beta))}.
\eeqn
For $\alpha+\beta\in\,]0,2\wedge k[$,  and up to a positive constant, the last integral 
is bounded by $\ep^{3-(\alpha+\beta)}$ (see Lemma \ref{la1}).
This finishes the proof. \hfill\qed 

\begin{lem}
\label{l6.2-bis}
Let $\mu_4(x,y)$ be defined by \eqref{mu4}. For each $\alpha\in\,]0,(2\wedge k)-\beta[$, there is a positive constant $C$ such that for all $x,y\in \rk$
\beqn
\mu_4(x,y)\le C |x-y|^\alpha \ep^{3-(\alpha+\beta)}.
\eeqn
\end{lem}
\noindent{\it Proof}. The same differences mentioned in the proof of Lemma \ref{l6.1-bis} apply here when comparing this statement with Lemma 6.2 in \cite{dss09}.

By applying Lemma 2.6(c) of \cite{dss09} with $d:=k$, $b:=\alpha$, $a:=k-(\alpha+\beta)$, $\alpha+\beta\in\,]0,2\wedge k[$, $u:=v-u$, $x:=x-y$ yield
\begin{align*}
\mu_4(x,y)&\le |x-y|^{\alpha}\\
&\quad\times \left (\sup_{x,y,w} \int_0^\ep dr\int_{\rk\times\rk} G(r,du)G(r,dv)\left\vert\  \vert y-x\vert w+u-v\right\vert^{-(\alpha+\beta)}\right)\\
&\quad \times \int_{\rk} dw \left\vert D^2\left(\vert \cdot\vert^{-(k-\alpha)}\right)(w,e)\right\vert,
\end{align*}
where $e$ denotes a unitary vector of $\rk$.

By Lemma \ref{la1}, Lemma 2.6(d) of \cite{dss09}, \eqref{a22a} and \eqref{basic}, the last expression is bounded by $|x-y|^{\alpha}\ep^{3-(\alpha+\beta)}$. \hfill\qed 

\begin{prop}
\label{la6-2}
Let $k \in \{1,2,3\}$ and let $Z=\{Z(t,x),\ (t,x)\in[0,T]\times \rk\}$ be a stochastic process satisfying the assumptions of Proposition \ref{la6-1}. 
Then for any $p\in[1,\infty[$, there exists a constant $C>0$ such that for all $0\le s<t\le T$ and $\ep\in\,]0,s[$, 
\begin{align}
& E\left(\int_{s-\ep}^s dr \left\Vert \left[G(t-r,x-\ast)-G(s-r,x-\ast)\right] Z(r,\ast)\right\Vert_\hac^2\right)^p\nonumber\\
&\quad \le C \Big\{
\ep^{[\frac{4-\beta}{2}]p}|t-s|^{2\rho p}|t-s+\ep|^{\frac{(2-\beta)}{2}p}\nonumber\\
&\left.\qquad\qquad+|t-s|^{\rho p}\Big[\ep|t-s||t-s+\ep|^{1-\beta}1_{\beta\in\,]0,1[}+\ep|t-s|^{2-\beta}1_{\beta\in[1,2[}\right.\nonumber\\
&\qquad\qquad+ \ep|t-s|^{\alpha_1}|t-s+\ep|^{2-(\alpha_1+\beta)}+\ep^{2-\frac{\alpha_1+\beta}{2}}|t-s|^{\alpha_1}|t-s+\ep|^{1-\frac{\alpha_1+\beta}{2}}\Big]^p\nonumber\\
&\qquad\qquad+ 1_{\beta\in[1,2[}\Big[\ep|t-s|^{2-\beta}+\ep |t-s|^{\alpha_2}|t-s+\ep|^{2-(\alpha_2+\beta)}\nonumber\\
&\qquad\qquad+\ep^{2-\beta}|t-s|\nonumber\\
&\qquad\qquad+|t-s|^{\alpha_2}\left[\ep^{2-(\alpha_2+\beta)} |t-s+\ep| + \ep^{2-\frac{\alpha_2+\beta}{2}}|t-s+\ep|^{1-\frac{\alpha_2+\beta}{2}}\right]\Big]^p\nonumber\\
&\qquad\qquad+1_{\beta\in\,]0,1[}\Big[\ep^{1-\beta}|t-s|^2\nonumber\\ \nonumber
&\qquad\qquad+|t-s|^{1+\alpha_3}\big[\big(|t-s+\ep|^{2-(\alpha_3+\beta)}-|t-s|^{2-(\alpha_3+\beta)}\big)\\
&\qquad\qquad+\ep^{1-\frac{\alpha_3+\beta}{2}}|t-s+\ep|^{1-\frac{\alpha_3+\beta}{2}}+\ep^{1-\alpha_3-\beta}|t-s+\ep|\big]\nonumber\\ \nonumber
&\qquad\qquad+|t-s|^{\alpha_4} \big[\ep^{2-(\alpha_4+\beta)}|t-s+\ep|\\
&\qquad\qquad+\ep^{2-\frac{\alpha_4+\beta}{2}}|t-s+\ep|^{1-\frac{\alpha_4+\beta}{2}}
+\ep^2|t-s+\ep|^{1-(\alpha_4+\beta)}\big]\Big]^p\Big\},
\label{fla6-2}
\end{align}
where $\alpha_i\in\,]0,1[$, $i=1,2,3$ and  $\alpha_1\in\,]0,(2\wedge k) - \beta[$, $\alpha_2 \in\,]0,2-\beta[$, $\alpha_3 \in\,]0,(2\wedge k)-\beta[$, $\alpha_4 \in\,]0,(2\wedge k)-\beta[$.
\end{prop}

\noindent{\it Proof}. We follow the scheme of the proof of Theorem 3.8 in \cite{dss09} (see in particular the analysis of the term $T^n_2(t,\bar t,x)$ in that Theorem) with suitable modifications of
the lemmas applied in that proof (see Lemmas \ref{la6-2-a1}, \ref{la6-2-a2}, \ref{la6-2-a3} below). According to \eqref{a22a} and the arguments of \cite{dss09}, page 28, we write
\begin{align*}
& E\left(\int_{s-\ep}^s dr \left\Vert \left[G(t-r,x-\ast)-G(s-r,x-\ast)\right] Z(r,\ast)\right\Vert_\hac^2\right)^p\\
&\quad \le C \sum_{i=1}^4 E\left(\left\vert R_i^\ep(s,t,x)\right\vert^p\right)
\end{align*}
where 
\beqn
R_i^\ep(s,t,x)=\int_{s-\ep}^s dr\int_{\rk}\int_{\rk} G(s-r,du)G(s-r,dv) r_i(s,t,r,x,u,v),
\eeqn
and
\begin{align*}
r_1(s,t,r,x,u,v)&= f\left(\frac{t-r}{s-r}v-u\right)\frac{t-r}{s-r}\\
&\quad\times\left[Z\left(r,x-\frac{t-r}{s-r}u\right)-Z(r,x-u)\right]\\
&\quad\times \left[Z\left(r,x-\frac{t-r}{s-r}v\right)-Z(r,x-v)\right],
\end{align*}
\begin{align*}
r_2(s,t,r,x,u,v)&=\left[\left(\frac{t-r}{s-r}\right)^2 f\left(\frac{t-r}{s-r}(v-u)\right)-\frac{t-r}{s-r}f\left(\frac{t-r}{s-r}v-u\right)\right]\\
&\quad\times Z\left(s,x-\frac{t-r}{s-r}u\right)\left[Z\left(r,x-\frac{t-r}{s-r}v\right)-Z(r,x-v)\right],
\end{align*}
\begin{align*}
r_3(s,t,r,x,u,v)&=\left[\left(\frac{t-r}{s-r}\right)^2 f\left(\frac{t-r}{s-r}(v-u)\right)-\frac{t-r}{s-r}f\left(v-\frac{t-r}{s-r}u\right)\right]\\
&\quad\times Z(r,x-v)\left(Z\left(r,x-\frac{t-r}{s-r}u\right)-Z(r,x-u)\right),
\end{align*}
\begin{align*}
r_4(s,t,r,x,u,v)&=\left[\left(\frac{t-r}{s-r}\right)^2 f\left(\frac{t-r}{s-r}(v-u)\right)-\frac{t-r}{s-r}f\left(\frac{t-r}{s-r}v-u\right)\right.\\
&\qquad\qquad -\left.\frac{t-r}{s-r}f\left(v-\frac{t-r}{s-r}u\right)+f(v-u)\right]\\
&\quad\times Z(r,x-u)Z(r,x-v),
\end{align*}
and $f(x) = \vert x\vert^{-\beta}$, $ \beta \in \,]0,2\wedge k[$.

Set
\beqn
\nu_1(s,t,\ep)=\int_{s-\ep}^s dr \int_{\rk\times \rk} G(s-r,du)G(s-r,dv)f\left(\frac{t-r}{s-r}v-u\right)\frac{t-r}{s-r}.
\eeqn
By Lemma \ref{la6-2-a1} below, we can apply first H\"older's inequality and then the Cauchy-Schwarz inequality to obtain
\begin{align*}
&E\left(\left\vert R_1^\ep(s,t,x)\right\vert^p\right)
\le  \left(\nu_1(s,t,\ep)\right)^{p-1}\\
&\quad\times \int_{s-\ep}^s dr\int_{\rk\times\rk} G(s-r,du)G(s-r,dv)f\left(\frac{t-r}{s-r}v-u\right)\frac{t-r}{s-r}\\
&\qquad\qquad\times \left(E\left(\left\vert Z\left(r,x-\frac{t-r}{s-r}u\right)-Z(r,x-u)\right\vert^{2p}\right)\right)^{\frac{1}{2}}\\
&\qquad\qquad\times \left(E\left(\left\vert Z\left(r,x-\frac{t-r}{s-r}v\right)-Z(r,x-v)\right\vert^{2p}\right)\right)^{\frac{1}{2}}.
\end{align*}
By the assumptions on the process $Z$, the product of the last two factors is bounded by
\beqn
\left\vert\frac{t-s}{s-r}\right\vert^{2\rho p}\vert u\vert^{\rho p} \vert v\vert^{\rho p},
\eeqn
and for $|u|\le s-r$, $|v|\le s-r$, this is bounded by $|t-s|^{2\rho p}$.
Consequently, by applying Lemma \ref{la6-2-a1}, we obtain
\begin{align}
\label{i4}
E\left(\left\vert R_1^\ep(s,t,x)\right\vert^p\right)&\le C \ep^{[\frac{3-\beta}{2}]p}|t-s|^{2\rho p}
\left((t-s+\ep)^{3-\beta}-(t-s)^{3-\beta}\right)^{\frac{p}{2}}\nonumber\\
&\le C  \ep^{[\frac{4-\beta}{2}]p}|t-s|^{2\rho p}|t-s+\ep|^{\frac{2-\beta}{2}p}.
\end{align}

Set
\begin{align*}
\nu_2(s,t,\ep)&=\int_{s-\ep}^s dr \int_{\rk\times\rk} G(s-r,du) G(s-r,dv)\\
&\quad\times\left\vert\left(\frac{t-r}{s-r}\right)^2 f\left(\frac{t-r}{s-r}(v-u)\right)-\frac{t-r}{s-r}f\left(\frac{t-r}{s-r}v-u\right)\right\vert.
\end{align*}
By applying Lemma \ref{la6-2-a2}, H\"older's inequality, the Cauchy-Schwarz inequality and the properties of $Z$, we have
\begin{align}
E\left(\left\vert R_2^\ep(s,t,x)\right\vert^p\right)&\le C |t-s|^{\rho p}\left(\nu_2(s,t,\ep)\right)^p\nonumber\\
&\le C |t-s|^{\rho p}\nonumber\\
&\qquad\times\left[\ep|t-s||t-s+\ep|^{1-\beta}1_{\beta\in\,]0,1[}+\ep|t-s|^{2-\beta}1_{\beta\in[1,2[}\right.\nonumber\\
&\qquad\qquad \left.+ \ep|t-s|^{\alpha}|t-s+\ep|^{2-(\alpha+\beta)}\right.\nonumber\\
&\qquad\qquad \left.+\ep^{2-\frac{\alpha+\beta}{2}}|t-s|^\alpha|t-s+\ep|^{1-\frac{\alpha+\beta}{2}}\right]^p,
\label{i5}
\end{align}
for any $\alpha\in\,]0,1[$ such that $\alpha+\beta\in\,]0,2\wedge k[$.
The same result holds for $E\left(\left\vert R_3^\ep(s,t,x)\right\vert^p\right)$.


For the study of the term $R_4^\ep(s,t,x)$, we use the same approach as for the preceding terms.
We define
\begin{align*}
\nu_3(s,t,\ep)&=\int_{s-\ep}^s dr \int_{\rk\times \rk} G(s-r,du) G(s-r,dv)\\
&\qquad\times \left[\left(\frac{t-r}{s-r}\right)^2 f\left(\frac{t-r}{s-r}(v-u)\right)-\frac{t-r}{s-r}f\left(\frac{t-r}{s-r}v-u\right)\right.\\
&\qquad\qquad -\left.\frac{t-r}{s-r}f\left(v-\frac{t-r}{s-r}u\right)+f(v-u)\right].
\end{align*}
Then, using Lemma \ref{la6-2-a3} and the hypotheses on $Z$, we have the following.

If $k\in\{2,3\}$ and $\beta\in[1,2[$,
\begin{align}
&E\left(\left\vert R_4^\ep(s,t,x)\right\vert^p\right) \nonumber\\
&\qquad \le C \Big[\ep|t-s|^{2-\beta}+\ep |t-s|^\alpha|t-s+\ep|^{2-(\alpha+\beta)}+\ep^{2-\beta}|t-s|\nonumber\\
&\qquad\qquad+|t-s|^\alpha\left[\ep^{2-(\alpha+\beta)} |t-s+\ep| + \ep^{2-\frac{\alpha+\beta}{2}}|t-s+\ep|^{1-\frac{\alpha+\beta}{2}}\right]\Big]^p,
\label{i9}
\end{align}
for any $\alpha\in\,]0,1[$ with $\alpha+\beta\in\,]0,2[$.

If $k\in\{1,2,3\}$ and $\beta\in\,]0,1[$, then
\begin{align}\nonumber
&E\left(\left\vert R_4^\ep(s,t,x)\right\vert^p\right)\\ \nonumber
&\qquad \le  C \Big[ \ep^{1-\beta}|t-s|^2\nonumber\\ \nonumber
&\qquad\qquad+|t-s|^{1+\alpha}\big[[|t-s+\ep|^{2-(\alpha+\beta)}-|t-s|^{2-(\alpha+\beta)}]\\ \nonumber
&\qquad\qquad\qquad\qquad\qquad +\ep^{1-\frac{\alpha+\beta}{2}}|t-s+\ep|^{1-\frac{\alpha+\beta}{2}} +\ep^{1-\alpha-\beta}|t-s+\ep|\big]\nonumber\\ \nonumber
&\qquad\qquad+|t-s|^{\bar \alpha} \big[\ep^{2-(\bar\alpha+\beta)}|t-s+\ep|+\ep^{2-\frac{\bar\alpha+\beta}{2}}|t-s+\ep|^{1-\frac{\bar\alpha+\beta}{2}}\\
&\qquad\qquad\qquad\qquad\qquad
+\ep^2|t-s+\ep|^{1-(\bar\alpha+\beta)}\big]\Big]^p,
\label{i91}
\end{align}
for any $\alpha \in\,]0,(k-\beta)\wedge1[$ and $\bar\alpha \in\,]0,(2\wedge k)-\beta[$.

The estimates \eqref{i4}, \eqref{i5}, \eqref{i9} and \eqref{i91} give \eqref{fla6-2}. \hfill\qed
\bigskip

   The three lemmas below were used in the proof of Proposition \ref{la6-2}.
   
\begin{lem}
\label{la6-2-a1} For any $0\le s<t\le T$, $\ep\in\,]0,s[$, set
\beqn
\nu_1(s,t,\ep)=\int_{s-\ep}^s dr \int_{\rk\times\rk}G(s-r,du) G(s-r,dv)f\left(\frac{t-r}{s-r}v-u\right)\frac{t-r}{s-r}.
\eeqn
There exists a constant $C>0$ such that
\beqn
\nu_1(s,t,\ep)\le C
\ep^{\frac{3-\beta}{2}}\left((t-s+\ep)^{3-\beta}-(t-s)^{3-\beta}\right)^{\frac{1}{2}}.
\eeqn
\end{lem}
\noindent{\it Proof}. We follow the proof of Lemma 6.3 in \cite{dss09}. Consider the change of variable 
$v\mapsto \frac{t-r}{s-r}v$, which maps $G(s-r,\cdot)$ into $\frac{s-r}{t-r}G(t-r, \cdot)$. Hence
\beqn
\nu_1(s,t,\ep)
= \int_{s-\ep}^s dr \int_{\rk}G(s-r,du)\int_{\rk}G(t-r,dv) f(u-v).
\eeqn
Using \eqref{a22a}, Schwarz's inequality yields
\begin{align*}
 &\nu_1(s,t,\ep)\\
 &\quad= \int_{s-\ep}^s dr\int_{\rk} d\xi \frac{\vert \tf G(s-r)(\xi)\vert\vert\tf G(t-r)(\xi)}{|\xi|^{k-\beta}}\\
 &\quad=\int_{s-\ep}^s dr\int_{\rk} d\xi \frac{\vert \tf G(s-r)(\xi)\vert}{|\xi|^{\frac{k-\beta}{2}}}\frac{\vert \tf G((t-r)(\xi)\vert}{|\xi|^{\frac{k-\beta}{2}}}\\
 &\quad\le \left(\int_{s-\ep}^s dr\int_{\rk} d\xi \frac{\vert \tf G(s-r)(\xi)\vert^2}{|\xi|^{k-\beta}}\right)^{\frac{1}{2}}
 \left(\int_{s-\ep}^s dr\int_{\rk} d\xi \frac{\vert \tf G(t-r)(\xi)\vert^2}{|\xi|^{k-\beta}}\right)^{\frac{1}{2}}.
 \end{align*}
 According to Lemma \ref{la1}, 
 \beqn
 \int_{s-\ep}^s dr\int_{\rk} d\xi\frac{\vert \tf G(s-r)(\xi)\vert^2}{|\xi|^{k-\beta}}=\int_0^\ep dr \Vert G(r,\ast)\Vert_{\hac}^2\le C \ep^{3-\beta}.
 \eeqn
 and 
 \begin{align*}
\int_{s-\ep}^s dr\int_{\rk} d\xi\frac{\vert \tf G(t-r)(\xi)\vert^2}{|\xi|^{k-\beta}}&=\int_{s-\ep}^sdr \Vert G(t-r,\ast)\Vert_{\hac}^2\\
&\le C\left((t-s+\ep)^{3-\beta}-(t-s)^{3-\beta}\right),
\end{align*}
which proves the lemma.
\hfill\qed


\begin{lem}
 \label{la6-2-a2}
  For any $0\le s<t\le T$ and $\ep\in\,]0,s[$, let 
 \begin{align*}
\nu_2(s,t,\ep)&=\int_{s-\ep}^s dr \int_{\rk\times \rk} G(s-r,du) G(s-r,dv)\\
&\quad\times\left\vert\left(\frac{t-r}{s-r}\right)^2 f\left(\frac{t-r}{s-r}(v-u)\right)-\frac{t-r}{s-r}f\left(\frac{t-r}{s-r}v-u\right)\right\vert.
\end{align*}
Then, for any $\alpha\in\,]0,1[$ such that $\alpha+\beta\in\,]0,2\wedge k[$, 
\begin{align*}
\nu_2(s,t,\ep)&\le C \left[\ep|t-s||t-s+\ep|^{1-\beta}1_{\beta\in\,]0,1[}+\ep|t-s|^{2-\beta}1_{\beta\in[1,2[}\right.\\
&\left.+ \ep|t-s|^{\alpha}|t-s+\ep|^{2-(\alpha+\beta)}+\ep^{2-\frac{\alpha+\beta}{2}}|t-s|^\alpha|t-s+\ep|^{1-\frac{\alpha+\beta}{2}}\right].
\end{align*}
 \end{lem}
 \noindent{\it Proof}. We follow the scheme of the proof of Lemma 6.4 in \cite{dss09}. Let $\nu_2^1(s,t,\ep)$ and $\nu_2^2(s,t,\ep)$ be defined in the same way as $\nu_2(s,t,\ep)$ but with
 the expression in absolute values replaced by 
 \beqn
 \left\vert\left(\frac{t-r}{s-r}\right)^2-\frac{t-r}{s-r}\right\vert f\left(\frac{t-r}{s-r}(v-u)\right),
 \eeqn
 and
 \beqn
 \frac{t-r}{s-r}\left\vert f\left(\frac{t-r}{s-r}(v-u)\right)-f\left(\frac{t-r}{s-r}v-u\right)\right\vert,
 \eeqn
 respectively.
 
 The change of variables $(u,v)\mapsto \frac{t-r}{s-r}(u,v)$, \eqref{a22a} and \eqref{basic} yield
 \begin{align*}
 \nu_2^1(s,t,\ep)&=|t-s|\int_{s-\ep}^s \frac{dr}{t-r}\int_{\rk\times\rk}G(t-r,du)G(t-r,dv) f(v-u)\\
 &\le C |t-s|\int_{s-\ep}^s \frac{dr}{t-r}|t-r|^{2-\beta}\\
 &\le C |t-s|\left((t-s+\ep)^{2-\beta}-(t-s)^{2-\beta}\right).
 \end{align*}
 By the Intermediate Value Theorem,
 \beq
 \label{incre}
  (t-s+\ep)^{2-\beta}-(t-s)^{2-\beta}\le
  \begin{cases}
  C\ep(t-s+\ep)^{1-\beta}, & {\rm if}\ \beta < 1,\\
  C\ep(t-s)^{1-\beta}, & {\rm if}\ \beta \geq 1.
  \end{cases}
\eeq
Thus,
\beq
\label{i7}
 \nu_2^1(s,t,\ep)\le C \left[\ep|t-s||t-s+\ep|^{1-\beta}1_{\beta\in\,]0,1[}+\ep|t-s|^{2-\beta}1_{\beta\in[1,2[}\right].
 \eeq
 Consider the change of variables $(u,v)\mapsto\left(u,\frac{t-r}{s-r}v\right)$. Then
 \beqn
 \nu_2^2(s,t,\ep)=\int_{s-\ep}^s dr \int_{\rk\times\rk} G(s-r,du) G(t-r,dv)\left\vert Df\left(v-u,-\frac{t-s}{s-r}u\right)\right\vert.
 \eeqn
 Let $\alpha>0$ be such that $\alpha+\beta\in\,]0,2\wedge k[$. We write $f(u)=\kappa_\beta(u) = \vert u \vert^{-\beta}$, $\beta \in \,]0,2\wedge k[$. By applying \cite[Lemma 2.6(a)]{dss09} with $d:=k$, $b:=\alpha$, $a:=k-(\alpha+\beta)$, $c:=t-s$, $u:=v-u$, $x=:-\frac{u}{s-r}$, we obtain
\begin{align}\nonumber
  \nu_2^2(s,t,\ep) &\le C|t-s|^{\alpha}\int_{s-\ep}^s dr \int_{\rk\times \rk} G(s-r,du) G(t-r,dv)\\
 &\qquad\qquad\times \int_{\rk} dw\, \kappa_{\alpha+\beta}((t-s)w-(v-u))\left\vert D\kappa_{k-\alpha}\left(w,-\frac{u}{s-r}\right)\right\vert.
 \label{i7rd1}
\end{align}
Next, we establish upper bounds for the integral above.  For this, we proceed as in the proof of  \cite[Lemma 6.4]{dss09} (see the analysis of the terms denoted by
$\nu_2^{n,2,1,1}$, $\nu_2^{n,2,1,2}$ in that proof).

   We consider first the domain $\vert w \vert \leq 3$, in which we split the absolute value in \eqref{i7rd1} into two separate terms. Following the calculations just mentioned in the proof of \cite[Lemma 6.4]{dss09} and using \eqref{a22a}, we find that
\begin{align*}
&\int_{s-\ep}^s dr \int_{\rk\times\rk} G(s-r,du) G(t-r,dv)\\
 &\qquad\times \int_{|w|\le 3} dw\, \kappa_{\alpha+\beta}((t-s)w-(v-u))\kappa_{k-\alpha}\left(\frac{u}{s-r}-w\right)\\
 &\le \left(\int_{|w|\le 4}\frac{dw}{|w|^{k-\alpha}}\right) \int_{s-\ep}^s dr \int_{\rk} d\xi \frac{\vert \tf G(t-r)(\xi)\vert^2}{|\xi|^{k-(\alpha+\beta)}}.
 \end{align*}
 We are assuming $\alpha+\beta\in\,]0,2\wedge k[$. Applying Lemma \ref{la1} with $\beta:=\alpha+\beta$, we see that this last expression is bounded by
 \begin{align*}
 &C \left(\int_{|w|\le 4}\frac{dw}{|w|^{k-\alpha}}\right) 
  \left[(t-s+\ep)^{3-(\alpha+\beta)}-(t-s)^{3-(\alpha+\beta)}\right]\\
 &\qquad \le C \ep(t-s+\ep)^{2-(\alpha+\beta)}.
 \end{align*}
 Similarly,
  \begin{align}
&\int_{s-\ep}^s dr \int_{\rk\times\rk} G(s-r,du) G(t-r,dv)\nonumber\\
 &\qquad\qquad\times \int_{|w|\le 4} dw\, \kappa_{\alpha+\beta}((t-s)w-(v-u))\kappa_{k-\alpha}(w)\nonumber\\
 &\quad \le \left(\int_{|w|\le 4}\frac{dw}{|w|^{k-\alpha}}\right)\int_{s-\ep}^s dr \int_{\rk} d\xi\frac{\vert\tf G(s-r)(\xi)||\tf G(t-r)(\xi)|}{|\xi|^{k-(\alpha+\beta)}}\label{i6}.
 \end{align}
 By applying the Cauchy-Schwarz inequality to the $d\xi$-integral, and using Lemma \ref{la1}, we find that
 \beq
 \label{cs}
 \int_{\rk} d\xi\frac{\vert\tf G(s-r)(\xi)||\tf G(t-r)(\xi)|}{|\xi|^{k-(\alpha+\beta)}}\le C(s-r)^{1-\frac{\alpha+\beta}{2}}(t-r)^{1-\frac{\alpha+\beta}{2}}.
 \eeq
 Hence \eqref{i6} is bounded by
 \begin{align*}  
&C\left(\int_{|w|\le 4}\frac{dw}{|w|^{k-\alpha}}\right)\int_{s-\ep}^s dr (s-r)^{1-\frac{\alpha+\beta}{2}}(t-r)^{1-\frac{\alpha+\beta}{2}}\\
&\qquad \le C \left(\int_{|w|\le 4}\frac{dw}{|w|^{k-\alpha}}\right) \left(\ep^{2-\frac{\alpha+\beta}{2}}(t-s+\ep)^{1-\frac{\alpha+\beta}{2}}\right).
\end{align*}

In order to finish the estimate of $\nu_2^2(s,t,\ep)$, we consider now the domain $|w|>3$ (see \cite{dss09}, page 59). For this, we apply the
inequality
\beq\label{rd6.23}
\left\vert D\kappa_{k-\alpha}\left(w,-\frac{u}{s-r}\right)\right\vert\le C \int_0^1d\lambda\, \kappa_{k+1-\alpha}\left(w-\lambda\frac{u}{s-r}\right),
\eeq
which yields, using \eqref{a22a},
\begin{align*}
&\int_{s-\ep}^s dr \int_{\rk\times\rk} G(s-r,du) G(t-r,dv)\\
 &\qquad\qquad\times \int_{|w|>3} dw\, \kappa_{\alpha+\beta}((t-s)w-(v-u))\left\vert D\kappa_{k-\alpha}\left(w,-\frac{u}{s-r}\right)\right\vert\\
 &\qquad\le C\left(\int_{|w|>3}\frac{dw}{|w|^{k+1-\alpha}}\right) \int_{s-\ep}^s dr \int_{\rk} d\xi \frac{\vert \tf G(s-r)(\xi)\vert \vert \tf G(t-r)(\xi)\vert}{|\xi|^{k-(\alpha+\beta)}}.
 \end{align*}
 By \eqref{cs}, this is bounded above by
\beq\label{rd6.23a}
 C \left(\int_{|w|>3}\frac{dw}{|w|^{k+1-\alpha}}\right)\ep^{2-\frac{\alpha+\beta}{2}}(t-s+\ep)^{1-\frac{\alpha+\beta}{2}}.
\eeq
Since $\alpha\in\,]0,1[$ by hypothesis, the $dw$-integral 
is finite. We have thus established that
 \beq
 \label{i8}
 \nu_2^2(s,t,\ep)\le C |t-s|^\alpha\left\{\ep |t-s+\ep|^{2-(\alpha+\beta)}+\ep^{2-\frac{\alpha+\beta}{2}}|t-s+\ep|^{1-\frac{\alpha+\beta}{2}}\right\}.
 \eeq
With \eqref{i7} and \eqref{i8}, the lemma is proved. 
\hfill\qed
\begin{lem}
\label{la6-2-a3} For any $0\le s<t\le T$ and every $0\le \ep<s$, let 
\begin{align*}
\nu_3(s,t,\ep)&=\int_{s-\ep}^s dr \int_{\rk}\int_{\rk} G(s-r,du) G(s-r,dv)\\
&\qquad\times \Big[\left(\frac{t-r}{s-r}\right)^2 f\left(\frac{t-r}{s-r}(v-u)\right)-\frac{t-r}{s-r}f\left(\frac{t-r}{s-r}v-u\right)\\
&\qquad\qquad -\frac{t-r}{s-r}f\left(v-\frac{t-r}{s-r}u\right)+f(v-u)\Big].
\end{align*}
We have the following.
\begin{enumerate}
\item For $k \in \{2,3\}$ and $\beta\in[1,2[$, 
\begin{align*}
\nu_3(s,t,\ep)&\le C \left[\ep|t-s|^{2-\beta}+\ep|t-s|^\alpha |t-s+\ep|^{2-(\alpha+\beta)}\right.+\ep^{2-\beta}|t-s|\nonumber\\
&\left.\qquad+|t-s|^\alpha\left[\ep^{2-(\alpha+\beta)} |t-s+\ep| + \ep^{2-\frac{\alpha+\beta}{2}}|t-s+\ep|^{1-\frac{\alpha+\beta}{2}}\right]\right],
\end{align*}
for any $\alpha\in\,]0,1[$ such that $\alpha+\beta\in\,]0,2[$.
\item For $k \in \{1,2,3\}$ and $\beta\in\,]0,1[$, then
\begin{align*}
\nu_3(s,t,\ep)&\le C \left[ \ep^{1-\beta}|t-s|^2\right.\nonumber\\
&\left.\qquad+|t-s|^{1+\alpha}\left[\left((t-s+\ep)^{2-(\alpha+\beta)}-(t-s)^{2-(\alpha+\beta)}\right)\right.\right.\\
&\qquad\qquad\qquad\qquad \left.\left.+\ep^{1-\frac{\alpha+\beta}{2}}|t-s+\ep|^{1-\frac{\alpha+\beta}{2}}+\ep^{1-\alpha-\beta}|t-s+\ep|\right]\right.\nonumber\\
&\qquad+|t-s|^{\bar\alpha} \left[\ep^{2-(\bar\alpha+\beta)}|t-s+\ep|+\ep^{2-\frac{\bar\alpha+\beta}{2}}|t-s+\ep|^{1-\frac{\bar\alpha+\beta}{2}} \right.\\
&\qquad\qquad\qquad\qquad\left.\left.
+\ep^2|t-s+\ep|^{1-(\bar\alpha+\beta)}\right]\right],
\end{align*}
for any $\alpha\in\,]0,((2\wedge k)-\beta)\wedge 1[$ and $\bar\alpha\in\,]0,(2\wedge k)-\beta[$.
\end{enumerate}
\end{lem}

\noindent{\it Proof}. We will follow the approach of \cite[Lemma 6.5]{dss09}. 
\bigskip

\noindent{\it Case 1}: $k \in \{2,3\}$ and $\beta\in[1,2[$. Set
\begin{align*}
\Delta^1f(r,s,t,u,v)&=\left\vert\left(\frac{t-r}{s-r}\right)^2 f\left(\frac{t-r}{s-r}(v-u)\right)-\frac{t-r}{s-r}f\left(\frac{t-r}{s-r}v-u\right)\right\vert,\\
\Delta^2f(r,s,t,u,v)&=\left\vert 1-\frac{t-r}{s-r}\right\vert f(v-u),\\
\Delta^3f(r,s,t,u,v)&=\frac{t-r}{s-r}\left\vert f(v-u)-f\left(v-\frac{t-r}{s-r}u\right)\right\vert.
\end{align*}
We consider the bound 
\beqn
\nu_3(s,t,\ep)\leq \nu_3^1(s,t,\ep)+\nu_3^2(s,t,\ep)+\nu_3^3(s,t,\ep),
\eeqn
where each $\nu_3^i(s,t,\ep)$, $i=1,2,3$, is defined in the same way as $\nu_3(s,t,\ep)$ with the integrand 
\begin{align}
\label{Delta}
\Delta f:&=\left(\frac{t-r}{s-r}\right)^2 f\left(\frac{t-r}{s-r}(v-u)\right)-\frac{t-r}{s-r}f\left(\frac{t-r}{s-r}v-u\right)\nonumber\\
&\quad \quad -\frac{t-r}{s-r}f\left(v-\frac{t-r}{s-r}u\right)+f(v-u)
\end{align}
replaced by $\Delta^if(r,s,t,u,v)$, $i=1,2,3$, respectively.

Notice that $\nu_3^1(s,t,\ep)$ coincides with $\nu_2(s,t,\ep)$ in Lemma \ref{la6-2-a2}, and we are assuming $\beta\in[1,2[$. Hence
\begin{align}
\nu_3^1(s,t,\ep)&\le C \left[\ep|t-s|^{2-\beta}+\ep|t-s|^\alpha |t-s+\ep|^{2-(\alpha+\beta)}\right.\nonumber\\
&\qquad\qquad\left.+\ep^{2-\frac{\alpha+\beta}{2}}|t-s|^\alpha|t-s+\ep|^{1-\frac{\alpha+\beta}{2}}\right], \label{i9bis}
\end{align}
 for any $\beta\in[1,2[$ and $\alpha\in\,]0,(2-\beta)\wedge 1[$.

By the definition of $\Delta^2f(r,s,t,u,v)$, using \eqref{a22a}, we have
\begin{align}\nonumber
\nu_3^2(s,t,\ep)&\le C |t-s| \int_{s-\ep}^s\frac{dr}{s-r}\int_{\rk} d\xi\frac{\vert\tf G(s-r)(\xi)\vert^2}{|\xi|^{k-\beta}}\\ \nonumber
&=C |t-s| \int_{s-\ep}^s dr (s-r)^{1-\beta}\\
&=C \ep^{2-\beta}|t-s| .
\label{i10}
\end{align}

Consider the change of variables $(u,v)\mapsto\left(\frac{t-r}{s-r}u,v\right)$. Then
\begin{align*}
\nu_3^3(s,t,\ep)&=\int_{s-\ep}^s dr\int_{\rk\times\rk} G(t-r,du)G(s-r,dv)\\
&\quad \times \left\vert Df\left(v-u,\frac{t-s}{t-r}u\right)\right\vert.
\end{align*}
This term is similar to $\nu_2^2(s,t,\ep)$ in the proof of Lemma \ref{la6-2-a2}, with the roles of $s$ and $t$ exchanged.
Thus, we will carry out similar calculations. For the sake of completeness, we give the details.

Let $\alpha>0$ be such that $\alpha+\beta\in\,]0,2[$. By applying \cite[Lemma 2.6(a)]{dss09} with $d:=k$, $b:=\alpha$, $a:=k-(\alpha+\beta)$, $c:=t-s$,
$u:=v-u$, $x:=\frac{u}{t-r}$, we obtain
\begin{align*}
\nu_3^3(s,t,\ep)&= |t-s|^\alpha\int_{s-\ep}^s dr\int_{\rk\times\rk} G(t-r,du)G(s-r,dv)\\
&\quad\times \int_{\rk} dw\, \kappa_{\alpha+\beta}((t-s)w-(v-u))\left\vert D\kappa_{k-\alpha}\left(w,\frac{u}{t-r}\right)\right\vert.
\end{align*}
We proceed as in the estimate of $\nu^2_2(s,t,\ep)$ (see \eqref{i7rd1}); however, note that the positions of the variables $s$ and $t$ are slightly different, and the bound that we will obtain is also different. We consider separately the domains $\vert w \vert \leq 3$ and $\vert w \vert > 3$. Consider the change of variables $w\mapsto w+\frac{u}{t-r}$, $u\mapsto \frac{s-r}{t-r}u$. Then, using \eqref{a22a},
\begin{align*}
&\int_{s-\ep}^s dr\int_{\rk\times\rk} G(t-r,du)G(s-r,dv)\\
&\qquad\qquad\qquad\times \int_{|w|\le 3} dw \, \kappa_{\alpha+\beta}((t-s)w-(v-u))\kappa_{k-\alpha}\left(w+\frac{u}{t-r}\right)\\
&\qquad \le C\left(\int_{|w|\le 4}\frac{dw}{|w|^{k-\alpha}}\right) \int_{s-\ep}^s dr \frac{t-r}{s-r}\int_{\rk} d\xi \frac{\vert\tf G(s-r)(\xi)\vert^2}{|\xi|^{k-(\alpha+\beta)}}\\
&\qquad \le C \int_{s-\ep}^s dr (t-r)(s-r)^{1-(\alpha+\beta)}
\\
&\qquad \le C \ep^{2-(\alpha+\beta)} |t-s+\ep|,
\end{align*}
with any $\alpha>0$ such that $\alpha+\beta\in\,]0,2[$ (we have used Lemma \ref{la1}).

By applying \eqref{i6} and \eqref{cs}, we see that 
\begin{align*}
&\int_{s-\ep}^s dr\int_{\rk\times\rk} G(t-r,du)G(s-r,dv)\\
&\qquad\qquad\qquad\times \int_{|w|\le 4} dw\, \kappa_{\alpha+\beta}((t-s)w-(v-u))\kappa_{k-\alpha}(w)\\
&\qquad \le C \ep^{2-\frac{\alpha+\beta}{2}}|t-s+\ep|^{1-\frac{\alpha+\beta}{2}},
\end{align*}
for any $\alpha>0$ with $\alpha+\beta\in\,]0,2[$.
Thus, for $\alpha$ and $\beta$ satisfying these conditions,
\begin{align*}
&\int_{s-\ep}^s dr\int_{\rk\times\rk} G(t-r,du)G(s-r,dv)\\
&\qquad\qquad\times \int_{|w|\le 3} dw\, \kappa_{\alpha+\beta}((t-s)w-(v-u))\left\vert D\kappa_{k-\alpha}\left(w,\frac{u}{t-r}\right)\right\vert\\
&\qquad \le C \left[\ep^{2-(\alpha+\beta)} |t-s+\ep|+\ep^{2-\frac{\alpha+\beta}{2}}|t-s+\ep|^{1-\frac{\alpha+\beta}{2}}\right].
\end{align*}
We continue the proof by dealing with the term
\begin{align*}
&\int_{s-\ep}^s dr\int_{\rk\times\rk} G(t-r,du)G(s-r,dv)\\
&\qquad\times \int_{|w|> 3} dw\, \kappa_{\alpha+\beta}((t-s)w-(v-u))\left\vert D\kappa_{k-\alpha}\left(w,\frac{u}{t-r}\right)\right\vert.
\end{align*}
Assume $\alpha\in\,]0,1[$ with $\alpha+\beta\in\,]0,2[$. Just as in the proof of Lemma \ref{la6-2-a2} (see \eqref{rd6.23a}), it is bounded above by $C\ep^{2-\frac{\alpha+\beta}{2}}|t-s+\ep|^{1-\frac{\alpha+\beta}{2}}$.

Hence, 
\beq
\label{i11}
\nu_3^3(s,t,\ep)\le C |t-s|^\alpha\left[\ep^{2-(\alpha+\beta)} |t-s+\ep|+\ep^{2-\frac{\alpha+\beta}{2}}|t-s+\ep|^{1-\frac{\alpha+\beta}{2}}\right].
\eeq

From \eqref{i9bis}--\eqref{i11} it follows that, for $\beta\in[1,2[$ and $\alpha+\beta\in\,]0,2[$, the conclusion of the lemma holds.
\medskip


\noindent{\it Case 2}: $k\in\{1,2,3\}$ and $\beta\in\,]0,1[$.
Here, we consider a new decomposition of $\Delta f=\sum_{i=1}^3\bar\Delta^if(r,s,t,u,v)$ (see \eqref{Delta}), as follows:
\begin{align*}
\bar\Delta^1f(r,s,t,u,v)&=\left(\left(\frac{t-r}{s-r}\right)^2-2\frac{t-r}{s-r}+1\right)f(v-u)\\
&=\left(\frac{t-s}{s-r}\right)^2f(v-u),
\end{align*}
\begin{align*}
\bar\Delta^2f(r,s,t,u,v)&=\left(\left(\frac{t-r}{s-r}\right)^2 - \frac{t-r}{s-r}\right)\\
&\quad\times \left(f\left(\frac{t-r}{s-r}(v-u)\right)-f(v-u)\right),
\end{align*}
\begin{align*}
\bar\Delta^3f(r,s,t,u,v)&=\frac{t-r}{s-r}\left(f\left(\frac{t-r}{s-r}(v-u)\right)-f\left(\frac{t-r}{s-r}v-u\right)\right.\\
&\left.\quad -f\left(v-\frac{t-r}{s-r}u\right)+f(v-u)\right).
\end{align*}
We define $\bar\nu_3^i(s,t,\ep)$ in the same way as $\nu_3^i(s,t,\ep)$, but with $\Delta^if(r,s,t,u,v)$ replaced by $\bar\Delta^if(r,s,t,u,v)$, $i=1,2,3$, respectively.

Using \eqref{a22a}, we see by \eqref{basic},
\begin{align}
\label{i12}
\bar\nu_3^1(s,t,\ep)&=\int_{s-\ep}^s dr \int_{\rk\times\rk} G(s-r,du)G(s-r,dv) \left(\frac{t-s}{s-r}\right)^2f(v-u)\nonumber\\
&\le C \ep^{1-\beta}|t-s|^2.
\end{align}
Set 
\begin{align*}
\bar\nu_{3}^{2,1}(s,t,\ep)&=\int_{s-\ep}^s dr \int_{\rk\times\rk} G(s-r,du)G(s-r,dv)\left(\left(\frac{t-r}{s-r}\right)^2 - \frac{t-r}{s-r}\right)\\
&\qquad\qquad\qquad\qquad\times \left\vert f\left(\frac{t-r}{s-r}(v-u)\right)-f\left(\frac{t-r}{s-r}v-u\right)\right\vert,\\
\bar\nu_{3}^{2,2}(s,t,\ep)&=\int_{s-\ep}^s dr \int_{\rk\times\rk} G(s-r,du)G(s-r,dv)\left(\left(\frac{t-r}{s-r}\right)^2 - \frac{t-r}{s-r}\right)\\
&\qquad\qquad\qquad\qquad\times\left\vert f\left(\frac{t-r}{s-r}v-u\right)-f(v-u)\right\vert,
\end{align*}
so that $\bar\nu_{3}^{2}(s,t,\ep) \leq \bar\nu_{3}^{2,1}(s,t,\ep) + \bar\nu_{3}^{2,2}(s,t,\ep)$.
With the change of variables $(u,v)\mapsto\left(u,\frac{t-r}{s-r}v\right)$,
\begin{align*}
\bar\nu_{3}^{2,1}(s,t,\ep)&=|t-s|\int_{s-\ep}^s \frac{dr}{s-r} \int_{\rk\times\rk} G(s-r,du)G(t-r,dv)\\
&\qquad\qquad\qquad\times \left\vert Df\left(v-u,-\frac{t-s}{s-r}u\right)\right\vert,
\end{align*}
while
\begin{align*}
\bar\nu_{3}^{2,2}(s,t,\ep)&=|t-s|\int_{s-\ep}^s dr \frac{t-r}{(s-r)^2}\int_{\rk\times\rk} G(s-r,du)G(s-r,dv)\\
&\qquad\qquad\qquad\times \left\vert Df\left(v-u,\frac{t-s}{s-r}v\right)\right\vert.
\end{align*}
Remember that $f(x):=\kappa_{\beta}(x)$. We handle $\bar\nu_{3}^{2,1}(s,t,\ep)$ in a manner similar to the term $\nu_2^2(s,t,\ep)$ (see \eqref{i7rd1}, and note the differences in the positions of the variables $s$ and $t$), as follows.
Firstly, we apply Lemma 2.6(a) of \cite{dss09} with $d:=k$, $b:=\alpha$, $a:=k-(\alpha+\beta)$, $c:=t-s$, $u:=v-u$, $x=-\frac{u}{s-r}$ to obtain
\begin{align*}
\bar\nu_{3}^{2,1}(s,t,\ep)&\le C|t-s|^{1+\alpha}\int_{s-\ep}^s \frac{dr}{s-r}\int_{\rk}\int_{\rk} G(s-r,du)G(t-r,dv)\\
&\quad\times \int_{\rk} dw\, \kappa_{\alpha+\beta}(v-u-(t-s)w)\left\vert D\kappa_{k-\alpha}\left(w,-\frac{u}{s-r}\right)\right\vert.
\end{align*}
Next, in the $dw$-integral, we consider first the domain $\{|w|\le 3\}$ and we split the term $\left\vert D\kappa_{k-\alpha}\left(w,-\frac{u}{s-r}\right)\right\vert$ into the sum of the absolute
values of its components, as we also did in the analysis of the term $\nu_2^2(s,t,\ep)$ in Lemma \ref{la6-2-a2}.
With the change of variables $w\mapsto w-\frac{u}{s-r}$,
$u\mapsto\frac{t-r}{s-r}u$, we obtain
\begin{align*}
&\int_{s-\ep}^s \frac{dr}{s-r}\int_{\rk\times\rk} G(s-r,du)G(t-r,dv)\\
&\qquad\qquad\qquad\times \int_{|w|\le 3} dw\, \kappa_{\alpha+\beta}(v-u-(t-s)w)\kappa_{k-\alpha}\left(w-\frac{u}{s-r}\right)\\
&\qquad\le C  \left(\int_{|w|\le 4} dw\,\kappa_{k-\alpha}(w)\right)\sup_{w\in\IR^k}\int_{s-\ep}^s \frac{dr}{t-r}\\
&\qquad\qquad\qquad\times \int_{\rk\times\rk} G(t-r,du)G(t-r,dv)\kappa_{\alpha+\beta}(v-u-(t-s)w).
\end{align*}
By \eqref{a22a} and Lemma \ref{la1}, this is bounded by
\begin{align*}
&C \left(\int_{|w|\le 4} dw\, \kappa_{k-\alpha}(w)\right)\int_{s-\ep}^s dr (t-r)^{1-(\alpha+\beta)}\\
&\qquad \le C [(t-s+\ep)^{2-(\alpha+\beta)}-(t-s)^{2-(\alpha+\beta)}].
\end{align*}
We also have
\begin{align*}
&\int_{s-\ep}^s \frac{dr}{s-r}\int_{\rk\times\rk} G(s-r,du)G(t-r,dv)\\
&\qquad\qquad\qquad\qquad\times \int_{|w|\le 3} dw\, \kappa_{\alpha+\beta}(v-u-(t-s)w)\kappa_{k-\alpha}(w)\\
&\qquad\le C \ep^{1-\frac{\alpha+\beta}{2}}(t-s+\ep)^{1-\frac{\alpha+\beta}{2}},
\end{align*}
where the last inequality has been obtained using \eqref{cs}.

To end with the analysis of $\bar\nu_{3}^{2,1}(s,t,\ep)$, we must now consider the term
\begin{align*} 
&\int_{s-\ep}^s \frac{dr}{s-r}\int_{\rk\times \rk} G(s-r,du)G(t-r,dv)\\
&\qquad\times \int_{|w|>3} dw\, \kappa_{\alpha+\beta}(v-u-(t-s)w)\left\vert D\kappa_{k-\alpha}\left(w,-\frac{u}{s-r}\right)\right\vert.
\end{align*}
We have already met a similar expression in Lemma \ref{la6-2-a2} and also in the first part of Lemma \ref{la6-2-a3} (see the study of the
term $\nu_3^3(s,t,\ep)$). Since $\alpha \in \,]0,1[$, this is bounded above by $C \ep^{1-\frac{\alpha+\beta}{2}}|t-s+\ep|^{1-\frac{\alpha+\beta}{2}}$. 

Considering the three preceding estimates, we obtain
\begin{align}\nonumber
\bar\nu_{3}^{2,1}(s,t,\ep)&\le C|t-s|^{1+\alpha} \big[\left((t-s+\ep)^{2-(\alpha+\beta)}-(t-s)^{2-(\alpha+\beta)}\right)\\
&\qquad\qquad\qquad\qquad +\ep^{1-\frac{\alpha+\beta}{2}}|t-s+\ep|^{1-\frac{\alpha+\beta}{2}}\big],\label{i13}
\end{align}
for $\alpha\in\,]0,(k-\beta)\wedge1[$.

Our next objective is to prove that
\beq
\label{i14}
\bar\nu_{3}^{2,2}(s,t,\ep)\le C |t-s|^{1+\alpha}\left[\ep^{1-\frac{\alpha+\beta}{2}}|t-s+\ep|^{1-\frac{\alpha+\beta}{2}}
+\ep^{1-\alpha-\beta}|t-s+\ep|\right].
\eeq
This will be achieved with similar arguments as those used to establish \eqref{i13}, but yet again with the variables $s$ and $t$ in slightly different positions (which lead in fact to a different bound).

Indeed, we apply Lemma 2.6(a) in \cite{dss09} with $d:=k$, $b:=\alpha$, $a:=k-(\alpha+\beta)$, $c:=t-s$, $u:=v-u$, $x:=\frac{v}{s-r}$ to obtain
\begin{align*}
\bar\nu_{3}^{2,2}(s,t,\ep)&=|t-s|^{1+\alpha}\int_{s-\ep}^s dr \frac{t-r}{(s-r)^2}\int_{\rk\times\rk} G(s-r,du)G(s-r,dv)\\
&\qquad\qquad\quad\times\int_{\rk} dw\, \kappa_{\alpha+\beta}(v-u-(t-s)w)\left\vert D\kappa_{k-\alpha}\left(w,\frac{v}{s-r}\right)\right\vert.
\end{align*}
We consider first the domain of $w$ where $|w|\le 3$ and split the term $\left\vert D\kappa_{k-\alpha}\left(w,\frac{v}{s-r}\right)\right\vert$ into the sum
of $\kappa_{k-\alpha}\left(w+\frac{v}{s-r}\right)$ and $\kappa_{k-\alpha}(w)$.

For the first one, we consider the change of variables $w\mapsto w+\frac{v}{s-r}$,
$v\mapsto \frac{t-r}{s-r}v$. Using \eqref{a22a}, we obtain
\begin{align*}
&\int_{s-\ep}^s dr \frac{t-r}{(s-r)^2}\int_{\rk\times\rk} G(s-r,du)G(s-r,dv)\\
&\qquad\qquad\times\int_{|w|\le 3} dw\, \kappa_{\alpha+\beta}(v-u-(t-s)w)\kappa_{k-\alpha}\left(w+\frac{v}{s-r}\right)\\
&\qquad\le C\int_{s-\ep}^s \frac{dr}{s-r}\int_{\rk} d\xi\,\frac{\vert \tf G(s-r)(\xi)\vert\vert\tf G(t-r)(\xi)\vert}{|\xi|^{k-(\alpha+\beta)}}\\
&\qquad\le C\ep^{1-\frac{\alpha+\beta}{2}}|t-s+\ep|^{1-\frac{\alpha+\beta}{2}},
\end{align*}
where in the last inequality, we have applied \eqref{cs}.

As for the term with the integrand $\kappa_{k-\alpha}(w)$, we have
\begin{align*}
&\int_{s-\ep}^s dr \frac{t-r}{(s-r)^2}\int_{\rk\times\rk} G(s-r,du)G(s-r,dv)\\
&\qquad\qquad\times\int_{|w|\le 3} dw \, \kappa_{\alpha+\beta}(v-u-(t-s)w)\kappa_{k-\alpha}(w)\\
&\qquad\le C \int_{s-\ep}^s dr\, (t-r)(s-r)^{-\alpha-\beta}\\
&\qquad\le C \ep^{1-\alpha-\beta}|t-s+\ep|.
\end{align*}
Finally, since $\alpha\in\,]0,1[$ and using \eqref{a22a}, the contribution of the term with the domain of integration  $|w|\ge 3$ is
\begin{align*}
&\int_{s-\ep}^s dr \frac{t-r}{(s-r)^2}\int_{\rk\times\rk} G(s-r,du)G(s-r,dv)\\
&\quad\times\int_{|w|\ge 3} dw\, \kappa_{\alpha+\beta}(v-u-(t-s)w)\left\vert D\kappa_{k-\alpha}\left(w,\frac{v}{s-r}\right)\right\vert\\
&\le C \ep^{1-\alpha-\beta}|t-s+\ep|,
\end{align*}
as can be checked by using the upper bound 
\beqn
\left\vert D\kappa_{k-\alpha}\left(w,\frac{v}{s-r}\right)\right\vert\le \int_0^1d\lambda\, \kappa_{k+1-\alpha}\left(w+\lambda\frac{v}{s-r}\right)
\eeqn
(see \eqref{rd6.23}). This finishes the proof of \eqref{i14}. 

   From \eqref{i13} and \eqref{i14}, we obtain
\begin{align}
\label{i15}
\bar\nu_{3}^2(s,t,\ep)&\le C|t-s|^{1+\alpha}\left[\left((t-s+\ep)^{2-(\alpha+\beta)}-(t-s)^{2-(\alpha+\beta)}\right)\right.\nonumber\\
&\left.\qquad\qquad\quad +\ep^{1-\frac{\alpha+\beta}{2}}|t-s+\ep|^{1-\frac{\alpha+\beta}{2}}+\ep^{1-\alpha-\beta}|t-s+\ep|\right].
\end{align}

The last part of the proof is devoted to the study of $\bar\nu_{3}^3(s,t,\ep)$. For this, we apply Lemma 2.6(e) in \cite{dss09} with the following choice
of parameters: $d:=k$, $b:=\bar\alpha$, $a:=k-(\bar\alpha+\beta)$, $c:=t-s$, $\tilde u:=v-u$, $x:=-\frac{u}{s-r}$, $y:=\frac{v}{s-r}$, with $\bar\alpha+\beta\in]0,k[$. Notice that, with the notation of that
Lemma,
\beqn
\bar\Delta^3f(r,s,t,u,v)=\frac{t-r}{s-r}\bar D^2\kappa_\beta(\tilde u,cx,cy),
\eeqn
where for a given function $g:{\rk}\to\mathbb{R}$, 
\beq
\label{i16}
\bar D^2g(u,x,y)=g(u+x+y)-g(u+x)-g(u+y)+g(u).
\eeq
By doing so, we obtain
\begin{align*}
\bar\nu_{3}^3(s,t,\ep)&= |t-s|^{\bar\alpha}\int_{s-\ep}^s dr \frac{t-r}{s-r}
\int_{\rk}\int_{\rk}G(s-r,du)G(s-r,dv)\\
&\quad\times \int_{\rk} dw\, \kappa_{\bar\alpha+\beta}(v-u-(t-s)w)\left\vert\bar D^2\kappa_{k-\bar\alpha}\left(w,-\frac{u}{s-r},\frac{v}{s-r}\right)\right\vert.
\end{align*}
Let us study this term. For this, we will first consider small values of $w$ ($|w|\le 5$) and split the term
$\left\vert\bar D^2\kappa_{k-\bar\alpha}\left(w,-\frac{u}{s-r},\frac{v}{s-r}\right)\right\vert$
into the sum of four terms, according to the definition \eqref{i16}.

Let
\begin{align*}
I_1(s,t,\ep):&= \int_{s-\ep}^s dr\frac{t-r}{s-r}\int_{\rk\times\rk}G(s-r,du)G(s-r,dv)\\
&\quad\times \int_{|w|<5} dw\, \kappa_{\bar\alpha+\beta}(v-u-(t-s)w)\kappa_{k-\bar\alpha}(w).
\end{align*}
By applying \eqref{basic}, we see that
\beqn
I_1(s,t,\ep)\le C \ep^{2-(\bar\alpha+\beta)}|t-s+\ep|.
\eeqn
Set 
\begin{align*}
I_2(s,t,\ep):&= \int_{s-\ep}^s dr\frac{t-r}{s-r}\int_{\rk\times\rk}G(s-r,du)G(s-r,dv)\\
&\qquad\times \int_{|w|<5} dw\, \kappa_{\bar\alpha+\beta}(v-u-(t-s)w)\kappa_{k-\bar\alpha}\left(w-\frac{u}{s-r}\right).
\end{align*}
By applying the change of variables $w\mapsto w-\frac{u}{s-r}$ and $u\mapsto \frac{t-r}{s-r}u$, we obtain
\begin{align*}
I_2(s,t,\ep)&\le C\left(\int_{|w|\le 6} dw\,\kappa_{k-\bar\alpha}(w)\right)\\
&\qquad\times \sup_{w\in\IR^k}\int_{s-\ep}^s dr \int_{\rk\times\rk}G(t-r,du)G(s-r,dv)\\
&\qquad \qquad\qquad\times \kappa_{\bar\alpha+\beta}(v-u-(t-s)w)\\
&\le C\ep^{2-\frac{\bar\alpha+\beta}{2}}|t-s+\ep|^{1-\frac{\bar\alpha+\beta}{2}},
\end{align*}
where the last inequality has been obtained using similar arguments as in \eqref{i6}, \eqref{cs}.

In a very similar way, we see that
\begin{align*}
I_3(s,t,\ep):&= \int_{s-\ep}^s dr\frac{t-r}{s-r}\int_{\rk\times\rk}G(s-r,du)G(s-r,dv)\\
&\qquad\times \int_{|w|<5} dw\, \kappa_{\bar\alpha+\beta}(v-u-(t-s)w)\kappa_{k-\bar\alpha}\left(w+\frac{v}{s-r}\right)\\
&\le C\ep^{2-\frac{\bar\alpha+\beta}{2}}|t-s+\ep|^{1-\frac{\bar\alpha+\beta}{2}}.
\end{align*}
Let 
\begin{align*}
I_4(s,t,\ep):&= \int_{s-\ep}^s dr\frac{t-r}{s-r}\int_{\rk\times\rk}G(s-r,du)G(s-r,dv)\\
&\qquad\times \int_{|w|<5} dw\, \kappa_{\bar\alpha+\beta}(v-u-(t-s)w)\kappa_{k-\bar\alpha}\left(w+\frac{v-u}{s-r}\right).
\end{align*}
With the change of variables $w\mapsto w+\frac{v-u}{s-r}$, $(u,v)\mapsto \frac{t-r}{s-r}(u,v)$, we obtain
\begin{align*}
I_4(s,t,\ep)&\le C\left(\int_{|w|<7}dw\, \kappa_{k-\bar\alpha}(w)\right)\int_{s-\ep}^s dr \frac{s-r}{t-r}\int_{\rk}d\xi
\frac{\vert\tf G(t-r)(\xi)\vert^2}{|\xi|^{k-(\bar\alpha+\beta)}}\\
&\le C\ep^2|t-s+\ep|^{1-(\bar\alpha+\beta)}.
\end{align*}
With the estimates obtained so far for $I_i(s,t,\ep)$, $i=1,\ldots,4$, we have proved that
\begin{align}
\label{i17}
&\int_{s-\ep}^s dr \frac{t-r}{s-r}
\int_{\rk\times\rk}G(s-r,du)G(s-r,dv)\nonumber\\
&\qquad\times \int_{|w|<5} dw\, \kappa_{\bar\alpha+\beta}(v-u-(t-s)w)\left\vert\bar D^2\kappa_{k-\bar\alpha}\left(w,-\frac{u}{s-r},\frac{v}{s-r}\right)\right\vert\nonumber\\
&\le C\left[\ep^{2-(\bar\alpha+\beta)}|t-s+\ep|+\ep^{2-\frac{\bar\alpha+\beta}{2}}|t-s+\ep|^{1-\frac{\bar\alpha+\beta}{2}}+\ep^2|t-s+\ep|^{1-(\bar\alpha+\beta)}\right].
\end{align}
To finish the proof, we have to study the remaining term
\begin{align*}
I_5(s,t,\ep):&=\int_{s-\ep}^s dr \frac{t-r}{s-r}
\int_{\rk\times\rk}G(s-r,du)G(s-r,dv)\nonumber\\
&\qquad\times \int_{|w|\ge5} dw\, \kappa_{\bar\alpha+\beta}(v-u-(t-s)w)\\
&\qquad\times \left\vert\bar D^2\kappa_{k-\bar\alpha}\left(w,-\frac{u}{s-r},\frac{v}{s-r}\right)\right\vert.
\end{align*}
Proceeding as in \cite[p.65]{dss09}, we obtain 
\begin{align*}
I_5(s,t,\ep)&\le C\int_{s-\ep}^s dr \frac{t-r}{s-r}\int_{\rk\times\rk}G(s-r,du)G(s-r,dv)\\
&\qquad\times \int_{|w|\ge5} dw\, \kappa_{\bar\alpha+\beta}(v-u-(t-s)w)\\
&\qquad\times \int_0^1 d\lambda\int_0^1 d\mu\,\kappa_{k+2-\bar\alpha}\left(w-\lambda\frac{u}{s-r}+\mu\frac{v}{s-r}\right)\\
&\le C \left(\int_{|w|\ge 3} \frac{dw}{|w|^{k+2-\bar\alpha}}\right)\sup_{w\in\IR^k}\int_{s-\ep}^s dr \frac{t-r}{s-r}\\
&\qquad\times \int_{\rk\times\rk}G(s-r,du)G(s-r,dv)\kappa_{\bar\alpha+\beta}(v-u-(t-s)w)\\
&\le C \left(\int_{|w|\ge 3} \frac{dw}{|w|^{k+2-\bar\alpha}}\right)\int_{s-\ep}^s dr \frac{t-r}{s-r}\int_{\rk} d\xi \frac{\vert\tf G(s-r)(\xi)\vert^2}{|\xi|^{k-(\bar\alpha+\beta)}}.
\end{align*}
The term $\int_{|w|\ge5} \frac{dw}{|w|^{k+2-\bar\alpha}}$ is finite for any $\bar\alpha\in\,]0,2[$. Moreover, for $\bar\alpha+\beta\in\,]0,2\wedge k[$,
\beqn
\int_{s-\ep}^s dr \frac{t-r}{s-r}\int_{\rk} d\xi\frac{\vert\tf G(s-r)(\xi)\vert^2}{|\xi|^{k-(\bar\alpha+\beta)}}\le C \ep^{2-(\bar\alpha+\beta)}|t-s+\ep|.
\eeqn
Thus,
\beq
\label{i18}
I_5(s,t,\ep)\le C \ep^{2-(\bar\alpha+\beta)}|t-s+\ep|.
\eeq
The estimates \eqref{i17} and \eqref{i18} yield
\begin{align}
&\bar\nu_3(s,t,\ep)\le C|t-s|^{\bar\alpha}\nonumber\\
&\qquad\left[\ep^{2-(\bar\alpha+\beta)}|t-s+\ep|+\ep^{2-\frac{\bar\alpha+\beta}{2}}|t-s+\ep|^{1-\frac{\bar\alpha+\beta}{2}}
+\ep^2|t-s+\ep|^{1-(\bar\alpha+\beta)}\right].\label{i19}
\end{align}

By considering \eqref{i12}, \eqref{i15} and \eqref{i19}, we obtain the conclusion of the lemma for $\beta\in\,]0,1[$, $\alpha \in \,]0,1\wedge (k-\beta)[$ and $\bar\alpha\in \,]0,(2\wedge k)-\beta[$.

This finishes the proof of the lemma.\hfill\qed
\medskip

As a consequence of Propositions \ref{la6-1} and \ref{la6-2}, we obtain the following.
\begin{cor}
\label{ca-1}
Let $k\in\{1,2,3\}$ and let $Z=\{Z(t,x), (t,x)\in [0,T]\times\rk\}$ be a stochastic process satisfying the hypotheses of Proposition \ref{la6-1}. Fix a compact
set $K\subset \rk$ and $p\in[1,\infty[$. Then there exists a constant $C>0$ such that for any $0\le s<t\le T$ and $\ep\in\,]0,s[$, 
\begin{align*}
&\sup_{x,y\in K\atop 0\le s<t\le T}\left\{E\left(\int_{s-\ep}^s dr \left\Vert\left[G(t-r,x-\ast)-G(s-r,y-\ast)\right] Z(r,\ast)\right\Vert^2_{\hac}\right)^p\right\}\\
&\qquad \le C \left[\ep^{p(1-\beta)^- }1_{\beta\in]0,1[}+ \ep^{p(2-\beta)^-}1_{\beta\in[1,2[}\right].
\end{align*}
\end{cor}
\noindent{\it Proof}. By taking the values of the parameters $\alpha$, $\bar\alpha$ in Proposition \ref{la6-1} and $\alpha_i$, $i=1,2,3$ in Proposition \ref{la6-2},
arbitrarily close to zero, and by bounding $\vert x-y\vert$ and $\vert t-s\vert$ by a constant, we obtain the result. Notice that for $\alpha_3=0^+$ and $\beta\in\,]0,1[$, 
\beqn
(t-s+\ep)^{2-(\alpha_3+\beta)}-(t-s)^{2-(\alpha_3+\beta)}\le C \ep(t-s+\ep)^{1-(\alpha_3+\beta)}.
\eeqn
This completes the proof.
\hfill\qed
\bigskip

In the proof of Theorem \ref{t5n}, we give an upper bound for the term $E\left(\left(\bar G_{33,\ep}\right)^p\right)$ by using Propositions \ref{la6-1} and \ref{la6-2} (see \eqref{280n}). In that case,  $Z(r,x):=\sigma_{i,j}(u(r,x))$, $\rho:=\frac{(2-\beta)^-}{2}$, $\ep:=\ep^{\theta_2}$, and it is assumed that $|x-y|\le \delta_0 \ep^{\theta}$, $|t-s|\le \ep^\theta$. 

Consider first the case $k=1$ and thus $\beta\in\,]0,1[$. By choosing $\alpha=\bar \alpha=(1-\beta)^-$, from Proposition \ref{la6-1} we obtain
\begin{align}
\label{kau1}
&E\left(\int_{s-\ep^{\theta_2}}^s dr \left\Vert\left[G(s-r,x-\ast)-G(s-r,y-\ast)\right] \sigma_{i,j}(u(r,\ast))\right\Vert_{\hac}^2\right)^p\nonumber\\
&\qquad \le C \left\{\ep^{\theta(2-\beta)^-+\theta_2}+\ep^{\theta(1-\beta)^-+2\theta_2}\right\}^p.
\end{align}

In the same context, and the choice of parameters $\alpha_1=\alpha_4=(1-\beta)^-$, $\alpha_3=0^+$, $\theta\geq \theta_2$, after tedious checking, Proposition \ref{la6-2} yields
\begin{align}
\label{kau2}
&E\left(\int_{s-\ep^{\theta_2}}^s dr \left\Vert\left[G(t-r,x-\ast)-G(s-r,x-\ast)\right] \sigma_{i,j}(u(r,\ast))\right\Vert_{\hac}^2\right)^p\nonumber\\
&\qquad\le C\left\{\ep^{\theta(2-\beta)^-+\theta_2}+\ep^{\theta^-+\theta_2(2-\beta)^-}
+\ep^{\theta(1-\beta)^-+2\theta_2} + \ep^{2\theta +\theta_2(1-\beta)}\right\}^p.
\end{align}
Consequently, from \eqref{kau1}, \eqref{kau2} we obtain
\begin{align}
\label{kau}
&E\left(\int_{s-\ep^{\theta_2}}^s dr \left\Vert\left[G(t-r,x-\ast)-G(s-r,y-\ast)\right] \sigma_{i,j}(u(r,\ast))\right\Vert^2_{\hac}\right)^p\nonumber\\
&\qquad\le C\left\{\ep^{\theta(2-\beta)^-+\theta_2}+\ep^{\theta^-+\theta_2(2-\beta)}
+\ep^{\theta(1-\beta)^-+2\theta_2} + \ep^{2\theta +\theta_2(1-\beta)}
\right\}^p,
\end{align}

 For $k \in \{2,3\}$, by choosing $\bar\alpha:=(2-\beta)^-$ and then $\alpha:=1^-$ if $\beta\in\,]0,1[$, and $\alpha:=(2-\beta)^-$ if $\beta\in[1,2[$, from Proposition \ref{la6-1} we obtain
\begin{align}
\label{ep1}
&E\left(\int_{s-\ep^{\theta_2}}^s dr \left\Vert\left[G(s-r,x-\ast)-G(s-r,y-\ast)\right] \sigma_{i,j}(u(r,\ast))\right\Vert_{\hac}^2\right)^p\nonumber\\
&\qquad \le C \ep^{[\theta(2-\beta)^-+\theta_2]p}.
\end{align}
In the same context of \eqref{280n}, for $\beta\in\,]0,1[$ and the choice of parameters $\alpha_1:=1^-$, $\alpha_3:=0^+$, $\alpha_4=1^-$, $\theta\geq \theta_2$, after tedious verification, Proposition \ref{la6-2} yields
\begin{align}
\label{ep2}
&E\left(\int_{s-\ep^{\theta_2}}^s dr \left\Vert\left[G(t-r,x-\ast)-G(s-r,x-\ast)\right] \sigma_{i,j}(u(r,\ast))\right\Vert_{\hac}^2\right)^p\nonumber\\
&\qquad\le C\left\{\ep^{\theta(2-\beta)^-+\theta_2}+\ep^{2\theta^-+\theta_2(1-\beta)}+\ep^{\theta^-+\theta_2(2-\beta)^-}
\right\}^p.
\end{align}

Consequently, from \eqref{ep1}, \eqref{ep2} we deduce 
\begin{align}
\label{ep3}
&E\left(\int_{s-\ep^{\theta_2}}^s dr \left\Vert\left[G(t-r,x-\ast)-G(s-r,y-\ast)\right] \sigma_{i,j}(u(r,\ast))\right\Vert^2_{\hac}\right)^p\nonumber\\
&\qquad\le C\left\{\ep^{\theta(2-\beta)^-+\theta_2}+\ep^{2\theta^-+\theta_2(1-\beta)}+\ep^{\theta^-+\theta_2(2-\beta)^-}
\right\}^p,
\end{align}
for $\beta\in\,]0,1[$ and $k\in\{2,3\}$.

Still for $k\in\{2,3\}$, let us now consider the case $\beta\in[1,2[$. By applying Proposition \ref{la6-2} with $\alpha_1:=(2-\beta)^-$ and $\alpha_2:=\frac{2-\beta}{2}$,
and assuming that $\theta>\theta_2$, we obtain
\begin{align}
\label{ep4}
&E\left(\int_{s-\ep^{\theta_2}}^s dr \left\Vert\left[G(t-r,x-\ast)-G(s-r,x-\ast)\right] \sigma_{i,j}(u(r,\ast))\right\Vert_{\hac}^2\right)^p\nonumber\\
&\qquad\le C\left\{\ep^{\theta(2-\beta)^-+\theta_2}+\ep^{\theta+\theta_2(2-\beta)}
+ \ep^{\theta\frac{2-\beta}{2}+\theta_2\frac{4-\beta}{2}}
\right\}^p.
\end{align}

Along with \eqref{ep1}, this yields
\begin{align}
\label{ep5}
&E\left(\int_{s-\ep^{\theta_2}}^s dr \left\Vert\left[G(t-r,x-\ast)-G(s-r,y-\ast)\right] \sigma_{i,j}(u(r,\ast))\right\Vert^2_{\hac}\right)^p\nonumber\\
&\qquad\le C\left\{\ep^{\theta(2-\beta)^-+\theta_2}+\ep^{\theta+\theta_2(2-\beta)}
+ \ep^{\theta\frac{2-\beta}{2}+\theta_2\frac{4-\beta}{2}}
\right\}^p.
\end{align}
\begin{lem}
\label{la7}
For any $\ep>0$ and $x\in \rk$,
\beqn
\inf_{s\in[0,T]}\int_0^\ep dr\, \Vert G(s+r,x-\ast)\Vert^2_{\hac}\ge C \ep^{3-\beta}.
\eeqn
\end{lem}
\noindent{\it Proof}. 
Let $c$ be the constant defined in Lemma \ref{la1}. With direct computations, we obtain
\begin{align*}
\int_0^\ep dr\, \Vert G(s+r,x-\ast)\Vert^2_{\hac}&=c\left(\frac{(s+\ep)^{3-\beta}-s^{3-\beta}}{3-\beta}\right)\\
&= c\ep\int_0^1 \left(\nu\ep+s\right)^{2-\beta}d\nu\\
&\ge c\ep\int_0^1(\nu\ep)^{2-\beta}d\nu= \frac{c}{3-\beta}\ep^{3-\beta}.
\end{align*}
Since the last expression does not depend on $s$, we obtain the result with the constant $C=\frac{c}{3-\beta}$.
\hfill\qed

  \begin{lem}
 \label{la10}
 Let $k\in\{1,2,3\}$. For $\ep>0$ and $0<h\le \ep$,
 \beqn
 \int_0^\ep dr\, \langle G(r,y-\ast),G(h+r,x-\ast)\rangle_{\hac}\le C \ep^{3-\beta}\varphi\left(\frac{x-y}{\ep},\frac{h}{\ep}\right),
 \eeqn
 where
 \beqn
 \lim_{|z|\to +\infty}\varphi(z,\lambda)=0,
 \eeqn
 uniformly in $\lambda\in[0,1]$.
 \end{lem}
 
 \noindent{\it Proof}. 
 Consider first the case $k=1$. In this case $\beta\in\,]0,1[$. Using \eqref{fundamental}, we have
 \begin{equation*}
 \int_0^\ep dr\, \langle G(r,y-\ast),G(h+r,x-\ast)\rangle_{\hac}=\int_0^\ep dr \int_{y-r}^{y+r} d\xi\int_{x-(h+r)}^{x+h+r}d\eta\, \frac{1}{|\xi-\eta|^\beta}.
 \end{equation*}
 Consider the change of variables $r\mapsto \frac{r}{\ep}$, $\xi\mapsto \frac{\xi-y}{\ep}$, $\eta\mapsto\frac{\eta-x}{\ep}$ to see that this is equal to
 \begin{equation*}
 \ep^3\int_0^1 dr \int_{-r}^{r} du\int_{-\frac{h}{\ep}-r}^{\frac{h}{\ep}+r} dv\,\frac{1}{|y-x-\ep(u-v)|^\beta},
 \end{equation*}
or even to
  \begin{equation*}
 \ep^{3-\beta}\int_0^1 dr \int_{-r}^{r} du\int_{-\frac{h}{\ep}-r}^{\frac{h}{\ep}+r} dv\,\frac{1}{|\frac{y-x}{\ep}-(u-v)|^\beta}.
 \end{equation*}
 Set
 \begin{equation*}
 \varphi(z,\lambda)= \int_0^1 dr \int_{-r}^{r} du\int_{-\lambda-r}^{\lambda+r} dv\,\frac{1}{|z+u-v|^\beta}.
 \end{equation*}
 We notice that for $\lambda \in [0,1]$,
 \begin{equation*}
  \varphi(z,\lambda)\le \int_{-1}^1 du\int_{-2}^2 dv\,
\frac{1}{|z+u-v|^\beta},
 \end{equation*}
 which clearly tends to zero as $|z|\to+\infty$.

 Next, we consider the case $k\in\{2,3\}$.
 By the definition of the inner product on $\hac$,
\begin{align*}
I&:=\int_0^\ep dr\left\langle G(r,y-\ast), G(h+r,x-\ast)\right\rangle_{\hac}\\
&\, =\int_0^\ep dr \int_{\rk}\frac{d\xi}{|\xi|^{k-\beta+2}}\, e^{i\langle\xi,x-y\rangle}\sin(r|\xi|)\sin((h+r)|\xi|).
\end{align*}
Use the trigonometric identity $\sin a\sin b=\frac{1}{2}\left[\cos(a-b)-\cos(a+b)\right]$ and integrate the obtained expression with respect to $dr$.
This gives
\begin{align*}
I&=\frac{1}{2} \int_{\rk}\frac{d\xi}{|\xi|^{k-\beta+2}}e^{i\langle\xi,x-y\rangle}\\
&\qquad\qquad\times\left[\ep\cos(h|\xi|)-\left(\frac{\sin((h+2\ep)|\xi|)}{2|\xi|}-\frac{\sin(h|\xi|)}{2|\xi|}\right)\right].
\end{align*}
Next, we make use of the identity $\sin a-\sin b=2\sin\left(\frac{a-b}{2}\right)\cos\left(\frac{a+b}{2}\right)$ and apply the change of variable
$\xi\to\ep\xi$, to write
\begin{align*}
I&=\frac{1}{2} \int_{\rk}\frac{d\xi}{|\xi|^{k-\beta+2}}e^{i\langle\xi,x-y\rangle}\\
&\qquad\qquad\times \ep\left[\cos(h|\xi|)-\frac{\sin(\ep|\xi|)\cos((h+\ep)|\xi|)}{\ep|\xi|}\right]\\
&=\frac{\ep^{3-\beta}}{2}\int_{\rk}\frac{d\xi}{|\xi|^{k-\beta+2}}e^{i\langle\xi,\frac{x-y}{\ep}\rangle}\\
&\qquad\qquad\times\left[\cos\left(\frac{h}{\ep}|\xi|\right)-\frac{\sin(|\xi|)\cos\left(\left(\frac{h}{\ep}+1\right)|\xi|\right)}{|\xi|}\right]\\
&=\frac{\ep^{3-\beta}}{2}\varphi\left(\frac{x-y}{\ep},\frac{h}{\ep}\right),
\end{align*}
where
\begin{align*}
\varphi(z,\lambda)&=\int_{\rk}\frac{d\xi}{|\xi|^{k-\beta+2}}e^{i\langle \xi,z\rangle}\\
&\qquad\qquad \times \left[\cos(\lambda|\xi|)- \frac{\sin(|\xi|)\cos((\lambda+1)|\xi|)}{|\xi|}\right].
\end{align*}
Let 
\beqn
f_0(r)=r^{\beta-k-2}\left[\cos(\lambda r)-\frac{\sin(r) \cos((\lambda+1)r)}{r}\right],
\eeqn
$r\in \mathbb{R}$, and $f(\xi)=f_0(|\xi|)$. According to \cite[p.429]{grafakos},
\beqn
\varphi(z,\lambda)=\tf f(z)=\frac{2\pi}{|z|^{\frac{k-2}{2}}}\int_0^\infty
f_0(r) J_{\frac{k}{2}-1}(2\pi r|z|)r^{\frac{k}{2}}\, dr,
\eeqn
where $J_{\frac{k}{2}-1}$ denotes the the Bessel function of order $\frac{k}{2}-1$ (see for instance \cite[p.425]{grafakos}).

Using the trigonometric formula $\cos(a+b)=\cos a\cos b-\sin a\sin b$, one can easily obtain
\beqn
\left\vert\cos(\lambda r)-\frac{\sin(r) \cos((\lambda+1)r)}{r}\right\vert\le C(r^2\wedge 1),
\eeqn
where $C$ does not depend on $\lambda\in[0,1]$.

Consequently,
\begin{align*}
 \tf f(z)&\le \frac{C}{|z|^{\frac{k-2}{2}}}\int_0^\infty r^{\beta-\frac{k}{2}-2}(r^2\wedge 1)\left\vert J_{\frac{k}{2}-1}(2\pi r|z|)\right\vert dr\\
 &=C|z|^{2-\beta}\int_0^\infty u^{\beta-\frac{k}{2}-2}\left(\frac{u^2}{|z|^2}\wedge 1\right)\left\vert J_{\frac{k}{2}-1}(2\pi u)\right\vert du,
 \end{align*}
where we have applied the change of variable $r\to r|z|$.
Set
\beqn
I(z):=|z|^{2-\beta}\int_0^\infty u^{\beta-\frac{k}{2}-2}\left(\frac{u^2}{|z|^2}\wedge 1\right)\left\vert J_{\frac{k}{2}-1}(2\pi u)\right\vert du.
\eeqn

We want to study the limit as $|z|\to\infty$ of $I(z)$. For this, we write $I(z)=I_1(z)+I_2(z)$, with
\begin{align*}
I_1(z)&= |z|^{-\beta}\int_0^{|z|}u^{\beta-\frac{k}{2}}\left\vert J_{\frac{k}{2}-1}(2\pi u)\right\vert du,\\
I_2(z)&= |z|^{2-\beta}\int_{|z|}^\infty u^{\beta-\frac{k}{2}-2}\left\vert J_{\frac{k}{2}-1}(2\pi u)\right\vert du.
\end{align*}
We consider two cases. 
\bigskip

\noindent{\em Case 1.} $k=2$ and $\beta<\frac{1}{2}$ or $k=3$ and $\beta<1$. 
For small values of $u>0$, 
$$
   \left\vert J_{\frac{k}{2}-1}(2\pi u)\right\vert\le C u^{\frac{k}{2}-1},
$$
while for
large $u$, $\left\vert J_{\frac{k}{2}-1}(2\pi u)\right\vert \leq C u^{-\frac{1}{2}}$ (see \cite{grafakos} pp. 430--433). Hence, $I_1(z) < +\infty$ for all $z \in \IR^k$, and with the indicated constraints on $k$ and $\beta$, we obtain
\beqn
\lim_{|z|\to\infty}\int_0^{|z|}u^{\beta-\frac{k}{2}}\left\vert J_{\frac{k}{2}-1}(2\pi u)\right\vert du=C<\infty.
\eeqn
Therefore,
\beq
\label{a14}
\lim_{|z|\to\infty}I_1(z)=0.
\eeq

 \noindent{\em Case 2.} $k=2$ and $\beta\ge\frac{1}{2}$ or $k=3$ and $\beta\ge1$. Then, using the expression for $J_\nu(r)$ given at the beginning of \cite[Appendix B.8]{grafakos}, we see that
\beqn
\lim_{|z|\to\infty}\int_0^{|z|}u^{\beta-\frac{k}{2}}\left\vert J_{\frac{k}{2}-1}(2\pi u)\right\vert du=\infty.
\eeqn
Therefore, \eqref{a14} follows from the Bernoulli-L'Hopital rule.

Using again the Bernoulli-L'Hopital rule, we can also check that for $k \in \{2,3\}$,
\beqn
\lim_{|z|\to\infty}I_2(z)=0.
\eeqn
The proof of the lemma is now complete.
\hfill\qed
\bigskip

\noindent{\sc Acknowledgements.} The authors are indebted to the {\it Isaac Newton Institute for Mathematical Sciences} for the hospitality and financial support during their visits there in Spring 2010, on the occasion of the scientific program {\it Stochastic Partial Differential Equations}. The first named author is pleased to thank the {\it University of Barcelona,} and the second named author thanks the {\it Ecole Polytechnique F\'ed\'erale de Lausanne} and its {\it Centre Interfacultaire Bernoulli,} for their hospitality and financial support during mutual visits during which part of this work was carried out.
\bigskip



\begin{thebibliography}{99}
\bibitem{adlermuller} R.J. Adler, P. M\"uller, B. Rozovskii (eds), {\em Stochastic modelling in physical oceanography.} Progress in Probability {\bf 39}, Birkh\"auser, Boston 1996
\bibitem{blc09} H. Bierm\'e, C. Lacaux, Y. Xiao (2009), Hitting probabilities and the Hausdorff dimension of the inverse images of anisotropic Gaussian random fields. {Bull. London. Math. Soc.} {\bf 41}, 253-273
\bibitem{biswas} S.K. Biswas, N.U. Ahmed  (1985), Stabilisation of systems governed by the wave equation in the presence of distributed white noise, {\em IEEE Trans.~on Aut.~Control,} Vol.~AC-30, 10, 1043-1045
\bibitem{cerraif} S. Cerrai, M. Freidlin  (2006), Smoluchowski-Kramers approximation for a general class of SPDEs.  {\em J. Evol. Equ.} {\bf 6}, 657-689
\bibitem{cd08} D. Conus, R.C. Dalang (2008), The non-linear stochastic wave equation in high dimensions. {\it Electron. J. Probab.} {\bf 13}, no. 22, 629-670
\bibitem{chss01} M. Chaleyat-Maurel, M. Sanz-Sol\'e (2003), Positivity of the density for the 
stochastic wave equation in two  spatial dimensions.  {\it ESAIM Probab. Stat.}
{\bf 7}  (2003), 89--114 (electronic)
\bibitem{d99} R.C. Dalang (1999), Extending the martingale measure stochastic integral with applications to spatially homegeneous s.p.d.e.'s. {\it  Electron. J. Probab.} {\bf 4}, no. 6, 1-29. With correction
\bibitem{dalmini}  R.C. Dalang (2009), The stochastic wave equation. In: {\em A minicourse on stochastic partial differential equations} (D. Khoshnevisan \& F. Rassoul-Agha, eds), Lecture Notes in Math. {\bf 1962}, Springer, Berlin, 39-71
\bibitem{df} R.C. Dalang, N.E. Frangos (1998), The stochastic wave equation in two spatial dimensions. {\it Ann. Probab.} {\bf 26}, 187-212
\bibitem{dkn07} R.C. Dalang, D. Khoshnevisan, E. Nualart (2007), Hitting probabilities for systems of non-linear stochastic heat equations
with additive noise. {\it Latin Amer. J. Probab. Statist. (ALEA)} {\bf 3}, 231-271
\bibitem{dkn09} R.C. Dalang, D. Khoshnevisan, E. Nualart (2009), Hitting probabilities for systems of non-linear stochastic heat equations
with multiplicative noise. {\it Probab. Theory and  Rel. Fields}, {\bf 144}, 371-427 
\bibitem{dkn10} R.C. Dalang, D. Khoshnevisan, E. Nualart (2013), Hitting probabilities for systems of non-linear stochastic heat equations in spatial dimensions $k\ge 1$. {\em Journal of SPDE's: Analysis and Computations} {\bf 1}-1, 94-151
\bibitem{dl1}  R.C. Dalang, O. L\'ev\^eque (2004), Second-order linear hyperbolic SPDEs driven by isotropic Gaussian noise on a sphere. {\em Ann. Probab.} {\bf 32}, 1068-1099
\bibitem{dl2}  R.C. Dalang, O. L\'ev\^eque (2006), Second-order hyperbolic S.P.D.E.'s driven by homogeneous Gaussian noise on a hyperplane. {\em Trans. Amer. Math. Soc.} {\bf 358}, 2123-2159
\bibitem{dmz06} R.C. Dalang, C. Mueller, L. Zambotti (2006), Hitting probabilities of s.p.d.e.'s with reflection. {\it Ann. Probab.} {\bf 34}, 1423-1450
\bibitem{dn04} R.C. Dalang, E. Nualart (2004), {\it Potential theory for hyperbolic SPDEs}. {Ann. Probab.} {\bf 32}, 2099-2148
\bibitem{dss09}  R.C. Dalang, M. Sanz-Sol\'e (2009), H\"older-Sobolev regularity of the solution to the stochastic wave equation in dimension three.
{\it Memoirs Amer. Math. Soc.} Vol. 199, Number 931
\bibitem{dss10} R.C. Dalang, M. Sanz-Sol\'e (2010), Criteria for hitting probabilities with applications to systems of stochastic wave equations. {\it Bernoulli} {\bf 16}(4), 1343-1368
\bibitem{dawson} D.A. Dawson, H. Salehi (1980), Spatially homogeneous random evolutions.
{\em J. Multivariate Anal.} {\bf 10}, 141-180
\bibitem{Evans} L.C. Evans, {\em Partial differential equations} (second ed). Graduate Studies in Mathematics {\bf 19}. American Math. Soc., Providence, RI, 2010
\bibitem{folland} G.B. Folland, {\it Introduction to Partial Differential Equations}. Princeton University Press, 1976
\bibitem{gonzMaddocks} O. Gonzalez, J.H. Maddocks (2001), Extracting parameters for base-pair level models of DNA from molecular dynamics simulations. {\em Theoretical Chemistry Accounts} {\bf 106}, 76-82
\bibitem{grafakos} L. Grafakos, {\it Classical Fourier Analysis}, 2nd Edition. Graduate Texts in Mathematics. Springer, New York, 2008
\bibitem{holley}  R.A. Holley, D.W. Stroock (1978), Generalized Ornstein-Uhlenbeck processes and infinite particle branching Brownian motions.{\em  Publ. Res. Inst. Math. Sci.} {\bf 14}, 741-788
\bibitem{k} D. Khoshnevisan, {\it Multiparameter Processes: An Introduction to Random Fields}. Springer, New York, 2002
\bibitem{kx09} D. Khoshnevisan, Y. Xiao (2009), Harmonic analysis of additive L\'evy processes. {\it Probab. Theory Relat. Fields} {\bf 145}, 459-515
\bibitem{mat} P. Mattila, {\it Geometry of sets and measures in Euclidean spaces}, Cambridge University Press, Cambridge 1995.
\bibitem{miller} R.N. Miller (1990), Tropical data assimilation with simulated data: The impact of the tropical ocean and global atmosphere thermal array for the ocean. {\em J. Geophysical Research} {\bf 95}, 11,461-11,482
\bibitem{mimo} A. Millet, J.-L. Morien (2000), On a stochastic wave equation in two space dimensions: regularity of the solution and its density, {\it Stoch. Proc. and their Appl.} {\bf 86} (2000), 141-162
\bibitem{mss99} A. Millet, M. Sanz-Sol\'e (1999), A stochastic wave equation in two space dimensions: smoothness of the law.  {\it Ann. Probab.} {\bf 27}, Vol. 2, 803-884
\bibitem{mt} C. Mueller, R. Tribe (2002), Hitting properties of the random string. {\it Electron. J. Probab}.  {\bf 7}, 1-29 (2002). With correction
\bibitem{n1} D. Nualart, {\it The Malliavin Calculus and Related Topics}, Springer, Heidelberg, 1995
\bibitem{n2} D. Nualart, {Analysis on Wiener space and anticipating stochastic calculus}, In: Ecole d'Et\'e de Probabilit\'es de Saint-Flour XXV. Lect. Notes in Math. Vol. 1690, pp. 123-227. Springer, Heidelberg, 1998
\bibitem{nq} D. Nualart, L. Quer-Sardanyons (2007), Existence and smoothness of the density for spatially homogeneous SPDEs, {\it Potential Analysis} {\bf 27},  281-299
\bibitem{en10} E. Nualart (2013), On the density of systems of non-linear spatially homogeneous SPDEs, {\it Stochastics}. {\bf 85}, 48-70
\bibitem{nv09} E. Nualart, F. Viens (2009), The fractional stochastic heat equation on the circle: Time regularity and potential theory, {\it  Stoch. Proc. and Their Appl.}, {\bf 119}, 1505-1540
\bibitem{pz} S. Peszat, J. Zabczyk (2000), Nonlinear stochastic wave and heat equations.  {\it Probab. Theory Related Fields}  {\bf 116} (3), 421-443 
\bibitem{qss05} L. Quer-Sardanyons, M. Sanz-Sol\'e (2004), A stochastic wave equation in dimension 3: smoothness of the law, {\it Bernoulli} {\bf 10} (1), 165-186
\bibitem{ss} M. Sanz-Sol\'e (2005),  {\it Malliavin Calculus, with Applications to Stochastic Partial Differential Equations},
EPFL Press. Fundamental Sciences, Mathematics. Distributed by CRC Press,
Taylor and Francis Group
\bibitem{ss08} M. Sanz-Sol\'e (2008), Properties of the density for a three-dimensional stochastic wave equation. {\it J. of Functional Analysis,} {\bf 255}, 255-281
\bibitem{w1} S. Watanabe, {\it Lectures on Stochastic Differential Equations and Malliavin Calculus}, Tata Institute of Fundamental Research. Bombay, 1984
\bibitem{zab} J. Zabczyk (2001),  {\it A mini course on stochastic partial differential equations}. In: {\it Stochastic climate models} (P.~Imkeller \& J.-S.~von Storch, eds), Progr. Probab. {\bf 49}, Birkh\"auser,  257-284

\end{thebibliography}
\end{document}